\documentclass[11pt,reqno]{amsart}

\usepackage{booktabs}
\usepackage{amssymb}
\usepackage{graphicx}
\usepackage{epstopdf}
\usepackage{url}
\usepackage{color}
\usepackage[usenames,dvipsnames,svgnames,table]{xcolor}
\usepackage[square,numbers,sort&compress]{natbib}

\usepackage{tikz}

\usepackage{pgfplots} 

\usetikzlibrary{positioning}

\def\visible<#1>{}  

\usepackage[utf8]{inputenc}

\usepackage[english]{babel}
\usepackage{amsfonts}
\usepackage{amsmath}
\usepackage{latexsym}
\usepackage{color}
\usepackage{subfigure}
\usepackage{enumerate}

\usepackage{hyperref}  

\usepackage{ifpdf}

\graphicspath{{figures/}}

\DeclareMathOperator    \rk                     {rk}

\newcommand{\old}[1]{{}}

\newcommand{\bb}{\mathbb}

\newcommand{\R}{\bb R}

\newcommand{\Z}{\bb Z}
\newcommand{\N}{\bb N}

\def\ve#1{\mathchoice{\mbox{\boldmath$\displaystyle\bf#1$}}
{\mbox{\boldmath$\textstyle\bf#1$}}
{\mbox{\boldmath$\scriptstyle\bf#1$}}
{\mbox{\boldmath$\scriptscriptstyle\bf#1$}}}

\newcommand{\setcond}[2]{\left\{\, #1 : #2 \,\right\}}

\renewcommand{\P}{\mathcal{P}}

\newcommand{\x}{{\ve x}}

\renewcommand{\b}{{\ve b}}

\newcommand{\B}{B}

\newenvironment{psmallmatrixbig}{\bigl(\smallmatrix}{\endsmallmatrix\bigr)}
\newcommand\InlineFrac[2]{#1/#2}  
\newcommand\ColVec[3][\relax]
{
  \ifx#1\relax
  \bgroup\let\frac=\InlineFrac\begin{psmallmatrixbig}#2\vphantom{/}\\#3\vphantom{/}\end{psmallmatrixbig}\egroup
  \else
  \bgroup\let\frac=\InlineFrac\begin{psmallmatrixbig}\ifx#200\else#2/#1\fi\\\ifx#300\else#3/#1\fi\end{psmallmatrixbig}\egroup
  \fi
}

\newtheorem{theorem}{Theorem}[section]

\makeatletter
\newcommand\MkNewTheorem[2]{%
  \newtheorem{#1}{#2}
  \expandafter\def\csname c@#1\endcsname{\c@theorem}
  \expandafter\def\csname p@#1\endcsname{\p@theorem}
  \expandafter\def\csname the#1\endcsname{\thetheorem}
  \expandafter\def\csname #1name\endcsname{#2}
}

\MkNewTheorem{corollary}{Corollary}
\MkNewTheorem{lemma}{Lemma}
\MkNewTheorem{proposition}{Proposition}
\MkNewTheorem{prop}{Proposition}
\MkNewTheorem{claim}{Claim}
\MkNewTheorem{observation}{Observation}
\MkNewTheorem{obs}{Observation}
\MkNewTheorem{conjecture}{Conjecture}
\MkNewTheorem{openquestion}{Open question}

\theoremstyle{definition}
\MkNewTheorem{example}{Example}
\MkNewTheorem{exercise}{Exercise}
\MkNewTheorem{notation}{Notation}
\MkNewTheorem{assumption}{Assumption}
\MkNewTheorem{definition}{Definition}
\MkNewTheorem{remark}{Remark}
\MkNewTheorem{goal}{Goal}

\makeatletter
\let\OurMathBbAux=\mathbb
\DeclareRobustCommand\OurMathBb{\OurMathBbAux}
\let\mathbb=\OurMathBb
\let\bfseries=\undefined
\DeclareRobustCommand\bfseries
{\not@math@alphabet\bfseries\mathbf
  \boldmath\fontseries\bfdefault\selectfont\let\OurMathBbAux=\mathbf}
\def\@thm#1#2#3{%
  \ifhmode\unskip\unskip\par\fi
  \normalfont
  \trivlist
  \let\thmheadnl\relax
  \let\thm@swap\@gobble
  \thm@notefont{\fontseries\mddefault\upshape\unboldmath}
  \thm@headpunct{.}
  \thm@headsep 5\p@ plus\p@ minus\p@\relax
  \thm@space@setup
  #1
  \@topsep \thm@preskip               
  \@topsepadd \thm@postskip           
  \def\@tempa{#2}\ifx\@empty\@tempa
    \def\@tempa{\@oparg{\@begintheorem{#3}{}}[]}%
  \else
    \refstepcounter{#2}%
    \def\@tempa{\@oparg{\@begintheorem{#3}{\csname the#2\endcsname}}[]}%
  \fi
  \@tempa
}
\makeatother

\makeatletter
\renewcommand{\pod}[1]
{\allowbreak\mathchoice{\mkern18mu}{\mkern8mu}{\mkern8mu}{\mkern8mu}(#1)}
\makeatother

\chardef\Myunderscore=`\_
\newcommand\underscore{\Myunderscore\allowbreak}

\newcommand\githubsearchurl{https://github.com/mkoeppe/infinite-group-relaxation-code/search}
\pgfkeyssetvalue{/sagefunc/california_ip}{\href{\githubsearchurl?q=\%22def+california_ip(\%22}{\sage{california\underscore{}ip}}}
\pgfkeyssetvalue{/sagefunc/automorphism}{\href{\githubsearchurl?q=\%22def+automorphism(\%22}{\sage{automorphism}}}
\pgfkeyssetvalue{/sagefunc/multiplicative_homomorphism}{\href{\githubsearchurl?q=\%22def+multiplicative_homomorphism(\%22}{\sage{multiplicative\underscore{}homomorphism}}}
\pgfkeyssetvalue{/sagefunc/projected_sequential_merge}{\href{\githubsearchurl?q=\%22def+projected_sequential_merge(\%22}{\sage{projected\underscore{}sequential\underscore{}merge}}}
\pgfkeyssetvalue{/sagefunc/restrict_to_finite_group}{\href{\githubsearchurl?q=\%22def+restrict_to_finite_group(\%22}{\sage{restrict\underscore{}to\underscore{}finite\underscore{}group}}}
\pgfkeyssetvalue{/sagefunc/restrict_to_finite_group_3}{\href{\githubsearchurl?q=\%22def+restrict_to_finite_group_3(\%22}{\sage{restrict\underscore{}to\underscore{}finite\underscore{}group\underscore{}3}}}
\pgfkeyssetvalue{/sagefunc/interpolate_to_infinite_group}{\href{\githubsearchurl?q=\%22def+interpolate_to_infinite_group(\%22}{\sage{interpolate\underscore{}to\underscore{}infinite\underscore{}group}}}
\pgfkeyssetvalue{/sagefunc/two_slope_fill_in}{\href{\githubsearchurl?q=\%22def+two_slope_fill_in(\%22}{\sage{two\underscore{}slope\underscore{}fill\underscore{}in}}}
\pgfkeyssetvalue{/sagefunc/ll_strong_fractional}{\href{\githubsearchurl?q=\%22def+ll_strong_fractional(\%22}{\sage{ll\underscore{}strong\underscore{}fractional}}}
\pgfkeyssetvalue{/sagefunc/hildebrand_2_sided_discont_1_slope_1}{\href{\githubsearchurl?q=\%22def+hildebrand_2_sided_discont_1_slope_1(\%22}{\sage{hildebrand\underscore{}2\underscore{}sided\underscore{}discont\underscore{}1\underscore{}slope\underscore{}1}}}
\pgfkeyssetvalue{/sagefunc/hildebrand_2_sided_discont_2_slope_1}{\href{\githubsearchurl?q=\%22def+hildebrand_2_sided_discont_2_slope_1(\%22}{\sage{hildebrand\underscore{}2\underscore{}sided\underscore{}discont\underscore{}2\underscore{}slope\underscore{}1}}}
\pgfkeyssetvalue{/sagefunc/hildebrand_discont_3_slope_1}{\href{\githubsearchurl?q=\%22def+hildebrand_discont_3_slope_1(\%22}{\sage{hildebrand\underscore{}discont\underscore{}3\underscore{}slope\underscore{}1}}}
\pgfkeyssetvalue{/sagefunc/dr_projected_sequential_merge_3_slope}{\href{\githubsearchurl?q=\%22def+dr_projected_sequential_merge_3_slope(\%22}{\sage{dr\underscore{}projected\underscore{}sequential\underscore{}merge\underscore{}3\underscore{}slope}}}
\pgfkeyssetvalue{/sagefunc/chen_4_slope}{\href{\githubsearchurl?q=\%22def+chen_4_slope(\%22}{\sage{chen\underscore{}4\underscore{}slope}}}
\pgfkeyssetvalue{/sagefunc/gmic}{\href{\githubsearchurl?q=\%22def+gmic(\%22}{\sage{gmic}}}
\pgfkeyssetvalue{/sagefunc/gj_2_slope}{\href{\githubsearchurl?q=\%22def+gj_2_slope(\%22}{\sage{gj\underscore{}2\underscore{}slope}}}
\pgfkeyssetvalue{/sagefunc/gj_2_slope_repeat}{\href{\githubsearchurl?q=\%22def+gj_2_slope_repeat(\%22}{\sage{gj\underscore{}2\underscore{}slope\underscore{}repeat}}}
\pgfkeyssetvalue{/sagefunc/dg_2_step_mir}{\href{\githubsearchurl?q=\%22def+dg_2_step_mir(\%22}{\sage{dg\underscore{}2\underscore{}step\underscore{}mir}}}
\pgfkeyssetvalue{/sagefunc/kf_n_step_mir}{\href{\githubsearchurl?q=\%22def+kf_n_step_mir(\%22}{\sage{kf\underscore{}n\underscore{}step\underscore{}mir}}}
\pgfkeyssetvalue{/sagefunc/gj_forward_3_slope}{\href{\githubsearchurl?q=\%22def+gj_forward_3_slope(\%22}{\sage{gj\underscore{}forward\underscore{}3\underscore{}slope}}}
\pgfkeyssetvalue{/sagefunc/drlm_backward_3_slope}{\href{\githubsearchurl?q=\%22def+drlm_backward_3_slope(\%22}{\sage{drlm\underscore{}backward\underscore{}3\underscore{}slope}}}
\pgfkeyssetvalue{/sagefunc/drlm_2_slope_limit}{\href{\githubsearchurl?q=\%22def+drlm_2_slope_limit(\%22}{\sage{drlm\underscore{}2\underscore{}slope\underscore{}limit}}}
\pgfkeyssetvalue{/sagefunc/drlm_2_slope_limit_1_1}{\href{\githubsearchurl?q=\%22def+drlm_2_slope_limit_1_1(\%22}{\sage{drlm\underscore{}2\underscore{}slope\underscore{}limit\underscore{}1\underscore{}1}}}
\pgfkeyssetvalue{/sagefunc/bhk_irrational}{\href{\githubsearchurl?q=\%22def+bhk_irrational(\%22}{\sage{bhk\underscore{}irrational}}}
\pgfkeyssetvalue{/sagefunc/bccz_counterexample}{\href{\githubsearchurl?q=\%22def+bccz_counterexample(\%22}{\sage{bccz\underscore{}counterexample}}}
\pgfkeyssetvalue{/sagefunc/drlm_3_slope_limit}{\href{\githubsearchurl?q=\%22def+drlm_3_slope_limit(\%22}{\sage{drlm\underscore{}3\underscore{}slope\underscore{}limit}}}
\pgfkeyssetvalue{/sagefunc/dg_2_step_mir_limit}{\href{\githubsearchurl?q=\%22def+dg_2_step_mir_limit(\%22}{\sage{dg\underscore{}2\underscore{}step\underscore{}mir\underscore{}limit}}}
\pgfkeyssetvalue{/sagefunc/rlm_dpl1_extreme_3a}{\href{\githubsearchurl?q=\%22def+rlm_dpl1_extreme_3a(\%22}{\sage{rlm\underscore{}dpl1\underscore{}extreme\underscore{}3a}}}
\pgfkeyssetvalue{/sagefunc/hildebrand_5_slope_22_1}{\href{\githubsearchurl?q=\%22def+hildebrand_5_slope_22_1(\%22}{\sage{hildebrand\underscore{}5\underscore{}slope\underscore{}22\underscore{}1}}}
\pgfkeyssetvalue{/sagefunc/hildebrand_5_slope_24_1}{\href{\githubsearchurl?q=\%22def+hildebrand_5_slope_24_1(\%22}{\sage{hildebrand\underscore{}5\underscore{}slope\underscore{}24\underscore{}1}}}
\pgfkeyssetvalue{/sagefunc/hildebrand_5_slope_28_1}{\href{\githubsearchurl?q=\%22def+hildebrand_5_slope_28_1(\%22}{\sage{hildebrand\underscore{}5\underscore{}slope\underscore{}28\underscore{}1}}}
\pgfkeyssetvalue{/sagefunc/kzh_7_slope_1}{\href{\githubsearchurl?q=\%22def+kzh_7_slope_1(\%22}{\sage{kzh\underscore{}7\underscore{}slope\underscore{}1}}}
\pgfkeyssetvalue{/sagefunc/kzh_28_slope_1}{\href{\githubsearchurl?q=\%22def+kzh_28_slope_1(\%22}{\sage{kzh\underscore{}28\underscore{}slope\underscore{}1}}}
\pgfkeyssetvalue{/sagefunc/extremality_test}{\href{\githubsearchurl?q=\%22def+extremality_test(\%22}{\sage{extremality\underscore{}test}}}
\pgfkeyssetvalue{/sagefunc/plot_2d_diagram}{\href{\githubsearchurl?q=\%22def+plot_2d_diagram(\%22}{\sage{plot\underscore{}2d\underscore{}diagram}}}
\pgfkeyssetvalue{/sagefunc/generate_example_e_for_psi_n}{\href{\githubsearchurl?q=\%22def+generate_example_e_for_psi_n(\%22}{\sage{generate\underscore{}example\underscore{}e\underscore{}for\underscore{}psi\underscore{}n}}}
\pgfkeyssetvalue{/sagefunc/chen_3_slope_not_extreme}{\href{\githubsearchurl?q=\%22def+chen_3_slope_not_extreme(\%22}{\sage{chen\underscore{}3\underscore{}slope\underscore{}not\underscore{}extreme}}}
\pgfkeyssetvalue{/sagefunc/psi_n_in_bccz_counterexample_construction}{\href{\githubsearchurl?q=\%22def+psi_n_in_bccz_counterexample_construction(\%22}{\sage{psi\underscore{}n\underscore{}in\underscore{}bccz\underscore{}counterexample\underscore{}construction}}}
\pgfkeyssetvalue{/sagefunc/gomory_fractional}{\href{\githubsearchurl?q=\%22def+gomory_fractional(\%22}{\sage{gomory\underscore{}fractional}}}
\pgfkeyssetvalue{/sagefunc/not_extreme_1}{\href{\githubsearchurl?q=\%22def+not_extreme_1(\%22}{\sage{not\underscore{}extreme\underscore{}1}}}
\pgfkeyssetvalue{/sagefunc/bhk_irrational_extreme_limit_to_rational_nonextreme}{\href{\githubsearchurl?q=\%22def+bhk_irrational_extreme_limit_to_rational_nonextreme(\%22}{\sage{bhk\underscore{}irrational\underscore{}extreme\underscore{}limit\underscore{}to\underscore{}rational\underscore{}nonextreme}}}
\pgfkeyssetvalue{/sagefunc/zhou_two_sided_discontinuous_cannot_assume_any_continuity}{\href{\githubsearchurl?q=\%22def+zhou_two_sided_discontinuous_cannot_assume_any_continuity(\%22}{\sage{zhou\underscore{}two\underscore{}sided\underscore{}discontinuous\underscore{}cannot\underscore{}assume\underscore{}any\underscore{}continuity}}}
\pgfkeyssetvalue{/sagefunc/drlm_gj_2_slope_extreme_limit_to_nonextreme}{\href{\githubsearchurl?q=\%22def+drlm_gj_2_slope_extreme_limit_to_nonextreme(\%22}{\sage{drlm\underscore{}gj\underscore{}2\underscore{}slope\underscore{}extreme\underscore{}limit\underscore{}to\underscore{}nonextreme}}}
\pgfkeyssetvalue{/sagefunc/not_minimal_2}{\href{\githubsearchurl?q=\%22def+not_minimal_2(\%22}{\sage{not\underscore{}minimal\underscore{}2}}}
\pgfkeyssetvalue{/sagefunc/drlm_not_extreme_1}{\href{\githubsearchurl?q=\%22def+drlm_not_extreme_1(\%22}{\sage{drlm\underscore{}not\underscore{}extreme\underscore{}1}}}
\pgfkeyssetvalue{/sagefunc/kzh_2q_example_1}{\href{\githubsearchurl?q=\%22def+kzh_2q_example_1(\%22}{\sage{kzh\underscore{}2q\underscore{}example\underscore{}1}}}

\DeclareRobustCommand\sage[1]{\texttt{#1}}
\DeclareRobustCommand\sagefunc[1]{\pgfkeys{/sagefunc/#1}}

\usepackage[normalem]{ulem}

 \newcommand{\tgreen}[1]{\textsf{\textcolor {ForestGreen} {#1}}}
 
 \newcommand{\tblue}[1]{\textcolor {blue} {#1}}

\pgfkeyssetvalue{/sagefunc/gj_2_slope}{\href{\githubsearchurl?q=\%22def+gj_2_slope(\%22}{\sage{gj\underscore{}2\underscore{}slope}}}
\pgfkeyssetvalue{/sagefunc/gj_forward_3_slope}{\href{\githubsearchurl?q=\%22def+gj_forward_3_slope(\%22}{\sage{gj\underscore{}forward\underscore{}3\underscore{}slope}}}
\pgfkeyssetvalue{/sagefunc/kzh_28_slope_1}{\href{\githubsearchurl?q=\%22def+kzh_28_slope_1(\%22}{\sage{kzh\underscore{}28\underscore{}slope\underscore{}1}}}
\pgfkeyssetvalue{/sagefunc/kzh_2q_example_1}{\href{\githubsearchurl?q=\%22def+kzh_2q_example_1(\%22}{\sage{kzh\underscore{}2q\underscore{}example\underscore{}1}}}
\pgfkeyssetvalue{/sagefunc/kzh_5_slope_fulldim_1}{\href{\githubsearchurl?q=\%22def+kzh_5_slope_fulldim_1(\%22}{\sage{kzh\underscore{}5\underscore{}slope\underscore{}fulldim\underscore{}1}}}
\pgfkeyssetvalue{/sagefunc/chen_4_slope}{\href{\githubsearchurl?q=\%22def+chen_4_slope(\%22}{\sage{chen\underscore{}4\underscore{}slope}}}
\pgfkeyssetvalue{/sagefunc/gmic}{\href{\githubsearchurl?q=\%22def+gmic(\%22}{\sage{gmic}}}
\pgfkeyssetvalue{/sagefunc/kzh_5_slope_fulldim_covers_1}{\href{\githubsearchurl?q=\%22def+kzh_5_slope_fulldim_covers_1(\%22}{\sage{kzh\underscore{}5\underscore{}slope\underscore{}fulldim\underscore{}covers\underscore{}1}}}
\pgfkeyssetvalue{/sagefunc/kzh_6_slope_fulldim_covers_1}{\href{\githubsearchurl?q=\%22def+kzh_6_slope_fulldim_covers_1(\%22}{\sage{kzh\underscore{}6\underscore{}slope\underscore{}fulldim\underscore{}covers\underscore{}1}}}
\pgfkeyssetvalue{/sagefunc/kzh_6_slope_1}{\href{\githubsearchurl?q=\%22def+kzh_6_slope_1(\%22}{\sage{kzh\underscore{}6\underscore{}slope\underscore{}1}}}
\pgfkeyssetvalue{/sagefunc/kzh_7_slope_4}{\href{\githubsearchurl?q=\%22def+kzh_7_slope_4(\%22}{\sage{kzh\underscore{}7\underscore{}slope\underscore{}4}}}
\pgfkeyssetvalue{/sagefunc/kzh_10_slope_1}{\href{\githubsearchurl?q=\%22def+kzh_10_slope_1(\%22}{\sage{kzh\underscore{}10\underscore{}slope\underscore{}1}}}
\pgfkeyssetvalue{/sagefunc/not_extreme_1}{\href{\githubsearchurl?q=\%22def+not_extreme_1(\%22}{\sage{not\underscore{}extreme\underscore{}1}}}
\pgfkeyssetvalue{/sagefunc/arithmetic_complexity}{\href{\githubsearchurl?q=\%22def+arithmetic_complexity(\%22}{\sage{arithmetic\underscore{}complexity}}}
\pgfkeyssetvalue{/sagefunc/number_of_slopes}{\href{\githubsearchurl?q=\%22def+number_of_slopes(\%22}{\sage{number\underscore{}of\underscore{}slopes}}}
\pgfkeyssetvalue{/sagefunc/number_of_components}{\href{\githubsearchurl?q=\%22def+number_of_components(\%22}{\sage{number\underscore{}of\underscore{}components}}}

\usepackage[norelsize,ruled]{algorithm2e}   
\usepackage{booktabs,array,ragged2e}
\usepackage{pdflscape}
\usepackage{blkarray}
\usepackage{etoolbox}
\usepackage{tikz}
\usetikzlibrary{calc}
\usetikzlibrary{spy, backgrounds, shadows}
\usetikzlibrary{shapes.symbols}
\usetikzlibrary{decorations.pathmorphing}

\newcommand{\matindex}[1]{\mbox{\scriptsize#1}}
\definecolor{mediumspringgreen}{rgb}{0.0, 0.98039215, 0.60392156}
\definecolor{darkgreen}{rgb}{0.0, 0.39215686274509803, 0.0}

\newcommand\CPL{\mathrm{CPL}}
\newcommand\DPL{\mathrm{DPL}}

\AtBeginDocument{ 

}

\renewcommand\tblue[1]{}
\renewcommand\tgreen[1]{}

\title[Computer-based search for extreme functions]{New computer-based search strategies for extreme functions of the
  Gomory--Johnson infinite group problem}

\thanks{The authors acknowledge partial support from the National Science
  Foundation through grant DMS-1320051 awarded to M.~K\"oppe.}

\author{Matthias K\"oppe}
\address{Matthias K\"oppe: Dept.\ of Mathematics, University of California, Davis}
\email{mkoeppe@math.ucdavis.edu}

\author{Yuan Zhou} 
\address{Yuan Zhou: Dept.\ of Mathematics, University of California, Davis}
\email{yzh@math.ucdavis.edu}

\date{$\relax$Revision: 2130 $ - \ $Date: 2016-10-01 14:32:34 -0700 (Sat, 01 Oct 2016) $ $\!\!\!}

\begin{document}

\begin{abstract}
  We describe new computer-based search strategies for extreme functions for
  the Gomory--Johnson infinite group problem.  They lead to the discovery of
  new extreme functions, whose existence settles several open questions.
\end{abstract}

\maketitle

\begin{tikzpicture} 
[spy scope={magnification=5, size=3cm},
   every spy in node/.style={
     magnifying glass, circular drop shadow,
     fill=white, draw=red, ultra thick, cap=round}]
    \node [anchor=south east] (image1) at (0,0) {%
      \pgfimage[height=0.58\textheight]{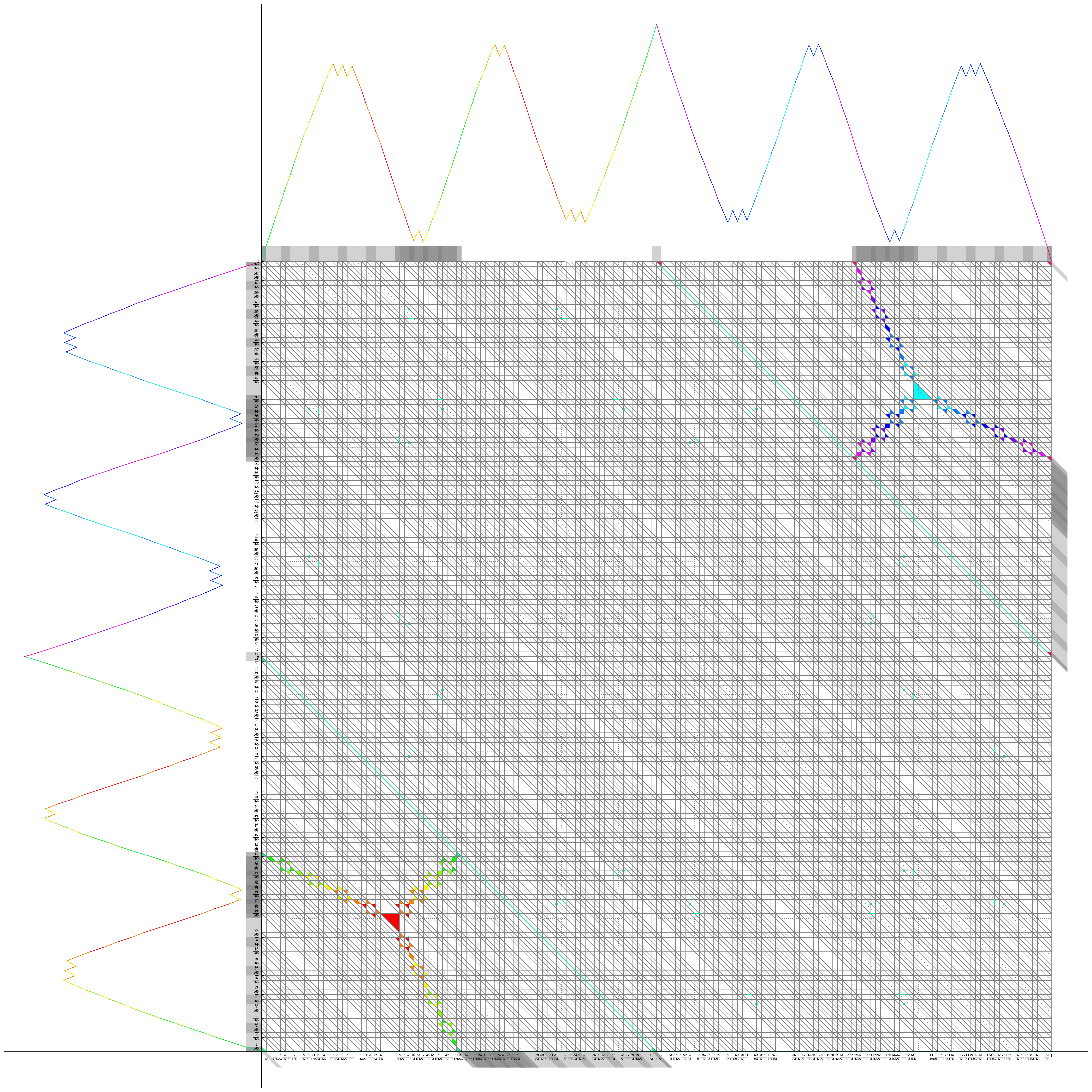}%
    };
  \spy on (-2.3,7.5) in node;
\end{tikzpicture}

\clearpage
\enlargethispage{1ex}
{\footnotesize
\tableofcontents}
\clearpage

\section{Introduction and statement of results}
\label{sec:introduction}
\subsection{Group relaxations and extreme functions}
Cutting planes are widely used in the state-of-the-art integer programming solvers. Important sources of general-purpose cutting planes are the master finite group relaxation of an integer program, which was introduced by Gomory in 1969 \cite{gom}, and the infinite group relaxation by Gomory and Johnson \cite{infinite, infinite2}. Due to the pressing need for effective cutting planes, the group problem has received renewed attention in the recent years, since it may be the key to new multi-row cutting plane approaches that have better performance than the ones in use today.

Computer-based search has been used for the investigations of Gomory's group problem and
Gomory--Johnson's infinite group problem since the very beginning, leading to the discovery of many cutting planes.
In this paper, we develop new computer-based search strategies to carry forward the discovery.

We restrict ourselves to the single-row (or, one-dimensional) problem. 
That is, we focus on only one row of a simplex tableau of an integer program. Suppose the row corresponding to some basic variable $x$ is of the form
\begin{equation}
\label{eq:single-row}
 \begin{aligned}
 x =& -f + \sum_{j = 1}^m r_j y_j,\\
 & x \in \Z_+,\\
 & y_j \in \Z_+, \forall j \in \{1,2,\dots,m\},
 \end{aligned}
\end{equation}
where $\{y_j\}_{j=1}^m$ denote the nonbasic variables.
We assume $f \in \R \setminus \Z$, that is, the basic variable $x$ is currently fractional.

When all data is rational, there exists some integer $q > 0$ such that $r_j \in \frac{1}{q} \Z$ for any $j \in \{1,2,\dots,m\}$ and $f \in \frac{1}{q}\Z$. 
\emph{Gomory's master finite (cyclic) group problem} of order $q$ is obtained by relaxing the basic variable $x \in \Z_+$ to $x \in \Z$ and by introducing variables $y(r) \in \Z_+$ for every $r \in \frac{1}{q}\Z$.
Using the quotient group $G/\Z$, i.e., reducing modulo $1$, and standard aggregation of variables whose coefficients are the same modulo $1$ (see \cite[Remark 2.1]{igp_survey}), the relaxation of \eqref{eq:single-row} takes the form
\begin{equation}
 \label{GP-finite}
  \begin{aligned}
    &\sum_{r \in G/\Z} r\, y(r) = f + \Z, \\
    & y(r) \in \Z_+, \forall r \in G/\Z,
  \end{aligned}
\end{equation}
where $G =  \frac{1}{q}\Z$ and $f$ is a given element of $G \setminus \Z$. 
The master finite group problem only depends on the parameters $f$ and  $q$, but not on any other problem data.

\emph{Gomory--Johnson's infinite group problem} is obtained by further introducing infinitely many new variables $y(r) \in \Z_+$ for every $r \in \R$. Formally, again by aggregation of variables, it can be written as
\begin{equation}
  \label{GP-infinite} 
  \begin{aligned}
    &\sum_{r \in G/\Z} r\, y(r) = f+\Z, \\
    &y\colon G /\Z \to\Z_+ \text{ is a function of finite support}, 
  \end{aligned}
\end{equation}
where $G=\R$ and $f$ is a given element of $G \setminus \Z$. 
The infinite group problem only depends on the parameter $f$.

We study the convex hull $R_{f}(G/\Z)$ of the set of functions $y\colon
G/\Z \to \Z_+$ satisfying the constraints in~\eqref{GP-finite} and
in~\eqref{GP-infinite} for the finite and infinite group problems
respectively.\footnote{For simplicity of notation in the $n$-dimensional
  ($n$-row) problem, \cite{igp_survey} and \cite{igp_survey_part_2} work with $R_f(\R^n, \Z^n)$ instead of
  the aggregated formulation $R_f(\R^n/\Z^n)$.  The aggregated formulation is
  of interest for the present paper, since for a master finite group problem where $G=\frac{1}{q}\Z$, the group $G/\Z$ is indeed finite.}
The elements of the convex hull are understood as functions $y\colon G/\Z \to \R_+$. 

A function $\pi\colon G/\Z \to\R$ is called a \emph{valid function} for $R_{f}(G/\Z)$ if 
\begin{equation}
\label{eq:valid-function}
\sum_{r \in G/\Z} \pi(r)y(r) \geq 1
\end{equation}
holds for any $y \in R_{f}(G/\Z)$.
\emph{Minimal (valid) functions} are those valid functions that are pointwise minimal. Let $\Pi_f(G/\Z)$ denote the set of minimal functions for $R_f(G/\Z)$. \emph{Extreme functions} are those valid functions that are not a proper convex combination of other valid functions. We focus on the extreme functions because they serve as strong cut-generating functions for general integer linear programs. 
(A related notion, \emph{facets}, has been studied in parts of the
literature.  For the finite group problem and for special cases of the
infinite group problem, it is known to be equivalent to that of extreme
functions; see \cite[section~2.2.4]{igp_survey} and the forthcoming
paper~\cite{koeppe-zhou:discontinuous-facets}.)

\subsection{Extreme functions and their slopes}
In this paper, we discuss how computer-based search can help in finding extreme functions.
The next two subsections (\ref{sec:intro_finite_literature}
 and \ref{sec:intro_infinite_literature}) are devoted to a short literature review on the
success of computer-based search in these problems.

\begin{figure}[t]
\centering
\includegraphics[width=.4\linewidth]{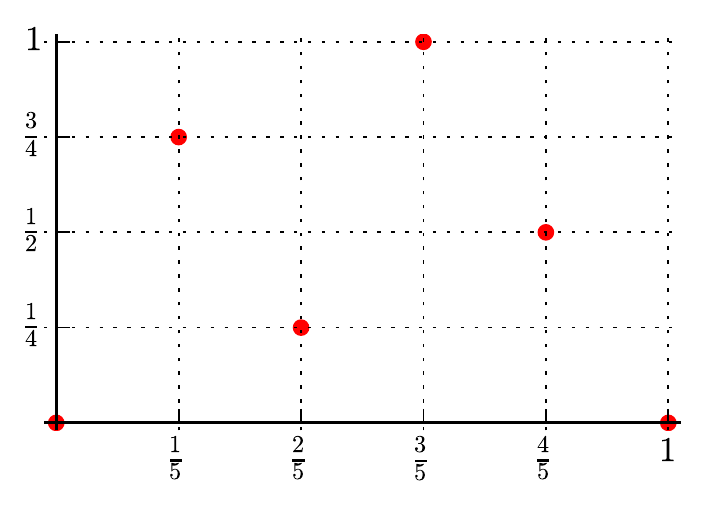}
\quad
\includegraphics[width=.4\linewidth]{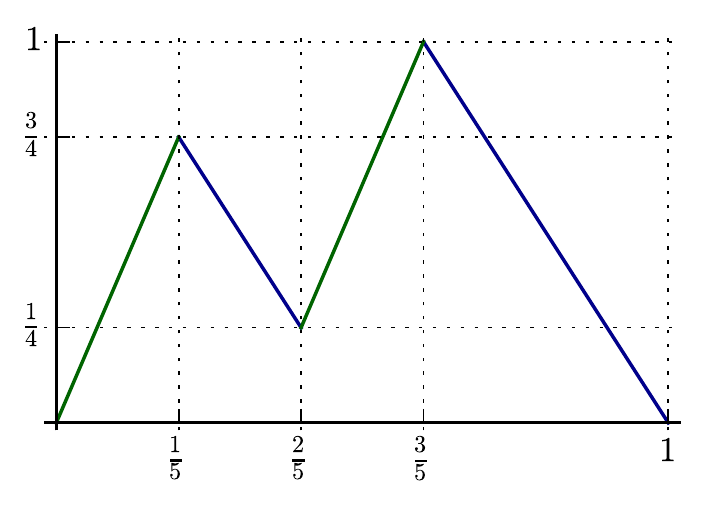}
\caption{The 2-slope extreme function
  \sagefunc{gj_2_slope}\protect\footnotemark, discovered by Gomory and Johnson
  \cite{tspace}. \textit{Left}, \sagefunc{gj_2_slope} for the finite group
  problem with $q=5$ and $f=\frac{3}{5}$, obtained by \sage{\sagefunc{restrict_to_finite_group}(\sagefunc{gj_2_slope}())}. It is a discrete function whose interpolation is the right subfigure. \textit{Right}, \sagefunc{gj_2_slope} for the infinite group problem with $f=\frac{3}{5}$. It is a continuous piecewise linear function with two slopes, although it has four pieces. Its restriction to $\frac{1}{5}\Z$ is the left subfigure.}
\label{fig:gj2s_restriction_interpolation}
\end{figure}

An important statistic that has received much attention in the literature is the number of slopes of an extreme function. For the infinite group problem, we use the term \emph{$k$-slope function} to refer to a continuous piecewise linear function with $k$ different slope values, whereas for the finite group problem, we use the same term to refer to a discrete function whose interpolation has $k$ different slope values.
\autoref{fig:gj2s_restriction_interpolation} shows a 2-slope function for the finite (left) and infinite (right) group problems respectively.
\footnotetext{Throughout this paper, we refer to an extreme function or a
  family of extreme functions by the name of the SageMath function in the
  Electronic Compendium \cite{electronic-compendium,zhou:extreme-notes} that constructs them;
  these names are shown in typewriter font.
  The reader is invited to explore these functions alongside reading this article.
  The Electronic Compendium is part of our software
  \cite{infinite-group-relaxation-code}, which allows to test extremality 
  and which has been used to make most of the diagrams in this paper. 
  The captions of some figures show SageMath code that can be used to reproduce
  the diagrams.
}
The number of slopes can be taken as a measure of complexity of an extreme
function.  Functions with only 2 slopes are the easiest to analyze; 
in fact, the Gomory--Johnson Two Slope Theorem \cite{infinite} states that any
continuous 2-slope function that is minimal valid for $R_f(\R/\Z)$ is already
extreme.  Gomory--Johnson in their 2003 paper \cite[section 6.2]{tspace},
after discussing the tools available for analyzing minimal functions for the
infinite group problem, explain:
\begin{quote}
  [T]heir application to generate facets is still somewhat ad~hoc.
  Also we don't have general theorems for three or more slope facets analogous to the
  Gomory--Johnson Two Slope Theorem.
\end{quote}
Indeed, functions with 2 or 3 slopes were the focus of much of the
literature, as we will discuss in more detail below. 

\subsection{Computer-based search used in the finite group problem}
\label{sec:intro_finite_literature}
Gomory's seminal paper \cite[Appendix 5]{gom}, introducing the group problem and corner polyhedra, listed all extreme functions up to automorphism and homomorphism for the finite group problems of order $q = 2,3, \dots, 11$.
Gomory proved that the set $\Pi_f(\frac{1}{q}\Z/\Z)$ of minimal functions for $R_f(\frac{1}{q}\Z/\Z)$ is a polytope, defined by linear inequalities that express subadditivity and certain equations that come from normalization; see \autoref{thm:finite-minimal} below for details. By~\cite[Theorem 18]{gom}\footnote{In~\cite{gom} and~\cite{evans-thesis} below, valid inequalities are not normalized to have the right hand side of \eqref{eq:valid-function} being $1$. We state their results in our unified notation.} and \cite[Theorem 2.2]{infinite}, the extreme functions for $R_f(\frac{1}{q}\Z/\Z)$ are the extreme points of the polytope $\Pi_f(\frac{1}{q}\Z/\Z)$. Gomory reported that the extreme points were computed, by enumerating simplex bases, using a computer code of Balinski and Wolfe.

During the revival of the interest in the group problem in the 2000s, Evans
\cite{evans-thesis} used her specialized implementation\footnote{\cite[Chapter
  4]{evans-thesis} used a variation of the double description method that
  includes a parallel implementation for maintaining the minimal system of
  generators. Evans reports that the parallel version achieved a speedup by a
  factor of $12.79$ using $32$ processors.} of the double description method
(see, e.g., \cite{fukuda1996dd}) to enumerate all extreme points of the
polytope $\Pi_f(\frac{1}{q}\Z/\Z)$, thereby obtaining all the extreme
functions for the finite group problems of order $q \leq 24$. 
By exploring the patterns of such functions, some parametric families of 2-slope and 3-slope extreme functions for finite group problems were constructed. Extreme functions from these families were generated by the Matlab code in \cite[Appendix B.1]{evans-thesis} for the finite group problems of order $q \leq 30$. Evans reported that these extreme functions received a large percentage of hits in the so-called shooting experiment \cite[Table 13]{evans-thesis}.

Gomory and Johnson \cite{infinite} showed that the number of extreme functions
grows exponentially with $q$. Hence it is impractical to enumerate all
extreme functions for $R_f(\frac{1}{q}\Z,\Z)$ when $q$ is large. The shooting
experiment was conducted in \cite{gomory2003corner} (more results appeared in
\cite{evans-thesis}) to identify the ``important'' extreme functions for the
finite group problems where the order is at most~$30$. This experiment was
extended to finite group problems of order up to~$90$, and then to problems of
order up to~$200$ with the so-called ``partial shooting'' variant in
\cite{dash06:shooting}; see also~\cite{shim2009large}. Extreme functions resulting from the shooting
experiments were expected to be important computationally in branch-and-cut
(see, e.g., \cite[Section 19.6.2]{Richard-Dey-2010:50-year-survey} for a
summary), though actual computational uses never seem to have
materialized. They are mostly \sagefunc{gmic} functions (up to multiplicative homomorphism and automorphism\footnote{These two procedures are available as \sagefunc{multiplicative_homomorphism} and \sagefunc{automorphism}, respectively, in the accompanying SageMath program.}), along
with some other 2-slope and 3-slope extreme functions. The shooting
experiment, however, is not suitable for finding functions with many slopes, 
as those appear to be extremely rare from the viewpoint of shooting, 
or functions with specific properties for finite group problems. It is not possible to perform the experiment for the infinite group problem either, according to \cite[Section 19.6.2.1]{Richard-Dey-2010:50-year-survey}.

Ar{\'a}oz, Evans, Gomory and Johnson \cite{AraozEvansGomoryJohnson03} demonstrated a close relation between the master finite group problem and the master knapsack problem. In particular, the convex hull $P(K_{q-1})$ of solutions to the master knapsack problem of size $q-1$ is a facet of the master finite group problem $R_f(\frac{1}{q}\Z/\Z)$, where $f = \frac{q-1}{q}$. Thus, extreme functions for $R_f(\frac{1}{q}\Z/\Z)$ where $f = \frac{q-1}{q}$ are all valid for the knapsack problem $K_{q-1}$. Furthermore, some extreme functions for $R_f(\frac{1}{q}\Z/\Z)$ 
may be obtained from facets for $P(K_{qf})$ through a process called \textit{tilting} (see \cite[Theorem 5.2]{AraozEvansGomoryJohnson03}). Examples of extreme functions derived by tilting a knapsack facet are listed in \cite[Table B.2]{AraozEvansGomoryJohnson03}.

A different approach was followed by Richard, Li and Miller \cite{Richard-Li-Miller-2009:Approximate-Liftings}, who proposed an approximate lifting scheme that converts certain superadditive functions into potentially strong valid inequalities. The superadditive functions that they studied were the $\DPL_n$ functions with $4n$ non-negative parameters, and the superadditive $\CPL_n$ functions as a special case of $\DPL_n$ where $2n$ parameters were fixed to $0$.
The parameters that define a superadditive $\DPL_n$ function belong to a certain polyhedron $P\Theta_n$ (or, in the case of  superadditive $\CPL_n$ function, to a simpler polyhedron that is a face of $P\Theta_n$). Several classes of well-known cutting planes can be generated by converting the $\DPL_n$ or $\CPL_n$ functions that correspond to the extreme points of $P\Theta_n$. However, the functions generated by the approximate lifting scheme are not always extreme.
By the lack of any automated extremality tests for a parametric family of functions, the study was restricted to so small $n$ that manual inspection of extremality became possible. 
The authors investigated the $\CPL_2$ functions for the finite group problem in \cite{Miller-Li-Richard2008} and a special case of $\CPL_3$ functions for both finite and infinite group problems in \cite{Richard-Li-Miller-2009:Approximate-Liftings}, all of which required extensive case analysis for extremality conditions by hand. They found the first parametric family of 4-slope extreme functions for the finite group problem in \cite{Richard-Li-Miller-2009:Approximate-Liftings}.

\subsection{Computer-based search used in the infinite group problem}
\label{sec:intro_infinite_literature}
Computer-based search has also been used in the study of the infinite group
problem.  One approach focuses on continuous piecewise linear functions $\pi$ 
with breakpoints in $\tfrac{1}{q} \Z$ for some fixed $q \in \Z_+$.  (Then, without
loss of generality, one can assume that also $f \in \tfrac{1}{q} \Z$ \cite[Lemma
2.4]{basu-hildebrand-koeppe:equivariant}.)  
Gomory and Johnson \cite{infinite} established a connection between finite and
infinite group problems by studying the restriction and interpolation of valid
functions.
They proved that a continuous piecewise linear function $\pi$ with breakpoints
and $f$ in $\frac{1}{q}\Z$ is minimal 
if and only if its
restriction to $\frac1q\Z$ 
is minimal for the finite group problem
$R_f(\tfrac{1}{q} \Z/ \Z)$.   
Moreover, if $\pi$ is extreme for $R_f(\R/\Z)$, then the restriction $\pi|_{\frac{1}{q} \Z}$
is extreme for $R_f(\tfrac{1}{q} \Z/\Z)$. 
Therefore all search approaches for the finite group problem, described above
in~\autoref{sec:intro_finite_literature},  yield \emph{candidates} for extreme
functions for the infinite group problem.

\begin{figure}[t]
\centering
\includegraphics[width=.6\linewidth]{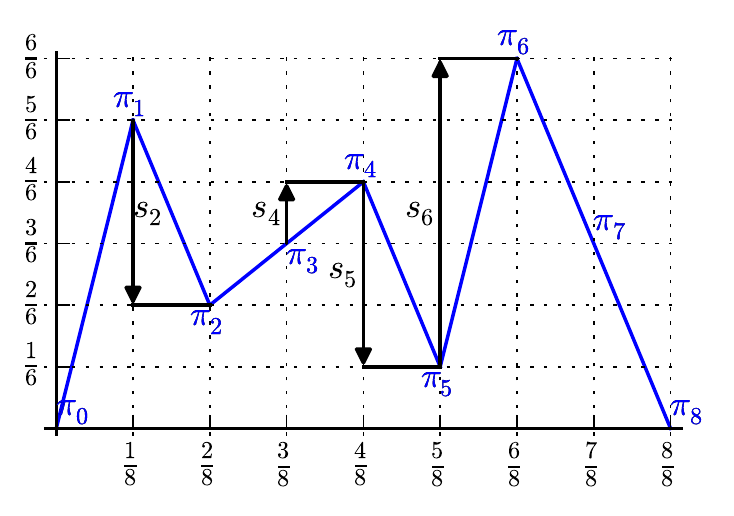}
\caption{The $q \times v$ grid discretization of the space of   continuous piecewise linear functions with rational data. Here $q=8$ and $v=6$.
}
\label{fig:q_v_grid}
\end{figure}

Using this connection with the finite group problem $R_f(\tfrac{1}{q} \Z/ \Z)$
amounts to discretizing the space of functions $\pi$, by discretizing the
breakpoints of~$\pi$ in $\frac{1}{q}\Z$ for some fixed~$q$.   
A function $\pi$ is then uniquely determined by its values at $\frac{i}{q}$ for
$i=0,1,\dots,q$, or by its slope values on the intervals
$[\frac{i-1}{q},\frac{i}{q}]$ for $i=1,\dots,q$. 
We introduce the following notation for use throughout the paper.
Denote the function value at $\frac{i}{q}$ by $\pi_i$ for $i \in \{0, 1, \dots, q\}$, where $\pi_0 = \pi_q = 0$.
Let $s_i$ be real numbers such that $q s_i$ are the slope values  on $[\frac{i-1}{q},\frac{i}{q}]$ for $i \in \{1, \dots, q \}$.
See \autoref{fig:q_v_grid} for an illustration.

Chen \cite{chen} designed an enumerative algorithm to find candidate
piecewise linear extreme functions for the infinite group problem using an
additional discretization.  In addition to fixing~$q$, also pick a natural
number~$v$. Then consider the continuous piecewise linear functions $\pi$ that
have \emph{$q \times v$ grid discretization}: the breakpoints being in $\{0, \frac{1}{q},
\dots, \frac{q-1}{q}, 1\}$ and  
\begin{enumerate}
\item $\pi_i \in \{0, \frac{1}{v}, \dots, \frac{v-1}{v}, 1\}$  for $i \in \{0, \dots, q \}$, or
\item $s_i \in \{-1, \dots, \frac{-1}{v}, 0, \frac{1}{v}, \dots, \frac{v-1}{v}, 1\}$  for $i \in \{1, \dots, q \}$.
\end{enumerate} 
In fact, these two natural ways of discretization are easily seen to be
equivalent (\autoref{lem:values-slopes-same-denominator} in \autoref{sec:limitation-of-grid}).
Chen's algorithm then enumerated every candidate function $\pi$ such that $\pi$ is symmetric, $\pi(0)=0$, $\pi(f)=1$, $\pi(1)=0$, and $\pi$ has the steepest positive and negative slopes at $0$ and $1$ respectively. With $q =10$, $v=9$, almost $500$ functions were found. However, no results were stated regarding the extremality of these candidate functions. In fact, it was not until the breakthrough algorithmic results  in \cite{basu-hildebrand-koeppe:equivariant} that an automated test for extremality for the infinite group problem became possible.  Despite the lack of automated test, Chen constructed the first parametric family of 4-slope extreme functions\footnote{The function is available in the Electronic Compendium \cite{electronic-compendium} as \sagefunc{chen_4_slope}.} for the infinite group problem in \cite{chen}. 

Gomory--Johnson, in their 2003 paper \cite[section 6.3]{tspace}, had written:
\begin{quote}
  [W]e have
  discussed mainly large facets, families of two-slope facets or three-slope
  facets. What about much smaller facets? It seems extremely likely that
  there are much smaller more complex facets all around the big ones. 
\end{quote}
However, for a brief period, it seemed that with Chen's functions the peak of
complexity had been reached, as no extreme functions with more than 4 slopes
were found. In 2009, Dey--Richard in \cite{dey-richard-2009-slides-gomory-functions}
stated as an open question whether there exist extreme functions with more
than 4~slopes.  The ``4-slope conjecture'', asserting that 4~slopes is the
maximum for continuous piecewise linear extreme functions, circulated among
researchers for a while.  The conjecture reflected the hope that, in contrast
to the ``unstructured and arithmetic'' finite group problem, the complexity of
the extreme functions of the infinite group problem would be under control.  
It is unclear how much support the conjecture had at any time
; Dey (2016, Personal Communication) remembers that he did not
believe in this conjecture.

The 4-slope conjecture had to be revised quickly when new computational tools
for testing extremality became available. 
Basu et al.~\cite[Theorem 1.5]{basu-hildebrand-koeppe:equivariant} (see also
\cite[Theorem 8.5]{igp_survey_part_2}) showed that in order to test extremality of $\pi$ for
$R_f(\R/\Z)$, one simply needs to test extremality of its restriction to $\frac{1}{4q} \Z$
for $R_f(\tfrac{1}{4q} \Z/\Z)$. (Later it was shown that restricting to $\frac{1}{3q}$
suffices \cite[Theorem 8.6]{igp_survey_part_2}.)  This extremality test can be
done by computing the rank of a finite-dimensional system of linear equations.

These algorithmic results in \cite{basu-hildebrand-koeppe:equivariant} enabled
Hildebrand to run a new computer-based search, using Matlab programs (2013,
unpublished, reported in \cite[Table 4]{igp_survey}).  Like Chen \cite{chen}, 
Hildebrand generated functions $\pi$ on the $q 
\times v$ grid.  However, instead of searching exhaustively, he generated the
functions randomly.  For each candidate function~$\pi$, he checked if
$\pi|_{\frac{1}{q}\Z}$ is minimal and extreme for 
$R_f(\tfrac{1}{q} \Z/ \Z)$ (a necessary condition), and finally tested if $\pi|_{\frac{1}{4q} \Z}$ is
extreme for $R_f(\tfrac{1}{4q} \Z/\Z)$ (Basu et al.'s necessary and sufficient
condition).
In this way, Hildebrand discovered the first 5-slope extreme
functions\footnote{The functions are available in the Electronic Compendium
  \cite{electronic-compendium} as
  \sage{hildebrand\underscore{}5\underscore{}slope...}}, thus refuting the
4-slope conjecture. 

Still it seemed possible that a version of this conjecture might be true;
either with $4$ replaced by another small number, or that the 4-slope
conjecture holds ``generically'' (see \cite[Open Question
2.16]{igp_survey_arxiv_v1}).


\subsection{New search strategies; outline of the paper}

In this paper, we develop new search strategies that aim to find extreme
functions for the infinite group problem with many different slope values or
with special properties.  Our goals are markedly different from those of some
earlier research surveyed above, in particular from studies using variants of the
shooting experiment \cite{gomory2003corner,dash06:shooting}. We merely wish to
settle theoretical questions regarding the structure of extreme functions.  We
make no claims whatsoever regarding the computational usefulness of the
functions that are found by our search.


Our implementation is based on the software
\cite{infinite-group-relaxation-code}, which 
implements an automated extremality test, following the ideas of the proof of
\cite[Theorem 1.3]{basu-hildebrand-koeppe:equivariant}. The practical
implementation, described in detail in
\cite{hong-koeppe-zhou:software-abstract,hong-koeppe-zhou:software-paper},
has an empirical running time that does not strongly depend on $q$, and therefore is suitable
for functions with extremely large $q$, even though no theoretical worst-case running time bound better than polynomial in $q$ is available.

Like Gomory \cite{gom} and Evans \cite{evans-thesis}, and in contrast to
Chen's and Hildebrand's approaches, our strategies only discretize the
breakpoints into $\frac{1}{q}\Z$. 
This is crucial to be able to reach larger values of~$q$, as we explain in
\autoref{sec:limitation-of-grid}, where we discuss the complexity of search
approaches based on $q \times v$ grid discretization.

\subsubsection{Section~\ref{sec:vertex-filtering-search}: Vertex filtering search}
As a first step of our study, we investigate how far we get by just combining
state-of-the-art vertex enumeration software 
with the automated extremality
test provided by the software \cite{infinite-group-relaxation-code}.
We build on state-of-the-art vertex enumeration software, Parma
Polyhedra Library (PPL) \cite{ppl-paper} and lrslib \cite{avis1998computational, avis2000revised}.
The vertex filtering search was able to enumerate
extreme functions with $q \leq 27$, among which \textbf{the first 6-slope extreme functions} were found, breaking the previous
record of $5$ slopes due to Hildebrand (2013, unpublished; reported in
\cite{igp_survey}). 
From the results obtained by this search, we observe:
\begin{itemize}
\item a diminishing fraction of vertices of $\Pi_f(\frac{1}{q}\Z/\Z)$ that correspond to extreme functions for $R_f(\R/\Z)$;
\item an exponential growth of time spent on vertex enumeration;
\end{itemize} when $q$ increases.

\subsubsection{Section~\ref{sec:mip_approach}: Search using MIP }
These factors suggest that one needs to consider other search strategies to reach
larger $q$.  We investigate \textbf{search strategies that, for the first time, are
guided directly by the subtle structure of minimal functions that were exposed
by the proof of the algorithmic extremality test}
in~\cite{basu-hildebrand-koeppe:equivariant}, rather than using the
extremality test merely as a black box.  

To this end, we review the notions of
the two-dimensional polyhedral complex $\Delta\P$ and of additive faces 
in \autoref{sec:2d-complex-painting}. Identifying the additive faces of
$\Delta\P$ is a crucial step in the algorithmic extremality test
\cite{basu-hildebrand-koeppe:equivariant}.  By means of Gomory--Johnson's 
celebrated Interval Lemma, additive faces give rise to ``affine-imposing''
intervals (in the terminology of \cite{basu-hildebrand-koeppe:equivariant}).
This reduces the infinite-dimensional test to a finite-dimensional one.
The combinatorics of the additive faces of the complex $\Delta\P$ has a central role in our new
approaches.  

In the remainder of \autoref{sec:mip_approach} we describe how the search based on the
combinatorics of the additive faces can be implemented using standard MIP
modeling techniques and running a commercial MIP solver.  This is easy to
implement and easy to tailor to a search for extreme functions with particular
properties.  However it is limited because floating-point implementations are
not a good match for finding functions of high arithmetic complexity, and
because MIP solvers are generally not the best tool for performing an
exhaustive search.

\subsubsection{Section~\ref{sec:backtracking-search}: Backtracking search}

To address the shortcomings of the MIP approach, we have developed a
specialized backtracking search algorithm, which we describe in
\autoref{sec:backtracking-search}. 
Like the MIP approach, it is based on the combinatorics of the additive faces.
It works best when combined with the vertex filtering search described in
\autoref{sec:vertex-filtering-search}. 
This \emph{combined search algorithm} looks for extreme functions by
backtracking on the additive faces of $\Delta\P$ in a first step
and vertex enumeration in a second step. The synergy and balance 
between these two steps to obtain the best computational performance is discussed in
\autoref{sec:combined-mode}.
Using the combined search algorithm, we discover new \textbf{extreme functions with up to 7
slopes}. 

\subsubsection{Section~\ref{sec:targeted_search}: Targeted search with patterns}

In the library of functions found by the above algorithms, we observe some
special combinatorial patterns on their two-dimensional
polyhedral complexes $\Delta\P$. In \autoref{sec:targeted_search}, we describe
how we use these patterns to make a targeted search for functions with a very
large number of slopes, which discovers piecewise linear \textbf{extreme
  functions with up to 28 slopes}.

\subsection{Summary of new results}

\begin{figure}
\centering
\includegraphics[width=\linewidth]{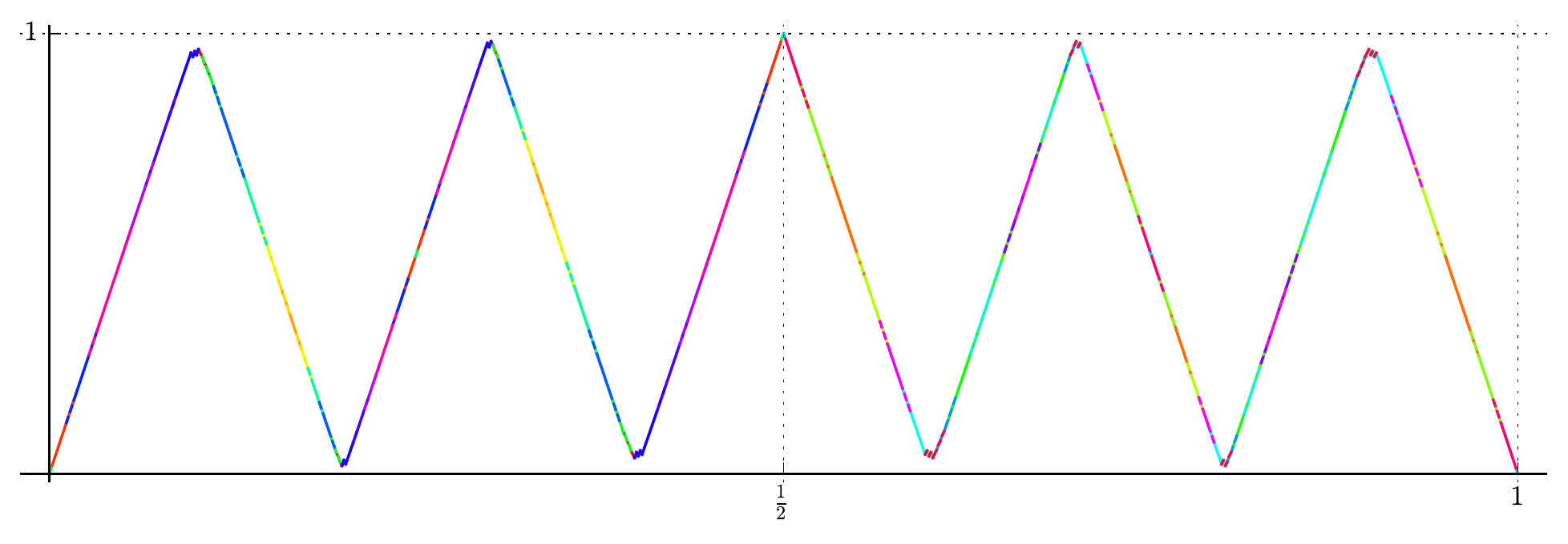}
\caption{A 28-slope extreme function \sagefunc{kzh_28_slope_1} found by our search code. Each color in the plotting corresponds to a different slope value.}
\label{fig:kzh_28_slope}
\end{figure}

\begin{theorem}
\label{thm:exists28slope}
There exist continuous piecewise linear extreme functions with $2$, $3$, $4$,
$5$, $6$, $7$, $8$, $10$, $12$, $14$, $16$, $18$, $20$, $22$, $24$, $26$, and $28$ slopes.
\end{theorem}
\autoref{fig:kzh_28_slope} shows one 28-slope extreme function found by our code, with $q=778$, out of reach for any previous study.
We remark that Basu et al.~\cite{bcdsp:arbitrary-slopes} have
improved~\autoref{thm:exists28slope} by constructing a family of extreme functions
with an arbitrary prescribed number of slopes.

\smallbreak

\begin{figure}
\centering
\includegraphics[width=\linewidth]{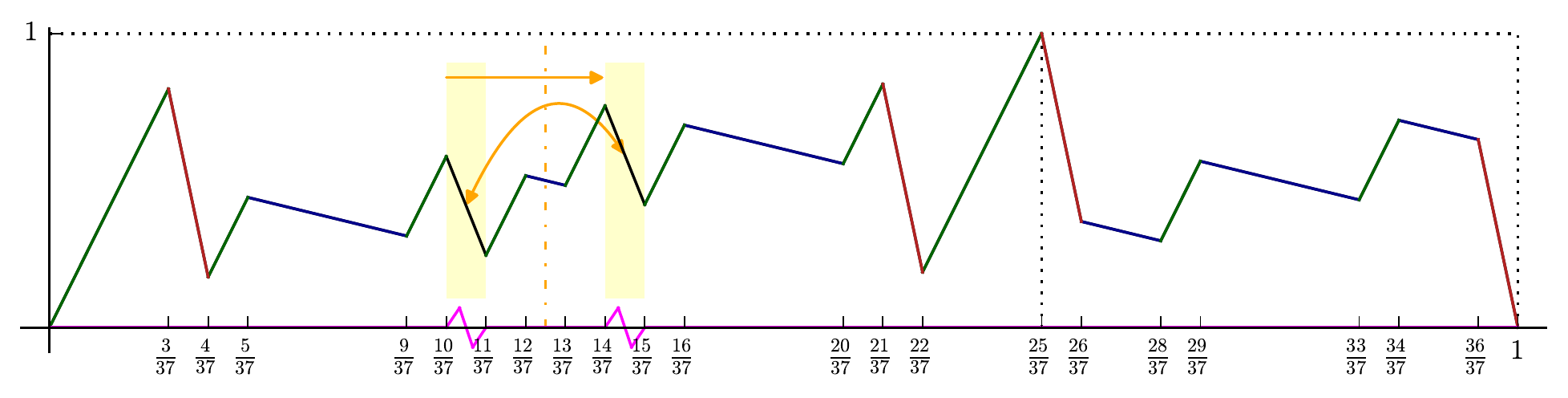}
\caption{The example \sagefunc{kzh_2q_example_1}, showing that an oversampling factor of $m=3$ in \cite[Theorem 8.6]{igp_survey_part_2} is best possible.}
\label{fig:kzh_2q_move}
\end{figure}

Our computer-based search also can be tailored to find extreme functions with certain properties. In particular, several open questions are resolved by such newly discovered extreme functions. Let $m\geq 3$ be a positive integer. \cite[Theorem 8.6]{igp_survey_part_2} states that $\pi$ is extreme for $R_f(\R/\Z)$ if and only if the restriction $\pi|_{\frac{1}{mq}\Z}$ is extreme for the finite group problem $R_f(\tfrac{1}{mq} \Z/ \Z)$. Our search found a function 
(see \autoref{fig:kzh_2q_move}) that is not extreme for $R_f (\R/\Z)$, but whose restriction to $\frac{1}{2q}\Z$ is extreme for $R_f (\frac{1}{2q}\Z/\Z)$.
This proves the following result, thereby answering the Open Question 8.7 in \cite{igp_survey_arxiv_v1}.
\begin{prop}
\label{prop:oversampling3optimal}
The hypothesis $m \geq 3$ in \cite[Theorem 8.6]{igp_survey_part_2} is best possible. The theorem does not hold for $m = 2$.
\end{prop}
The search also found piecewise linear extreme functions of $R_f(\R/\Z)$ to
answer the Open Question 2.16 in \cite{igp_survey_arxiv_v1}, which we
mentioned in \autoref{sec:intro_infinite_literature}.


\subsection{Available software}
We have made all of the discovered functions mentioned in this paper available
as part of the Electronic Compendium \cite{electronic-compendium}.  The reader
is invited to investigate the functions using our software
\cite{infinite-group-relaxation-code}. 
The articles
\cite{hong-koeppe-zhou:software-abstract,hong-koeppe-zhou:software-paper}
describe the software in detail. 
The computer-based search code will be released
shortly as part of a new version of the software
\cite{infinite-group-relaxation-code}.

\section{Restriction to $q$ grid -- Vertex filtering search}
\label{sec:vertex-filtering-search}
Recall that we are looking for continuous piecewise linear functions $\pi
\colon \R \rightarrow \R_+$ with breakpoints in $\tfrac{1}{q} \Z$ that are
extreme for the single-row Gomory--Johnson infinite group problem. The
construction of parametric families of extreme functions (a focus of many
previous studies), 
extreme functions with irrational breakpoints (for example, the function \sagefunc{bhk_irrational} in \cite{basu-hildebrand-koeppe:equivariant}), and non--piecewise linear extreme functions such as \sagefunc{bccz_counterexample} \cite{bccz08222222} are beyond the scope of this paper. 

\subsection{Restriction to grid}
Our approach is based on the discretization of the breakpoints of $\pi$. More precisely, we only focus on the functions $\pi$ with rational breakpoints  in $\tfrac1q\Z$ for some $q \in \N$. Suppose without loss of generality \cite[Lemma 2.4]{basu-hildebrand-koeppe:equivariant} that $f \in \tfrac{1}{q} \Z$.  Under such hypotheses, $\pi$ is uniquely determined by its values at points in  $\tfrac1q\Z$. We say that $\pi$ is the \emph{(continuous) interpolation} of $\pi|_{\frac{1}{q}\Z}$, while $\pi|_{\frac{1}{q}\Z}$ is the \emph{restriction} of $\pi$ to the grid $\tfrac1q\Z$. \autoref{fig:gj2s_restriction_interpolation} in \autoref{sec:introduction} illustrates the interpolation and restriction of a \sagefunc{gj_2_slope} function with $q=5$.

Gomory and Johnson proved the following relations between $\pi$ and $\pi|_{\frac{1}{q}\Z}$:
\begin{theorem}[{\cite{infinite}; see also \cite[Theorem 8.3]{igp_survey_part_2}}]
\label{thm:GJ-restrictions}
Let $\pi$ be a continuous piecewise linear function with breakpoints in
$\tfrac{1}{q} \Z$ for some $q \in \Z_+$ and let $f \in \tfrac{1}{q} \Z$.
 Then the following hold:
\begin{enumerate}[\rm(1)]
\item $\pi$ is minimal for $R_f(\R/\Z)$ if and only if $\pi|_{\frac{1}{q}\Z}$ is minimal for $R_f(\tfrac{1}{q} \Z/\Z)$.
\item If $\pi$ is extreme for $R_f(\R/\Z)$, then $\pi|_{\frac{1}{q} \Z}$ is extreme for $R_f(\tfrac{1}{q} \Z/\Z)$.
\end{enumerate}
\end{theorem}
Hence the interpolations of those $\pi|_{\frac{1}{q} \Z}$ that are extreme for
$R_f(\tfrac{1}{q} \Z/\Z)$ are the only possible candidates of continuous
piecewise linear extreme functions with breakpoints in $\frac1q\Z$ for $R_f(\R/\Z)$.
Extreme functions are clearly minimal functions; the latter have a
characterization given by a theorem by Gomory and Johnson.  In the case of the
finite group problem $R_f(\tfrac1q\Z/\Z)$, it can be stated as follows.
\begin{theorem}[{\cite{infinite}; see also \cite[Theorem 2.6]{igp_survey}}]
\label{thm:finite-minimal}
Let $\pi$ and $f$ be as above. $\pi|_{\frac{1}{q} \Z}$ is minimal for $R_f(\tfrac1q\Z/\Z)$ if and only if
\begin{enumerate}[\rm(1)]
\item $\pi_0 = 0$,
\item \label{itm:subadditive} $\pi|_{\frac{1}{q} \Z}$ is \emph{subadditive}: $\pi_{(x+y) \bmod q} \leq \pi_x + \pi_y$ for $x,y \in \Z$,
\item \label{itm:symmetric} $\pi|_{\frac{1}{q} \Z}$ is \emph{symmetric}: $\pi_x + \pi_{qf - x} = 1$ for $x \in \Z$,
\end{enumerate}
where $\pi_i = \pi(\frac{i}{q})$ for $i \in \Z$.
\end{theorem}
Since $\pi \colon \R/\Z \to \R$ is periodic modulo $1$, a minimal function $\pi|_{\frac{1}{q}\Z}$ for the finite group problem is specified by its values $(\pi_0, \pi_1, \dots, \pi_{q-1})$ on the grid points $\tfrac1q\Z \cap [0,1)$. The following statement immediately follows from the observation that the above conditions are all linear constraints.
\begin{prop}[{\cite[Theorem 2.2]{infinite}}]
\label{prop:minimal-polytope}
The set $\Pi_f(\frac{1}{q}\Z/\Z)$ of minimal functions for $R_f(\tfrac1q\Z/\Z)$ is a convex polytope.
Furthermore, extreme functions for $R_f(\tfrac1q\Z/\Z)$ are the extreme points (i.e., vertices) of this polytope.
\end{prop}

By \autoref{thm:GJ-restrictions} and \autoref{prop:minimal-polytope}, all continuous piecewise linear extreme functions for $R_f(\R/\Z)$ with breakpoints in $\frac{1}{q}\Z$ can be found by interpolating the vertices of the polytope $\Pi_f(\tfrac1q\Z/ \Z)$. However, in general, extremality of $\pi|_{\frac{1}{q} \Z}$ for $R_f(\frac{1}{q} \Z/ \Z)$ does not imply extremality of $\pi$ for $R_f(\R/\Z)$.
This makes further filtering necessary, which we explain below.


\subsection{Preprocessing}
\label{sec:preprocessing}
We begin with a remark on the minimal H-represen\-ta\-tion of $\Pi_f(\frac{1}{q}\Z/\Z)$.
The H-representation of the polytope $\Pi_f(\frac{1}{q}\Z/\Z)$ defined by \autoref{thm:finite-minimal} has asymptotically $\frac{1}{2}q^2$ constraints, many of which are redundant. Indeed, \cite[Corollary 2.7]{shim2013cyclic} gives a minimal representation of $\Pi_f(\frac{1}{q}\Z/\Z)$ that only has asymptotically $\frac{1}{6} q^2$ constraints, mainly by replacing the subadditivity constraints \eqref{itm:subadditive} of \autoref{thm:finite-minimal}:
\[\pi_i + \pi_j \geq \pi_{(i+j) \bmod q} \text{ for } 0 \leq i \leq j < q \]
with the triple system: 
\[\pi_i + \pi_j + \pi_k \geq 1 \text{ for } 0 \leq i \leq j \leq k < q,\; i + j + k = qf \pmod q.\]
A minimal H-representation is of interest for vertex enumeration, because having many redundant inequalities may greatly slow down the vertex enumeration process. Although a minimal H-representation is known for the polytope $\Pi_f(\frac{1}{q}\Z/\Z)$, the search strategies described in later sections of the present paper also need to deal with other polytopes whose minimal H-representations are not known. 
For this purpose, computational preprocessing is used to remove the redundant inequalities
from the H-representation of a polytope before enumerating its
vertices. Namely, we apply the preprocessing program \sage{redund} provided by
lrslib (version 5.0\footnote{Available from \url{http://cgm.cs.mcgill.ca/~avis/C/lrs.html}.}),
which removes redundant inequalities using Linear Programming. The third
columns of Tables \ref{tab:vertex-enumeration-without-prep}
and~\autoref{tab:vertex-enumeration-with-prep} show the number of
inequalities in the H-representation of the polytope
$\Pi_{f}(\frac{1}{q}\Z/\Z)$ for $f=\frac{1}{q}$ before and after preprocessing, respectively. 
The number after preprocessing is roughly $\frac{1}{6} q^2$, which is
consistent with \cite[Corollary 2.7]{shim2013cyclic}. 

\subsection{Performance of various vertex enumeration codes}
\label{sec:vertex_enumeration_codes}
Various software packages are available for computing the vertices of a polytope
given by an H-representation.  We considered the following popular packages.
\begin{itemize}
\item cddlib (version 094g\footnote{Available from \url{http://www.inf.ethz.ch/personal/fukudak/cdd\_home/cdd.html}.}),
  a well-known implementation of the double description method
  \cite{fukuda1996dd}. 
\item Parma Polyhedra Library (version
  1.1\footnote{Available from \url{http://bugseng.com/products/ppl/}.}).
  PPL is a C++ library for the manipulation and computation of rational convex
  polyhedra \cite{ppl-paper}.  Like cddlib, polyhedral computations in PPL are based on the double
  description method.

  An extensive computational study \cite[section 4]{ppl-paper} showed that the
  double description method implementation in PPL has a better performance
  (on the vertex/facet enumeration problem) compared with cddlib 
  and with other polyhedra libraries that are popular in PPL's primary application domain,
  New Polka and PolyLib.
\item Porta (version
  1.4.1\footnote{Available from \url{http://www.iwr.uni-heidelberg.de/groups/comopt/software/PORTA/}.}),
  based on the Fourier--Motzkin elimination method. 
\item lrslib (version 5.0
),  a C
  implementation of the lexicographic reverse search algorithm for vertex
  enumeration and convex hull problems \cite{avis1998computational, avis2000revised}. 
  This algorithm uses little memory space during the computation; vertices are
  generated as a stream and are not stored in memory, which makes it suitable
  for vertex enumeration problems for polytopes with a large number of
  vertices. 

  The computational study \cite[Table 2]{ppl-paper} showed that
  lrslib outperforms PPL for large problems, whereas PPL outperforms lrslib for easy problems
  \cite[Table 1]{ppl-paper}.
\end{itemize}
Based on the study \cite{ppl-paper}, we expected PPL to be the best choice for
low dimensions and lrslib to be the best choice for high dimensions.  
We decided to verify this using a computational study for our polytopes $\Pi_f(\tfrac{1}{q}\Z/\Z)$
given in \autoref{prop:minimal-polytope} for various values of $q$ and
$f=\tfrac{1}{q}$, recording the running times for vertex enumeration.

In addition to the four packages listed above, we also included the following
package in our experiments, which was developed very recently and had not been
investigated in a major computational study.
\begin{itemize}
\item PANDA (version 2015-02-24\footnote{Available from 
    \url{http://comopt.ifi.uni-heidelberg.de/software/PANDA/}.})
  based on the parallel adjacency decomposition algorithm
  \cite{Loerwald-Reinelt-2015-PANDA}.
\end{itemize}
At the time of the revision of the present paper, we also became aware of the high
performance vertex enumeration code in the following package.
\begin{itemize}
\item Normaliz (version
  3.1.1\footnote{Available from \url{https://www.normaliz.uni-osnabrueck.de/}.}), a software
  tool for the computation of Hilbert bases and enumerative data of rational
  cones and affine monoids. Its vertex enumeration code is based on the
  Fourier--Motzkin elimination method and the pyramid decomposition described
  in~\cite{Bruns:2016:PPD}.
\end{itemize} 

We used PPL via its Cython-based library interface within the SageMath
\cite{sage} computer algebra system.  The other systems were tested using
their command-line executables because no library interfaces were available
for them in SageMath.\footnote{A dataset with the input files in the formats
  of the systems in our
  study is available at
  \url{https://www.math.ucdavis.edu/~mkoeppe/art/infinite-group/dataset_gomory_polyhedra_cyclic_group.zip}.}
Hence in very low dimension, there is a slight bias in 
favor of PPL due to the overhead in using the other systems via executables
and file passing. All systems were tested using exact computation modes and
using a single thread. 

Tables \ref{tab:vertex-enumeration-without-prep} and~\ref{tab:vertex-enumeration-with-prep} report for each test the size of the polytope and the running times measured in CPU seconds\footnote{The tests have been performed on a virtual machine running under the QEMU hypervisor, which reports to have access to $12$ processors running at $2.0$ GHz. However, due to the virtualization, the measured running times have a large variance between runs, though all algorithms are deterministic.}, without and with preprocessing, respectively.  
The preprocessing in \autoref{tab:vertex-enumeration-with-prep} consists of
removing redundant inequalities from the H-representation using the command
\sage{redund} provided by lrslib. 
We also measured the computational overhead of interfacing to \sage{redund} in
Python, which exceeds the actual \sage{redund} running times. 
Comparing \autoref{tab:vertex-enumeration-without-prep} and
\autoref{tab:vertex-enumeration-with-prep} shows that for dimension at most
$6$, it is best to just run PPL on the original inequalities without
preprocessing.   
For dimensions at least $7$, Normaliz is the clear winner.
The total time for preprocessing pays off when the dimension of the polytope
is at least $10$ for Normaliz, or is at least $8$ for the other systems. 
For dimension at least $13$, Normaliz is faster than the second-best code, lrslib, by more than an
order of magnitude. 
Normaliz is also the least dependent on preprocessing.\footnote{W. Bruns (Personal Communication, 2016)
  explains that Normaliz performs initial transformations that remove a large
  part of the redundancy present in these problems.} 

For the computational experiments in the remainder of this paper, we used a
combination of PPL (for low dimensions) and lrslib (for high dimensions),
using preprocessing when the dimension is at least~$8$.  We did not use
Normaliz, as we only became aware of its high performance code for vertex
enumeration at the time of revision of the present paper.


\subsection{Filtering}
As mentioned earlier, extremality of $\pi|_{\frac{1}{q} \Z}$ for $R_f(\frac{1}{q} \Z/\Z)$ does not always imply extremality of $\pi$ for $R_f(\R/\Z)$.
We call the interpolation $\pi$ of a vertex $\pi|_{\frac{1}{q} \Z}$ of the polytope $\Pi_f(\tfrac1q\Z/\Z)$ a \emph{vertex-function}. 
Once the vertices of $\Pi_f(\frac{1}{q} \Z/\Z)$ are enumerated, we can use the
automated extremality test 
implemented in the software \cite{infinite-group-relaxation-code,hong-koeppe-zhou:software-abstract,hong-koeppe-zhou:software-paper} to filter out those $\pi|_{\frac{1}{q} \Z}$ whose interpolations are not extreme for $R_f(\R/\Z)$. We will discuss in \autoref{sec:2d-complex-painting} the two-dimensional polyhedral complex $\Delta\P$ and the notion of \emph{covered intervals} that make this automated extremality test possible.
\begin{remark}
  Given that $\pi|_{\frac{1}{q} \Z}$ is a vertex of
  $\Pi_f(\frac{1}{q} \Z/\Z)$, the extremality test of $\pi$ for $R_f(\R/\Z)$
  reduces to testing whether all intervals are covered, according
  to \autoref{thm:extreme-finite-covered}.  This can be determined using the
  function \sage{generate\underscore{}uncovered\underscore{}intervals}
  of~\cite{infinite-group-relaxation-code}.
\end{remark}

We now summarize the above ideas in the following algorithm, which is referred to as ``vertex filtering mode'' in our code. The implementation uses Parma Polyhedra Library and lrslib as described in~\autoref{sec:dd-ppl}.
\begin{algorithm}
\caption{vertex filtering mode}
\label{alg:vertex-filtering-mode}
\begin{enumerate}
\item Consider the restriction of $\pi$ to the grid $\tfrac1q\Z$. \\
Define $\pi_0, \pi_1, \dots, \pi_q$ as variables, where $\pi_i = \pi(\frac{i}{q})$.
\item \label{itm:initial-polytope} Construct the polytope $\Pi_f(\tfrac1q\Z/\Z)$ of minimal functions for $R_f(\tfrac1q\Z/\Z)$, defined by \autoref{thm:finite-minimal}.
\item Enumerate the vertices $\pi|_{\frac{1}{q} \Z}$ of this polytope.
\item For every vertex $\pi|_{\frac{1}{q} \Z}$, do:
\begin{enumerate}
\item Interpolate to get $\pi$, a minimal valid function for $R_f(\R/\Z)$.
\item If the intervals $[\frac{i}{q}, \frac{i+1}{q}]$ for $i = 0, 1, \dots,
  q-1$ are all covered, then $\pi$ is extreme for
  $R_f(\R/\Z)$. \textbf{Output} function $\pi$.
\end{enumerate}
\end{enumerate}
\end{algorithm}

\subsection{Performance of the vertex filtering search}

Our search code is implemented in SageMath \cite{sage}, an open-source mathematics software system that uses Python and Cython as its primary programming languages and interfaces with various existing packages. 

Our vertex filtering search code uses the strategies described in \autoref{sec:vertex_enumeration_codes} to decide whether preprocessing is needed and which software to use for vertex enumeration. We test its performance for $q=10,11,\dots, 27$ and $f = \frac{x}{q}$ for $x=1,2,\dots, \lfloor\frac{q}{2}\rfloor$. 

\begin{figure}[tp]
\centering
\includegraphics[width=.3\linewidth]{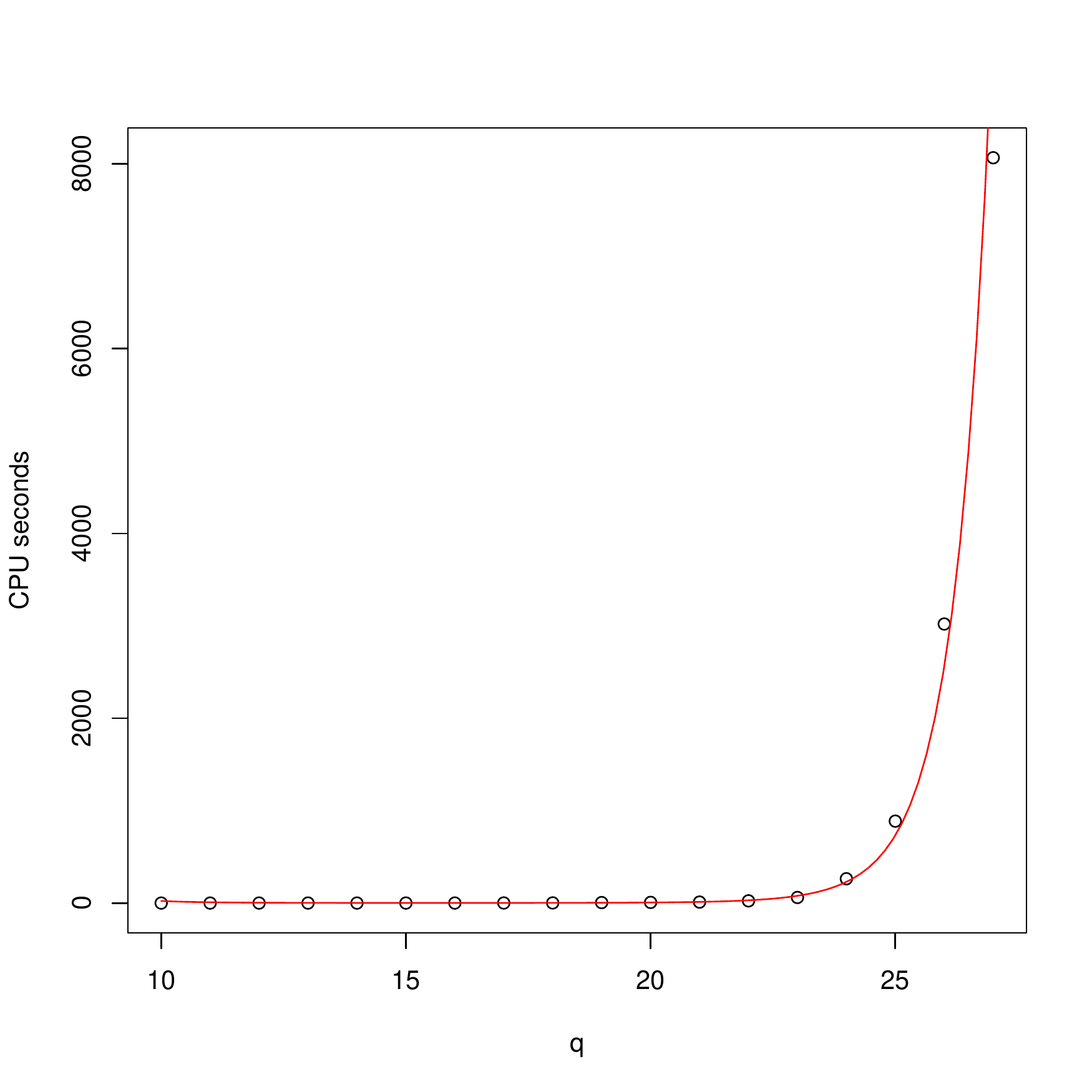}
\quad
\includegraphics[width=.3\linewidth]{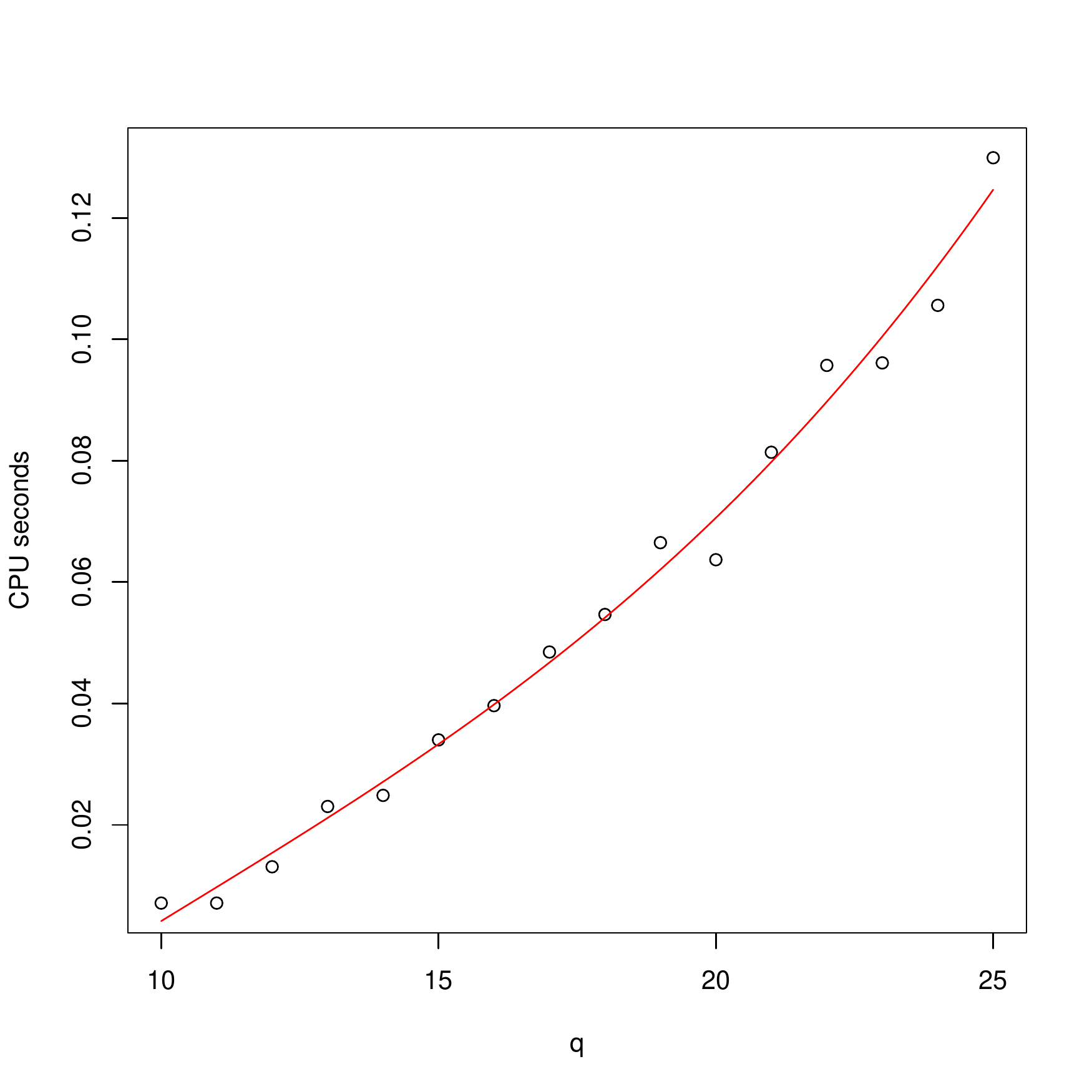}
\quad
\includegraphics[width=.3\linewidth]{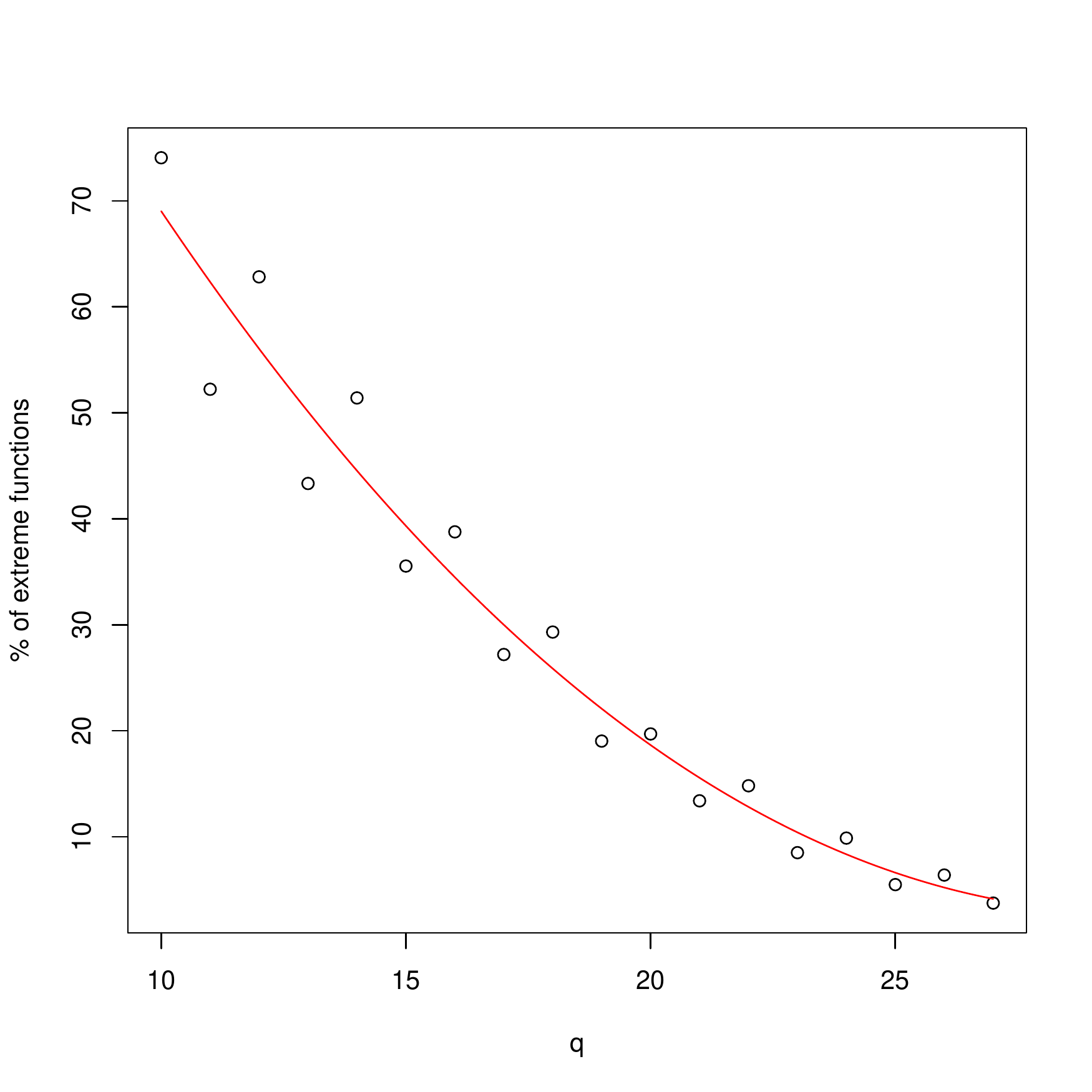}
\caption{Vertex enumeration time (not including checking extremality of vertex-functions), mean extremality checking time for a vertex-function, and percentage of extreme functions}
\label{fig:vertex-filtering-performance}
\end{figure}

Observe that as $q$ increases, the dimension and the number of vertices of the polytope increase. In particular, it results in an exponential growth of running time for vertex enumeration (cf.~\autoref{fig:vertex-filtering-performance}--left). In addition to vertex enumeration, the vertex filtering search has to run extremality tests for the vertex-functions once they are found, which consumes non-negligible extra time (cf.~\autoref{fig:vertex-filtering-performance}--middle). 
Furthermore, \autoref{fig:vertex-filtering-performance}--right illustrates a decrease in the percentage of extreme functions to vertex-functions. It suggests that when $q$ is large, vertex filtering search does enumeration in high dimension and throws away many non-extreme functions. Therefore, it is not surprising that the vertex filtering search is only suitable for small~$q$ ($ q \leq 27$). 

\subsection{Results and discussion}
\label{sec:vertex-filtering-results}
Despite its limitations, vertex filtering search finds up to 6-slope extreme functions with $q\leq 27$, breaking the previous record of 5 slopes\footnote{These functions are available in \cite{electronic-compendium} as \sage{hildebrand\underscore{}5\underscore{}slope...}} due to Hildebrand (2013, unpublished). 

The extreme functions obtained from the vertex filtering search with
breakpoints in $\frac1q\Z/\Z$ ($q \leq 27$) for the infinite and finite group
problems are plotted in \autoref{fig:denominators}, using blue and red
dots. They are placed from left to right according to the value $q$ and the
number of slopes. The log-scale $y$-axis refers to the \emph{arithmetic
  complexity} of the extreme functions, which we define as follows. 
\begin{definition}
  The \emph{arithmetic complexity}\footnote{This function is available as \sagefunc{arithmetic_complexity}.} of a function $\pi \colon \frac1q\Z/\Z \to [0,1]$
  (or of a continuous piecewise linear function
  $\pi \colon \R/\Z \to [0,1]$ with breakpoints in $\tfrac{1}{q} \Z$) is
  defined as the least common denominator of the values $\pi_i =
  \pi(\frac{i}{q})$ for $i \in \{ 0, \dots, q\}$.
\end{definition}
For example, it is easy to see that the \sagefunc{gmic} function with
$f\in\tfrac{1}{q} \Z$ has an arithmetic complexity of $O(q^2)$. 

We observe from \autoref{fig:denominators} an empirical exponential growth of
the arithmetic complexity as $q$ increases. This suggests that huge values of
$v$ would be needed in a search based on the $q\times v$ grid discretization, 
like Chen's and Hildebrand's described in
\autoref{sec:intro_infinite_literature}), in order to make them exhaustive.
This renders this approach unsuitable for large $q$.

In the following, we give an \emph{upper bound} on the arithmetic complexity of an
extreme function in terms of $q$.  We will refer to this bound in the
following sections.  Let $\pi$ be an extreme function for $R_f(\frac1q\Z/\Z)$.
\autoref{prop:minimal-polytope} states that $(\pi_0, \pi_1, \dots, \pi_{q-1})$
is a vertex of the polytope $\Pi_f(\frac{1}{q}\Z/\Z)$ defined in
\autoref{thm:finite-minimal}. 
By introducing slack variables, the constraint system of $\Pi_f(\frac{1}{q}\Z/\Z)$ can be written in the standard form using matrix notation as $A \x = \b, \x \geq \ve0$, where $A$ and $\b$ have all integer entries. Then by Cramer's rule, the denominators of
$\{ \pi_i \}_{i = 0, 1,\dots, q-1}$ come from the inverse of simplex basis matrices. 
The next lemma investigates the determinants of simplex basis matrices of $A$. The proof of this lemma appears in \autoref{sec:limitation-of-grid}.
\begin{lemma}
\label{lem:denominators-upper-bound}
Let $q \in \Z_+$ and $f \in \frac{1}{q}\Z$, $0<f<1$. Let $A \x = \b$, $\x\geq\ve0$ be the constraint system of \autoref{thm:finite-minimal} written in the standard form.
Let $\B$ be a basis matrix of $A$. Then $\mathopen|\det B| \leq 10^{{q}/{4}}$.
\end{lemma}
Given $q \in \Z_+$, let $d_{\text{ver}}$ and $d_{\text{ext}}$ denote the maximum arithmetic complexities of any extreme function $\pi$ for $R_f(\frac1q\Z/\Z)$ and of any extreme function $\pi$ for $R_f(\R/\Z)$ with breakpoints in $\frac1q\Z$,  respectively.  By \autoref{thm:GJ-restrictions} and \autoref{prop:minimal-polytope}, it is clear that the number $d_{\text{ext}}$ is well-defined and $d_{\text{ext}} \leq  d_{\text{ver}}$. Since the entries of the right-hand side $\b$ of $A \x = \b$ are integers, we have that
$$d_{\text{ext}} \leq d_{\text{ver}} \leq d_{\text{bas}} 
:= \max \{\,\mathopen|\det B| : B 
\text{ is a basis matrix of }A\,\} \leq 10^{q/4}.$$ 
We state this upper bound on $d_{\text{ver}}$ and on $d_{\text{ext}}$ as a
corollary. 
\begin{corollary}
\label{cor:arithmetic-complexity-upper-bound}
Let $\pi$ be an extreme function for $R_f(\R/\Z)$ with breakpoints in $\frac1q\Z$ (or an extreme function for $R_f(\frac1q\Z/\Z)$). Then the arithmetic complexity of $\pi$ (i.e., the least common denominator of $(\pi_0, \pi_1, \dots, \pi_q)$) is at most $10^{q/4}$.
\end{corollary}

We do not have a matching \emph{lower bound} for $d_{\text{ext}}$; however,
\autoref{lem:denominators} in \autoref{sec:limitation-of-grid} offers 
an exponential lower bound on $d_{\text{bas}}$ when an arithmetic condition
on $q$ and $f$ is satisfied.

\begin{landscape}
\begin{table}[h]
  \caption{Efficiency of various vertex enumeration codes without preprocessing}
  \label{tab:vertex-enumeration-without-prep}
  \centering 
  \begin{minipage}{\linewidth}
    \let\footnoterule=\relax
  \centering
  \begin{tabular}[t]{cccrrrrrrr}
    \toprule
    & 
    &
    &
    & \multicolumn{6}{c}{Running time (s)}
    \\
    \cmidrule(lr){5-10}
	\multicolumn{1}{c}{$q$} 
	& \multicolumn{1}{c}{dimension}
	& \multicolumn{1}{c}{inequalities}
	& \multicolumn{1}{c}{vertices}
	& \multicolumn{1}{c}{PPL}
	& \multicolumn{1}{c}{Porta}
	& \multicolumn{1}{c}{cddlib}
	& \multicolumn{1}{c}{lrslib}
	& \multicolumn{1}{c}{Panda}
	& \multicolumn{1}{c}{Normaliz}
	\\
    \midrule
	\hphantom{0}5
	& \hphantom{0}1 
	& \hphantom{0}21
	& \hphantom{0000}2
	& 0.001
	& 0.018 
	& 0.009 
	& 0.008
	& 0.026
	& 0.003
	\\
	\hphantom{0}7
	& \hphantom{0}2
	& \hphantom{0}36 
	& \hphantom{0000}4
	& 0.001 
	& 0.012
	& 0.011 
	& 0.005 
	& 0.026
	& 0.004
    \\
    \hphantom{0}9
    & \hphantom{0}3 
    & \hphantom{0}55
    & \hphantom{0000}7
    & 0.002
    & 0.016
    & 0.018 
    & 0.004
    & 0.065
    & 0.005
    \\
    11 
    & \hphantom{0}4
    & \hphantom{0}78
    & \hphantom{000}18
    & 0.003
    & 0.016 
	& 0.031
    & 0.009 
    & 23\hphantom{.000}
    & 0.007
    \\
    13
    & \hphantom{0}5
    & 105 
    & \hphantom{000}40 
    & 0.007 
    & 0.018 
    & 0.11\hphantom{0}
    & 0.021 
    & 4604\hphantom{.000} 
    & 0.011
    \\
    15
    & \hphantom{0}6
    & 136 
    & \hphantom{000}68
    & 0.017
    & 0.037
    & 0.21\hphantom{0}
    & 0.14\hphantom{0}
    & \relax
    & 0.017
    \\
    17
    & \hphantom{0}7
    & 171 
    & \hphantom{00}251
    & 0.14\hphantom{0} 
    & 0.20\hphantom{0} 
    & 1.2\hphantom{00}
    & 0.71\hphantom{0}
    & \relax
    & 0.047
    \\
    19
    & \hphantom{0}8
    & 210
    & \hphantom{00}726
    & 0.91\hphantom{0} 
    & 1.6\hphantom{00} 
    & 5.0\hphantom{00}
    & 2.3\hphantom{00} 
    & \relax
    & 0.16\hphantom{0}
    \\
    21
    & \hphantom{0}9
    & 253 
    & \hphantom{0}1661
    & 6.6\hphantom{00} 
    & 13\hphantom{.000}
    & 24\hphantom{.000} 
    & 13\hphantom{.000} 
    & \relax
    & 0.67\hphantom{0}
    \\
    23\rlap{\footnote{By isomorphism, this vertex enumeration problem ($q = 23$, $f = \frac{1}{23}$) is the same as the problem with $q = 23$ and $f = \frac{22}{23}$. The latter was tested by L.~Evans \cite[Table 6]{evans-thesis} using her parallel C implementation of the double description method, reporting a running time of $9.58$ hours (ca.\ 34500~s) on one processor and $0.75$ hours on $32$ processors, each a 550MHz Pentium III Xeon, on the Jedi cluster of the Interactive High Performance Computing Cluster at Georgia Tech.}}
    & 10
    & 300
    & \hphantom{0}7188
    & 166\hphantom{.000}
    & 558\hphantom{.000}
    & 785\hphantom{.000}
    & 74\hphantom{.000} 
    & \relax
    & 4.9\hphantom{00}
    \\
    25
    &11
    & 351
    & 23214 
    & 1854\hphantom{.000}
    & 10048\hphantom{.000}
    & 12129\hphantom{.000}
    & 471\hphantom{.000}
    & \relax
    & 21\hphantom{.000}
    \\
    26
    & 12
    & 378
    & 54010
    & \relax
    & \relax
    & \relax
    & 2167\hphantom{.000}
    & \relax
    & 62\hphantom{.000}
    \\
    27
    & 12
    & 406
    & 68216
    & \relax
    & \relax
    & \relax
    & \relax
    & \relax
    & 89\hphantom{.000}
    \\
    28
    & 13
    & 435
    & 195229
    & \relax
    & \relax
    & \relax
    & \relax
    & \relax
    & 326\hphantom{.000}
    \\
    29
    & 13
    & 465
    & 317145 
    & \relax
    & \relax
    & \relax
    & \relax
    & \relax
    & 644\hphantom{.000}
    \\
    30
    & 14
    & 496
    & 576696 
    & \relax
    & \relax
    & \relax
    & \relax
    & \relax
    & 1693\hphantom{.000}
    \\
    31
    & 14
    & 528
    & 1216944 
    & \relax
    & \relax
    & \relax
    & \relax
    & \relax
    & 3411\hphantom{.000}
    \\
    \bottomrule
  \end{tabular}
\end{minipage}
\end{table}
\end{landscape}

\begin{landscape}
\begin{table}[h]
  \caption{Efficiency of various vertex enumeration codes with preprocessing}
  \label{tab:vertex-enumeration-with-prep}
  \centering 
  \begin{minipage}{\linewidth}
    \let\footnoterule=\relax
  \centering
  \begin{tabular}[t]{cccrrrrrrrrr}
    \toprule
    & 
    &
    &
    & \multicolumn{8}{c}{Running time (s)}
    \\
    \cmidrule(lr){5-12}
	\multicolumn{1}{c}{$q$} 
	& \multicolumn{1}{c}{dim} 
	& \multicolumn{1}{c}{ineqs}
	& \multicolumn{1}{c}{vertices} 
	& \multicolumn{1}{c}{PPL}
	& \multicolumn{1}{c}{Porta}
	& \multicolumn{1}{c}{cddlib}
	& \multicolumn{1}{c}{lrslib}
	& \multicolumn{1}{c}{Panda}
	& \multicolumn{1}{c}{Normaliz}
	& \multicolumn{1}{c}{\sage{redund}}
	& \multicolumn{1}{c}{overhead}
	\\
    \midrule
	\hphantom{0}5
	& \hphantom{0}1 
	& \hphantom{00}7
	& 2
	& 0.003
	& 0.009
	& 0.006
	& 0.010
	& 0.019
	& 0.003
	& 0.006 
	& 0.040
    \\
	\hphantom{0}7
	& \hphantom{0}2
	& \hphantom{0}10
	& 4
	& 0.001
	& 0.010
	& 0.009
	& 0.007
	& 0.015
	& 0.003
	& 0.006 
	& 0.029
    \\
    \hphantom{0}9
    & \hphantom{0}3 
    & \hphantom{0}14
    & 7
    & 0.001
    & 0.008
    & 0.009
    & 0.008
    & 0.021
    & 0.004
    & 0.009 
    & 0.010
    \\
    11 
    & \hphantom{0}4
    & \hphantom{0}20
    & 18
    & 0.002 
    & 0.008
	& 0.015
    & 0.010
    & 0.017
    & 0.006
    & 0.012 
    & 0.049
    \\
    13
    & \hphantom{0}5
    & \hphantom{0}27
    & 40 
    & 0.003
    & 0.007
    & 0.021
    & 0.012
    & 0.039
    & 0.009
    & 0.022 
    & 0.050
    \\
    15
    & \hphantom{0}6
    & \hphantom{0}35
    & 68
    & 0.004
    & 0.012
    & 0.032
    & 0.025
    & 0.040
    & 0.012
    & 0.041 
    & 0.14\hphantom{0}
    \\
    17
    & \hphantom{0}7
    & \hphantom{0}45
    & 251
    & 0.016
    & 0.030
    & 0.22\hphantom{0}
    & 0.10\hphantom{0}
    & 0.16\hphantom{0}
    & 0.029
    & 0.041 
    & 0.21\hphantom{0}
    \\
    19
    & \hphantom{0}8
    & \hphantom{0}56
    & 726
    & 0.061
    & 0.087
    & 0.34\hphantom{0}
    & 0.48\hphantom{0}
    & 0.44\hphantom{0}
    & 0.085
    & 0.16\hphantom{0}
    & 0.40\hphantom{0}
    \\
    21
    & \hphantom{0}9
    & \hphantom{0}68
    & 1661
    & 0.25\hphantom{0}
    & 0.25\hphantom{0}
    & 1.1\hphantom{00}
    & 2.5\hphantom{00}
    & 3.1\hphantom{00}
    & 0.22\hphantom{0}
    & 0.25\hphantom{0} 
    & 0.72\hphantom{0}
    \\
    23
    & 10
    & \hphantom{0}82
    & 7188
    & 4.0\hphantom{00} 
    & 4.1\hphantom{00}
    & 8.0\hphantom{00}
    & 15\hphantom{.000}
    & 9.0\hphantom{00}
    & 1.3\hphantom{00}
    & 0.46\hphantom{0} 
    & 1.1\hphantom{00}
    \\
    25
    &11
    & \hphantom{0}97
    & 23214 
    & 69\hphantom{.000}
    & 43\hphantom{.000}
    & 31\hphantom{.000}
    & 94\hphantom{.000}
    & 15 h\hphantom{00} 
    & 6.0\hphantom{00}
    & 0.75\hphantom{0}
    & 1.8\hphantom{00}
    \\
    26 
    & 12 
    & 115 
    & 54010 
    & 511\hphantom{.000}
    & 350\hphantom{.000}
    & 692\hphantom{.000}
    & 594\hphantom{.000}
    & \relax 
    & 21\hphantom{.000}
    & 0.95\hphantom{0} 
    & 2.6\hphantom{00}
    \\
    27 
    & 12 
    & 113 
    & 68216 
    & 433\hphantom{.000}
    & 493\hphantom{.000}
    & 672\hphantom{.000}
    & 543\hphantom{.000} 
    & \relax
    & 23\hphantom{.000}
    & 1.0\hphantom{00} 
    & 3.0\hphantom{00}
    \\
    28 
    & 13 
    & 133 
    & 195229 
    & 8399\hphantom{.000}
    & 5796\hphantom{.000}
    & 9550\hphantom{.000}
    & 3617\hphantom{.000}
    & \relax
    & 102\hphantom{.000}
    & 1.6\hphantom{00} 
    & 4.0\hphantom{00} 
    \\
    29 
    & 13 
    & 131 
    & 317145 
    & 18361\hphantom{.000}
    & 11341\hphantom{.000}
    & \relax
    & 3366\hphantom{.000}
    & \relax
    & 158\hphantom{.000}
    & 1.9\hphantom{00}
    & 4.9\hphantom{00}
    \\
    30 
    & 14 
    & 152 
    & 576696 
    & $>1$ d\hphantom{.000} 
    & 66747\hphantom{.000}
    & \relax
    & 22743\hphantom{.000} 
    & \relax
    & 488\hphantom{.000}
    & 2.5\hphantom{00} 
    & 6.1\hphantom{00}
    \\
    31 
    & 14 
    & 150 
    & 1216944 
    & $>3$ d\hphantom{.000}
    & $>2$ d\hphantom{.000}
    & \relax 
    & 20407\hphantom{.000}
    & \relax
    & 734\hphantom{.000}
    & 2.8\hphantom{00}
    & 7.8\hphantom{00}
    \\
    \bottomrule
  \end{tabular}
\end{minipage}
\end{table}
\end{landscape}

\begin{landscape}
\begin{table}[p]
  \caption{\small The arithmetic complexity of the search based on $q \times v$ grid discretization. $d_{\text{ext}}$ and $d_{\text{ver}}$ are the empirical values of $v$ for infinite and finite group problems, and $\mathopen|\det B|$ is the estimated value from \autoref{lem:denominators}.}
  \label{tab:denomiators}
  \vspace{-0.5em}
  \begin{minipage}{\linewidth}
  \centering \small
  \begin{tabular}[t]{c*{18}r}
    \toprule 
	$q$  &  10   &  11  &  12  &  13  &  14  & 15  &  16  & 17 & 18& 19 & 20 & 21 & 22 & 23 & 24 & 25 & 26 & 27\\
    \midrule
    $d_{\text{ext}}$ &  21 & 30 & 35 & 48 & 51 & 64 & 63 & 120 & 91 & 168 & 165 & 208 & 255 & 348 & 289 & 504 & 459 & 800 	\\
    $d_{\text{ver}}$ &  21 & 30 & 35 & 48 & 51 & 70 & 65 & 138 & 95 & 210 & 165 & 250 & 315 & 570 & 425 & 768 & 651 & 1120 \\
    $\mathopen|\det B|$ & & 32 & & 64 & & 128 & & 256 & & 512 & & 1024 & & 2048 & & 4096 & & 8192 \\
    \bottomrule
  \end{tabular}
  \end{minipage}
\end{table}

\begin{figure}[p]
\centering
\begin{tikzpicture}
    \draw (0,0) node [anchor=north west]{
    {\resizebox{.9\linewidth}{6.0cm}{%
    \pgfimage[]{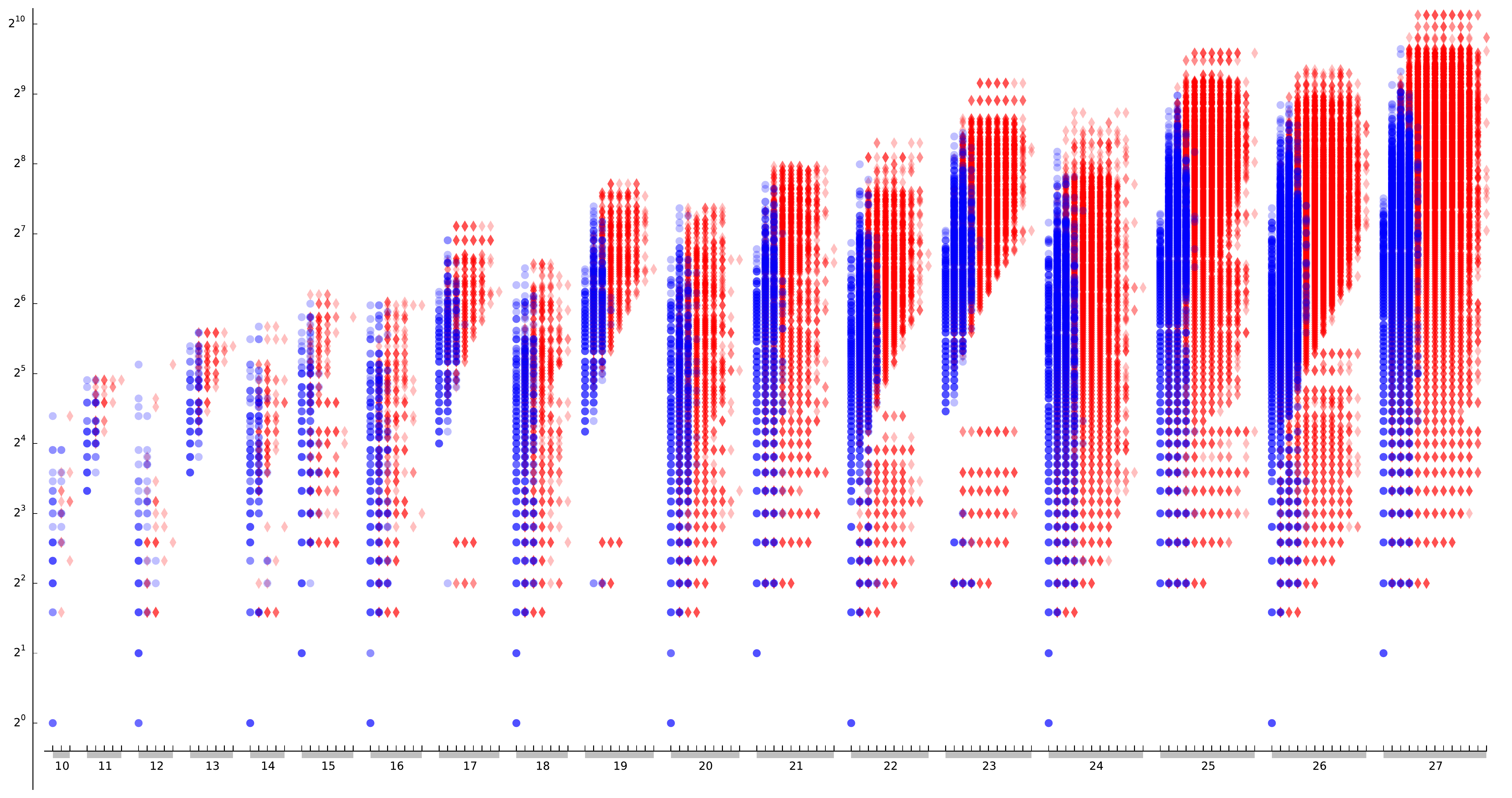}}}};
    \draw (0.8,-0.5) node [anchor=north west]{\footnotesize \textcolor{blue}{$\bullet$} extreme for $R_f(\R/\Z)$};
    \draw (0.8,-1) node [anchor=north west]{\footnotesize \textcolor{red}{\scriptsize $\blacklozenge$} extreme for $R_f(\frac{1}{q}\Z/\Z)$};
\end{tikzpicture}
  \caption{\small Arithmetic complexity and number of slopes depending on
    $q$. Extreme functions $\pi$ with breakpoints in $\frac{1}{q}\Z/\Z$ for
    $R_f(\R/\Z)$ and for $R_f(\frac{1}{q}\Z/\Z)$ are plotted in \emph{dark
      blue} and \emph{bright red},
    respectively. The $x$-axis refers to the value $q$. Within the same value
    $q$, extreme functions $\pi$ are placed in ascending order by their number
    of slopes, from left ($2$ slopes) to right. The log-scale $y$-axis refers
    to the arithmetic complexity of the extreme functions~$\pi$, i.e., the
    least common denominator of $\{\pi_0, \pi_1, \dots, \pi_{q-1}\}$, showing
    the complexity of an exhaustive search based on $q \times v$ grid
    discretization.}
  \label{fig:denominators}
\end{figure}

\end{landscape}



\section{MIP approach}
\label{sec:mip_approach}

In this section we present a new approach for computer-based search for extreme
functions.  It uses standard mixed integer linear programming (MIP) modeling
techniques to obtain a MIP that mimics the algorithmic extremality test due to
Basu et al.~\cite{basu-hildebrand-koeppe:equivariant}.

\subsection{The two-dimensional polyhedral complex $\Delta\P$}
\label{sec:2d-complex-painting}

\begin{figure}[t]
  \centering
  \includegraphics[width=.44\linewidth]{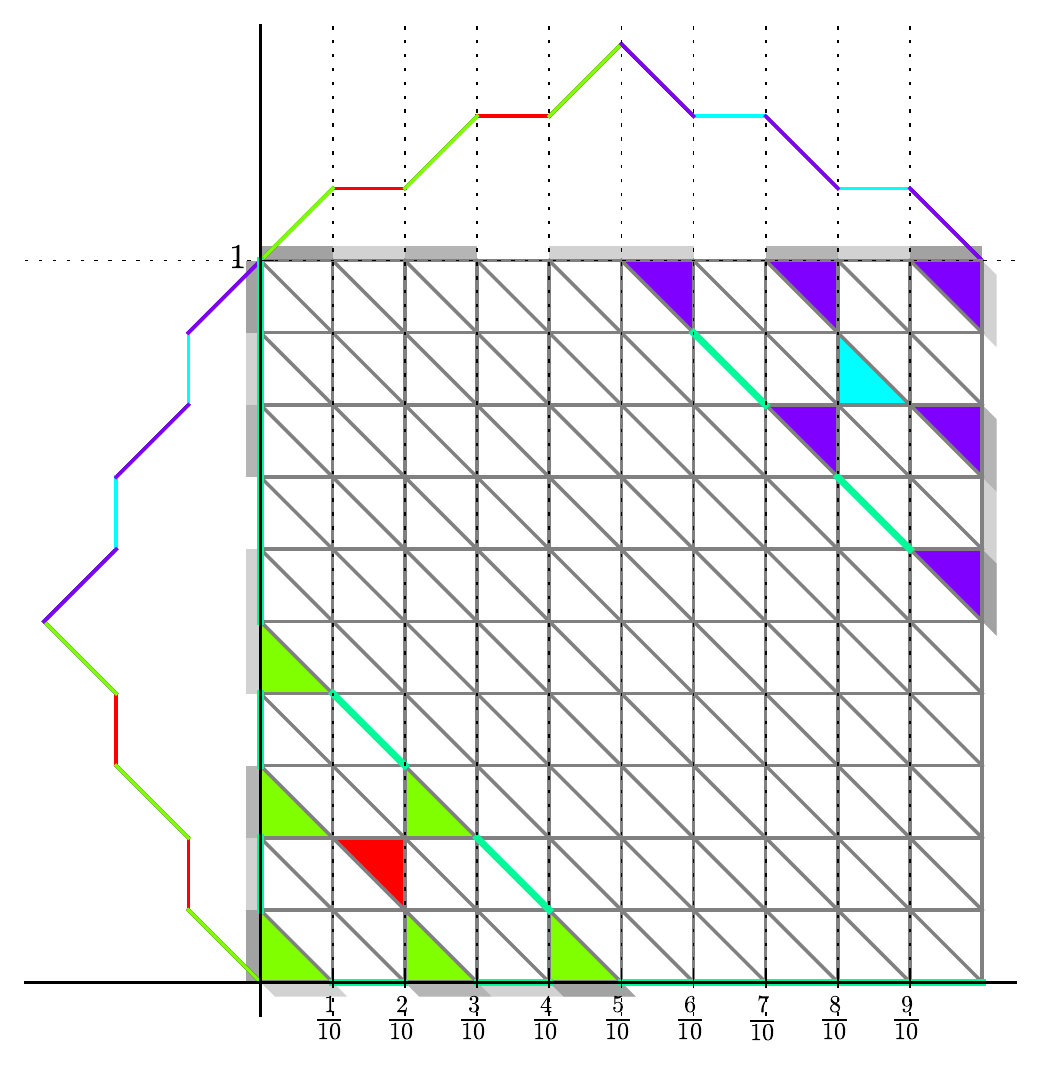}
  \caption{
    Diagram of a minimal valid function (\emph{graphs on the top and the
      left}) on the grid $\frac{1}{10}\Z$ and the corresponding painting on
    the two-dimensional polyhedral complex $\Delta\P_{\frac1{10}\Z}$
    (\emph{gray solid lines}), 
    as plotted by the command \sage{\sagefunc{plot_2d_diagram}($\pi$, colorful=True)}, where \sage{$\pi$ = \sagefunc{not_extreme_1}()}.
    Faces of $\Delta\P_{\frac1{10}\Z}$ on which $\Delta\pi = 0$, i.e.,
    additivity holds, are \emph{shaded}
    in colors that correspond to the 4 connected components of this function.
    The \emph{heavy diagonal green lines} $x + y = f$ and $x + y = 1+f$ correspond to the symmetry
    condition.  
    At the borders, the projections $p_i(F)$ of two-dimensional
    additive faces are shown as \emph{gray shadows}: $p_1(F)$ at the top border, $p_2(F)$ at
    the left border, $p_3(F)$ at the bottom and the right borders.
  }
  \label{fig:uniform-spacing}
\end{figure}

We first review the notion of a two-dimensional polyhedral complex, which serves as a tool for studying additivity relations and covered (affine imposing) intervals of piecewise linear functions. 
We follow \cite[Section 3]{igp_survey}, but define the notions in our case where the function $\pi$ is continuous, piecewise linear and has all its breakpoints in $\tfrac{1}{q}\Z$. 
This matches the setting of \cite{basu-hildebrand-koeppe:equivariant}.
Since a minimal function is periodic modulo 1, we can restrict the study to the domain $[0,1]$ only. We define
the evenly spaced one-dimensional polyhedral complex~$\P_{\frac1q\Z}$ to be the collection of singletons and elementary closed intervals on the grid $\tfrac1q\Z$ by
\[
\P_{\frac1q\Z} := \left\{\,\emptyset, \{\tfrac0q \}, \{ \tfrac1q \}, \dots, \{\tfrac{q}{q} \}, [\tfrac0q, \tfrac1q], [\tfrac1q, \tfrac2q], \dots, [\tfrac{q-1}q, 1]\,\right\}.
\]
For any $I,J,K \in \P_{\frac1q\Z}$, let
\[
F(I,J,K) := \setcond{(x,y) \in I \times J}{x \oplus y \in K} \subseteq [0, 1] \times [0, 1],
\]
where $x \oplus y = (x + y) \bmod 1$. Then the set \[\Delta\P_{\frac1q\Z} := \setcond{F(I,J,K)}{I, J, K \in \P_{\frac1q\Z}}\] is a two-dimensional polyhedral complex. It is the collection of the elementary upper or lower triangles on the grid  $\tfrac1q\Z \times \tfrac1q\Z$ with the vertices (zero-dimensional faces) and edges (one-dimensional faces) that arise as intersections of these triangles (two-dimensional faces).  See \autoref{fig:uniform-spacing} for an illustration.

Define the \emph{subadditivity slack} $\Delta\pi\colon [0,1]\times [0,1]\rightarrow \R $ of~$\pi$ by 
\[
\Delta \pi(x, y) := \pi(x) + \pi(y) -\pi(x \oplus y)
\]
for $x,y\in [0, 1]$.
Note that $\Delta \pi$ is non-negative if $\pi$ is minimal, since minimality implies subadditivity.
A face $F$ of the two-dimensional polyhedral complex $\Delta\P_{\frac1q\Z}$ is
said to be \emph{additive} if $\Delta\pi = 0$ on $F$.  Together the additive
faces form the \emph{additivity domain} of the function~$\pi$. 
Since $\pi$ is linear on the intervals of~$\P_{\frac1q\Z}$, 
the function~$\Delta\pi$ is linear on the faces of the two-dimensional complex $\Delta\P_{\frac1q\Z}$.
Hence, the condition above is equivalent to $\Delta\pi = 0$ on the set of vertices of $F$. 

In the diagrams of the polyhedral complex $\Delta\P_{\frac1q\Z}$ such as
\autoref{fig:uniform-spacing}, we indicate additive faces by various colors.
Isolated additive points and additive edges are always drawn in green;
two-dimensional additive faces (triangles) are painted in various colors, the
significance of which we shall explain below.  If a triangle is white, this means that strict subadditivity holds in the interior of this
face.
We will often refer to the diagram of the additive faces as the
\emph{painting} of~$\pi$ on the complex $\Delta\P_{\frac1q\Z}$.


An additive face implies, among other things, the important covering (\emph{affine
imposing} in the terminology of \cite{basu-hildebrand-koeppe:equivariant}) property that we outline here. We refer the reader to~\cite[Section 4]{igp_survey} for details on the Interval Lemma and its generalizations.

Define the projections $p_1,p_2,p_3\colon [0,1] \times [0,1] \to [0,1]$ by
\[
p_1(x,y) = x, \quad p_2(x,y) = y, \quad  p_3(x,y) = x \oplus y.
\]

Let $F$ be a two-dimensional additive face of $\Delta\P_{\frac1q\Z}$ (i.e.,
$F$ is an elementary upper or lower triangle in the two-dimensional polyhedral
complex such that $\Delta\pi = 0$ on $F$). By the convex additivity domain lemma for continuous functions \cite[Corollary 4.9]
{igp_survey} (a consequence of the celebrated Gomory--Johnson Interval Lemma), 
$\pi$ is affine imposing with the same slope on the projection
intervals $p_1(F)$, $p_2(F)$ and $p_3(F)$. We say that these three intervals
are \emph{(directly) covered} and \emph{connected} to each other.

Let $F$ be a one-dimensional additive face of $\Delta\P_{\frac1q\Z}$ (i.e.,
$F$ is an elementary horizontal, vertical or diagonal edge in the
two-dimensional polyhedral complex such that $\Delta\pi = 0$ on $F$). Then two
of the projections $p_1(F), p_2(F)$ and $p_3(F)$ are one-dimensional. These
two intervals are said to be \emph{connected by an edge}. An interval that is
connected to a covered interval is also said to be \textit{(indirectly)
  covered}; see \cite[Lemma 4.5]{basu-hildebrand-koeppe:equivariant}. 

The covered intervals of $\pi$ are computed in two steps. Start with directly covered intervals as $p_1(F), p_2(F)$ and $p_3(F)$ of two-dimensional additive faces $F$. Then continue transferring indirectly covered properties 
using one-dimensional additive faces until no new covered intervals are found. (This saturation process clearly ends after a finite number of steps.)

The set of covered intervals is partitioned into \emph{connected components}\footnote{The connected components are understood in a graph-theoretic sense. We refer the reader to~\cite{basu-hildebrand-koeppe:equivariant} for details on the graph of intervals that is used.}. 
In the diagrams of the polyhedral complex $\Delta\P_{\frac1q\Z}$, colors are used to indicate membership in a
connected component.  The function in \autoref{fig:uniform-spacing}, for
example, has 4~connected components, though it only has 3~slopes.  
Within a connected component, the function~$\pi$ has the same slopes. 
Thus the number of connected components gives an upper bound on the number of
slopes of the function $\pi$\footnote{Available as
  \sagefunc{number_of_components} and \sagefunc{number_of_slopes} in \cite{infinite-group-relaxation-code}.}.
\begin{remark}
  Though the number of different slopes of a function has attracted the
  attention in the past, it appears that the number of connected
  components is a more fundamental notion.
\end{remark}

A painting is called a \emph{covering painting} if all
intervals $[\frac{x}{q}, \frac{x+1}{q}]$ ($0 \leq x \leq q-1$) are
covered. 
By~\autoref{thm:extreme-finite-covered} below, every extreme function
$\pi$ has a covering painting. This property will be used as an important
ingredient in the MIP approach described in this section and also 
the backtracking search approach to be
discussed in 
\autoref{sec:backtracking-search}. 

\begin{theorem}[rephrased from results in \cite{basu-hildebrand-koeppe:equivariant}]
Let $\pi$ be a continuous piecewise linear function with breakpoints in $\tfrac{1}{q} \Z$ for some $q \in \Z_+$ and let $f \in \tfrac{1}{q} \Z$. Then $\pi$ is extreme for $R_f(\R/\Z)$ if and only if $\pi|_{\frac{1}{q} \Z}$ is extreme for $R_f(\frac{1}{q} \Z/\Z)$ and all intervals $[\frac{x}{q}, \frac{x+1}{q}]$ for $x = 0, 1, \dots, q-1$ are covered.
\label{thm:extreme-finite-covered}
\end{theorem}
Instead of giving a proof, we point the reader to the results from \cite{basu-hildebrand-koeppe:equivariant} that are rephrased as \autoref{thm:extreme-finite-covered}.
The ``if'' direction follows directly from \cite[Corollary 3.4]{basu-hildebrand-koeppe:equivariant}.
If $\pi|_{\frac{1}{q} \Z}$ is not extreme for $R_f(\frac{1}{q} \Z/\Z)$, then $\pi$ is not extreme for $R_f(\R/\Z)$ by \autoref{thm:GJ-restrictions}. If the intervals $[\frac{x}{q}, \frac{x+1}{q}]$ for $x = 0, 1, \dots, q-1$ are not all covered, then \cite[Lemma 4.8]{basu-hildebrand-koeppe:equivariant} implies the nonextremality of $\pi$ by equivariant perturbation.  This shows the ``only if'' direction by contraposition.

\subsection{Additivity variables and prescribed partial paintings}

\begin{figure}
\begin{tikzpicture}[scale=1/2,sharp corners]
\fill[draw = black, fill = darkgreen, thick]  (0,0) node[anchor=north east]{$p_{x, y}$} -- (0,1.5)node[anchor=east]{$v_{x,y}$} -- (0,3) node[anchor=south east]{$p_{x, y+1}$} -- (3,0) node[anchor=north west]{$p_{x + 1, y}$} -- (1.5, 0)node[anchor=north]{$h_{x,y}$}-- (0,0) node[circle, fill=black, inner sep=1pt]{};
\draw[](1,1) node[white]{$\ell_{x,y}$};
\draw[black,thin]  (3,3) node[anchor=south west]{$p_{x + 1, y + 1}$} -- (3,0 )node[circle, fill=black, inner sep=1pt]{}-- (0,3) node[circle, fill=black, inner sep=1pt]{} -- (3,3) node[circle, fill=black, inner sep=1pt]{};
\draw[](2,2) node[]{$u_{x,y}$};
\draw[-latex, thin](4,1.5) node[anchor=west]{$d_{x,y}$}
       to[out=180, in =0] (2,1);
\end{tikzpicture}
\caption{Binary variables for color: additive (colored \emph{dark green} in the diagram) =
  0; strictly subadditive (colored \emph{white} in the diagram) = 1}
  \label{fig:color-variables}
\end{figure}
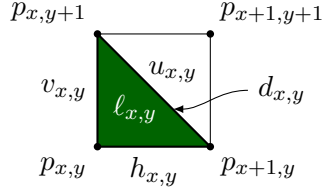

In the MIP approach we use binary variables to control the additivity
(coloredness) of a face, i.e., a triangle, an edge, or a vertex of the
two-dimensional complex $\Delta\P_{\frac1q\Z}$. 
Let $x, y$ be integers between $0$ and $q-1$. We use variables $\ell_{x,y}$ for the lower triangle whose lower left corner
is vertex $(\frac{x}{q},\frac{y}{q})$, $v_{x,y}$ for its
vertical edge, $h_{x,y}$ for its horizontal edge, $d_{x,y}$ for its diagonal
edge, $u_{x, y}$ for the upper triangle with the same diagonal edge, 
and $p_{x, y}$ for the vertex $(x,y)$. The value $0$ for these variables
represents additivity (colored face) and the value $1$ represents strict
subadditivity in the relative interior of the face (white face);  see \autoref{fig:color-variables}. These
variables are subject to the invariance of the subadditivity slack under
exchanging $x$ and $y$, hence we have 
\begin{equation}  
p_{x,y} = p_{y,x},\ \ell_{x,y}=\ell_{y,x},\ u_{x,y}=u_{y,x},\ h_{x,y}=v_{y,x} \text{ and } d_{x,y}=d_{y,x}, 
\label{eq:xy_swapped_constraints}
\end{equation}
for $x,y \in \{0, 1, \dots, q-1\}$.
Therefore, it suffices to consider only the upper left triangular part of the complex $\Delta\P_{\frac1q\Z}$ where $x \leq y$.
These binary variables are also subject to inclusion constraints: for any $x,y \in \{0, 1, \dots, q-1\}$,
\begin{subequations}%
\begin{align}
\max\{p_{x,y}, p_{x, y+1}, p_{x+1, y}\} &\leq \ell_{x,y} \leq p_{x,y} + p_{x, y+1} + p_{x+1, y}; \\
\max\{p_{x,y+1}, p_{x+1, y+1}, p_{x+1, y}\} &\leq
  u_{x,y} \leq p_{x,y+1} + p_{x+1, y+1} + p_{x+1, y}; \\ 
 \max\{p_{x,y}, p_{x, y+1}\} &\leq v_{x,y} \leq p_{x,y} + p_{x, y+1}; \\ 
 \max\{p_{x,y}, p_{x+1, y}\} &\leq h_{x,y} \leq p_{x,y} + p_{x+1, y}; \\ 
 \max\{p_{x,y+1}, p_{x+1, y}\} &\leq d_{x,y} \leq p_{x,y+1} + p_{x+1, y},
\end{align}
\label{eq:logical_constraints} 
\end{subequations}
where $p_{x, q} = p_{0,q}$, $p_{q, y}=p_{q,0}$ and $p_{q,q}=p_{0,0}$.

\begin{remark}\label{ppp}
  The subproblems obtained from branching on the values of these binary
  variables have an obvious interpretation in terms of paintings on the
  two-dimensional complex $\Delta\P_{\frac1q\Z}$, which we shall refer to as
  \emph{prescribed partial paintings}: If an additivity variable is not fixed
  yet in a branching node, the corresponding face is shown in light grey.  If it
  has been fixed to 0, the face is painted dark green.  If it has been fixed to 1, the
  face is painted white. 
 
  We shall say that a function \emph{satisfies} a prescribed partial painting
  when it would be a feasible solution to the corresponding node subproblem,
  i.e., if it satisfies all additivity conditions and strict subadditivity conditions
  corresponding to faces that have been painted dark green and white, respectively.
\end{remark}

\begin{example}
\autoref{fig:partial-painting-tree} shows a tree of prescribed partial paintings of $\Delta\P_{\frac1q\Z}$ with $q=4$ and $f=\frac14$, obtained from branching on the variables $\ell_{x,y}$ and $u_{x,y}$. In (the upper left triangular part of) the partial painting at its root node (a), the values of the following binary variables are set to $0$.
\[
 p_{0,0}, \; p_{0,1}, \; p_{0,2}, \; p_{0,3}, \; p_{2,3},\; 
 v_{0,0}, \;  v_{0,1}, \;  v_{0,2}, \;  v_{0,3}, \; 
 d_{0,0}, \;  d_{1,3}, \;  d_{2,2}.  
\]
\end{example}

\subsection{Function value variables}

The values of candidate functions are modeled by continuous variables $\pi_0,
\pi_1, \dots, \pi_{q-1} \in [0, 1]$. They are subject to the constraints in \autoref{thm:finite-minimal}.
\begin{subequations}
\begin{gather}
\pi_0 =  0 \\
\epsilon \, p_{x,y} \leq \pi_x + \pi_y - \pi_{(x + y) \bmod q} \leq 2 p_{x,y} \label{eq:strict-subadditivity}\\
\pi_x +\pi_{(qf-x) \bmod q} = 1,
\label{eq:fn_minimality_test}
\end{gather}
\end{subequations}
where $\epsilon$ is a small positive number used to enforce the strict subadditivity 
of vertices $(x, y)$ with $p_{x,y} = 1$. See \autoref{sec:epsilon-discussion} for a further discussion on the value of $\epsilon$.

\subsection{Slope value variables and assignment}
\label{sec:slope-assignment}
In a search for functions with a prescribed number $k$ of different slopes, 
we introduce $k$ continuous variables $s_1, s_2, \dots, s_k \in [-q, q]$ for the different
slope values of $\pi$. We can enforce $s_1 > s_2 > \dots > s_k$ by another artificial lower bound $\epsilon'$ (see \autoref{sec:epsilon-discussion}):
\begin{subequations}
\begin{equation}
s_j - s_{j+1} \geq \epsilon', \text{ for } 1 \leq j \leq k-1. \label{eq:slope-difference}
\end{equation}
Then binary variables $\delta_{x,j}$ ($0 \leq x \leq q-1, 1 \leq j \leq k$) are used to assign intervals to slope values:
\begin{equation}
\sum_{j = 1}^{k}{\delta_{x, j} = 1}, \text{ for } 0 \leq x \leq q-1,
\end{equation}
\[s_j = q(\pi_{x+1} - \pi_x) \text{ if and only if } \delta_{x,j} = 1, \text{
    for } 0 \leq x \leq q-1 \text{ and } 1 \leq j \leq k.\]
The last condition can be written as the linear inequality
\begin{equation}
\left| s_j + q (\pi_x - \pi_{x+1}) \right| \leq 2q (1 - \delta_{x,j}).
\end{equation}
We add the linear constraints 
\begin{equation}
\sum_{x=0}^{q-1} \delta_{x,j} \geq 1, \text{ for } 1 \leq j \leq k
\end{equation} 
to the MIP formulation, so as to ensure that every slope value $s_j$ is used by at least one interval.
\label{eq:print_slope_constraints}
\end{subequations}

\subsection{Variables for directly and indirectly covered intervals}

We use binary variables $c_{z,0}$ ($0 \leq z \leq q-1$) to control whether the
interval $[\frac{z}{q}, \frac{z+1}{q}]$ is directly covered or not: $0$ for
covered and $1$ for uncovered. They are subject to combinatorial conditions of
being directly covered by dark green triangles presented in
\autoref{sec:2d-complex-painting}.
For $0 \leq z \leq q-1$, we have that $c_{z,0} = 0$ if and only if at least one of the following binary variables has value $0$.
\begin{subequations}
\label{eq:directly_covered_constraints}
\begin{gather}
\ell_{z,y} \; (0 \leq y \leq q-1), \quad u_{z,y} \; (0 \leq y \leq q-1), \\
\ell_{x,z} \; (0 \leq x \leq q-1),  \quad u_{x,z} \; (0 \leq x \leq q-1),\\ 
\ell_{x,y} \; (0 \leq x, y \leq q-1 \text{ such that } x+y \equiv z \pmod q), \\
u_{x,y} \;  (0 \leq x, y \leq q-1 \text{ such that } x+y+1 \equiv z \pmod q)
\end{gather}
\end{subequations}

We assume that the saturation process of
transferring indirectly covered properties described in \autoref{sec:2d-complex-painting}
ends in a finite number
$\sage{maxstep}$ of steps. (Though no theoretical bound better than
$\sage{maxstep} \leq q$ is known, in practice a small value of \sage{maxstep}
such as 2 is sufficient.)
For $1 \leq i \leq
\sage{maxstep}$, we define binary variables $c_{z,i}$ ($0 \leq z \leq q-1$) to
model whether the interval $[\frac{z}{q}, \frac{z+1}{q}]$ is covered in the
first $i$ steps of the saturation process. If the interval $[\frac{z}{q}, \frac{z+1}{q}]$ was already
covered in the step $i-1$, then it remains covered in the step $i$. Thus, we have the constraint
\begin{subequations}
\begin{equation}
c_{z,i-1}=0 \Rightarrow c_{z,i}=0.
\end{equation}
If the interval $[\frac{z}{q}, \frac{z+1}{q}]$ is connected by a green edge to an interval $[\frac{x}{q}, \frac{x+1}{q}]$ that was already covered in the step $i-1$, then the interval $[\frac{z}{q}, \frac{z+1}{q}]$ becomes covered in the step $i$.
\begin{equation}
 \min\{d_{x,z}, h_{x, (z-x) \bmod  q}, v_{(x-z)\bmod q, z}\} = 0 \text{ and } c_{x,i-1}=0
 \Rightarrow c_{z,i}=0.
\end{equation}
\label{eq:undirectly_covered_i_constraints}
\end{subequations}
Otherwise $c_{z,i}=1$.

Note that these constraints can all be expressed using linear equations or inequalities
over binary variables. 

\subsection{Trade-off between strictness and basicness}
\label{sec:epsilon-discussion}

Assume that the variables and the constraints defined in \autoref{sec:slope-assignment} were not introduced to our MIP, and that the strict subadditivity constraints were not enforced (i.e., set $\epsilon=0$ in \eqref{eq:strict-subadditivity}). 
\begin{remark}
\label{rk:smaller-polytope}
Painting faces dark green in $\Delta\P_{\frac1q\Z}$  (i.e., setting binary variables $p_{x,y}$ to $0$) amounts to restricting the corresponding subadditivity inequalities \eqref{itm:subadditive} of \autoref{thm:finite-minimal} to equations in the constraint system of the polytope $\Pi_f(\tfrac1q\Z/\Z)$ of minimal functions. 
Thus, the set of restricted functions $\pi|_{\frac{1}{q}\Z}$ that satisfy the new constraint system is a smaller polytope, which is a face of the polytope $\Pi_f(\frac{1}{q} \Z/ \Z)$. A vertex of the smaller polytope is also a vertex of the polytope $\Pi_f(\frac{1}{q} \Z/ \Z)$.
\end{remark}
After fixing $0/1$ variables, a basic feasible solution $\pi$ of the system
with $c_{z, \sage{maxstep}} = 0$ for $0\leq z \leq q-1$ is extreme for $R_f(\R/\Z)$,
according to \autoref{thm:extreme-finite-covered}. 
Unfortunately, in this case we cannot expect the solutions returned by the MIP solver to always have $k$ different slope values; indeed, they often degenerate to 2-slope or 3-slope functions. 
The same phenomenon was observed in the shooting
experiments in the literature
\cite{gomory2003corner,evans-thesis,shim2009large}, where the extreme
functions that received a large percentage of 
hits are 2-slope or 3-slope functions.

Because of this, we use the slope value variables and constraints from
\autoref{sec:slope-assignment} in the MIP, always enforcing the strict inequality in slope variables by making a practical choice of $\epsilon'>0$ for \eqref{eq:slope-difference}. In this way we can set the number of slopes $k$ of the resulting functions explicitly. We also use $\epsilon>0$ in \eqref{eq:strict-subadditivity}, which allows for a significant speedup.

\begin{remark}
\label{rk:choice-of-epsilon}
\autoref{cor:arithmetic-complexity-upper-bound} shows that we do not lose any
extreme functions by setting $0<\epsilon, \epsilon' \leq 10^{-{q}/{4}}$. 
Our code instead uses the heuristic choice $\epsilon =\epsilon'= 1/4$ (or $1/12$), which allows for
stronger pruning based on linear programming 
at the expense of losing some functions.  
\end{remark}

However, now a basic feasible solution of the system after fixing $0/1$ variables is no
longer guaranteed to be an extreme function, since the smaller polytope with the inequalities $\pi_x + \pi_y - \pi_{(x + y) \bmod q} \geq \epsilon$ is not a face of the polytope $\Pi_f(\tfrac1q\Z/\Z)$. 
 In this case, a further extremality test needs to be applied to the returned
 candidate function~$\pi$. Since all intervals are covered (i.e.,  $c_{z,
   \sage{maxstep}} = 0$ for $0\leq z \leq q-1$), it suffices by
 \autoref{thm:extreme-finite-covered} to test whether $\pi|_{\frac{1}{q}\Z}$ is a vertex of $\Pi_f(\frac{1}{q}\Z/\Z)$. 

\subsection{Objective function}

A tailored objective function can be used to steer the optimum away from equality of slopes 
and from the lower bounds of the subadditivity slacks $\Delta\pi_{x,y}$, 
ensuring that the basic optimal
solution returned by the MIP solver will
correspond to an extreme function. However, there is no a priori best choice
of such an objective function that will guarantee success.  
In our computations, we always maximize the
difference of slopes $s_1 - s_k$.  Other objective functions, including the
following, are plausible:
\begin{itemize}
\item maximize a weighted difference of slopes $\displaystyle\sum_{j=1}^{k-1}{\lambda_j (s_j - s_{j+1})}$,
  for some suitable weights~$\lambda_j$;
\item maximize a weighted sum of subadditivity slacks $\displaystyle\sum_{0 \leq x \leq y \leq q-1}{\omega_{x, y}\Delta\pi_{x,y}}$;
\item maximize a weighted sum of subadditivity points $\displaystyle\sum_{0 \leq x \leq y \leq q-1}{\omega_{x, y}p_{x,y}}$; 
\item minimize the covering count $\displaystyle\sum_{0 \leq x, y \leq q-1}{u_{x,y}+\ell_{x,y}}$.
\end{itemize}
\medbreak

\subsection{Implementation and discussion}

The MIP approach is easy to implement and also easy to tailor to a search for
extreme functions with particular properties.  

However the approach is limited because floating-point implementations are not
a good match for finding functions of high arithmetic complexity.  If $q$ is
large, the difference between two slope values $s_i$ often becomes extremely
small and may completely disappear in floating point fuzz, making the choice
of the parameter~$\epsilon'$ difficult. 

Moreover, MIP solvers are generally not the best tool for performing an
exhaustive search.  While listing several solutions should be
possible by varying the objective function, this is rather difficult to do in
practice.  We have used the solver Gurobi (version 5.6.3) to solve the MIP problem. We
resorted to setting the Gurobi parameter \sage{SolutionLimit=1} and calling
\sage{optimize()} repeatedly, so that the feasible solutions found by Gurobi
before reaching the global optimal solution are recorded.

Despite these limitations, we have obtained new results
with the MIP approach, which we report in the next two subsections.

\subsection{Result: Optimality of the oversampling factor $3$}

Using the MIP approach described above, we found an example to answer an open
question in \cite{igp_survey_arxiv_v1}.  The open question relates to the
following theorem.
\begin{theorem}[{\cite[Theorem 8.6]{igp_survey}}]
\label{thm:extreme-restriction-m}
  Let $m\geq 3$, the \emph{oversampling factor}, be a positive integer.  Let $\pi$ be a continuous piecewise linear  minimal
  valid function for $R_f(\R/\Z)$ with breakpoints in $\tfrac{1}{q} \Z$ and suppose $f \in \tfrac{1}{q} \Z$.    The following are equivalent:
\begin{enumerate}
\item $\pi$ is a facet for $R_f(\R/\Z)$,
\item $\pi$ is extreme for $R_f(\R/\Z)$,
\item $\pi|_{\frac{1}{mq}\Z}$ is extreme for $R_f(\tfrac{1}{mq} \Z/\Z)$.
\end{enumerate}
\end{theorem}
Dey et al.~\cite{dey1} gave a function $\pi = \sage{\sagefunc{drlm_not_extreme_1}()}$ that is not 
extreme for $R_f(\R/\Z)$, but $\pi|_{\frac1q\Z}$ is extreme $R_f(\tfrac{1}{q}
\Z/\Z)$. Thus, \autoref{thm:extreme-restriction-m} does not hold with the
oversampling factor $m=1$.  Basu et
al.~\cite{basu-hildebrand-koeppe:equivariant} proved
\autoref{thm:extreme-restriction-m}
for an oversampling factor $m=4$.  It was strengthened to any $m\geq3$
in~\cite{igp_survey}.
The question whether $m\geq3$ is best possible or can be improved to $m=2$ was
stated in \cite[Open Question 8.7]{igp_survey_arxiv_v1}.
We resolve this question by the following result, which was stated in the
introduction as \autoref{prop:oversampling3optimal}.
\begin{prop}
\label{prop:oversampling3optimal-repeat}
  The lower bound $m\geq 3$ for the oversampling factor in
  Theorem~\ref{thm:extreme-restriction-m} is best possible.
  Theorem~\ref{thm:extreme-restriction-m} does not hold when $m=2$.
\end{prop}
Unable to construct an example function by hand
, we established this result by using the MIP approach.  
Before we explain the details of the search strategy, we describe the
structure of the example function that we found, \sagefunc{kzh_2q_example_1}. 
It is a continuous $4$-slope function with $q=37$ and $f=\tfrac{25}{37}$; see
\autoref{fig:kzh_2q_move}.  One can verify, simply using the automated
extremality test implemented in the software
\cite{infinite-group-relaxation-code}, that \sagefunc{kzh_2q_example_1} is a
non-extreme function\footnote{As proved by \sage{extremality\underscore
    test(kzh\underscore 2q\underscore example\underscore 1())} returning
  \sage{False}.}, 
whose restriction to $\frac{1}{2q}\Z$ is extreme\footnote{As proved by
  \sage{simple\underscore finite\underscore dimensional\underscore
    extremality\underscore test(kzh\underscore 2q\underscore
    example\underscore 1(), oversampling=2)} returning \sage{True}.}. This
proves that an oversampling factor of 3 is optimal. 

In the following, we provide a brief justification of the extremality result
given by the code.
It connects to the theory of equivariant perturbations developed
in~\cite{basu-hildebrand-koeppe:equivariant}.

\begin{proof}[Proof of \autoref{prop:oversampling3optimal-repeat}]
We show that the function $\pi$ is not extreme for $R_f(\R/\Z)$, by computing
the covered intervals. It can be verified\footnote{One can type
  \sage{plot\underscore 2d\underscore diagram(kzh\underscore 2q\underscore
    example\underscore 1(),colorful=True)} to visualize the painting on the
  complex $\Delta\P_{\frac1q\Z}$.} 
on the complex $\Delta\P_{\frac1q\Z}$ that
\begin{enumerate}[(i)]
\item the intervals $[\tfrac{10}{37}, \tfrac{11}{37}]$ and
  $[\tfrac{14}{37}, \tfrac{15}{37}]$, indicated by yellow strips in
  \autoref{fig:kzh_2q_move}, are uncovered; 
\item the others intervals are covered.
\end{enumerate}
There exists thus a non-zero perturbation function $\bar{\pi}$, such that
$\pi=\tfrac{1}{2}\pi^1+\tfrac{1}{2}\pi^2$ where $\pi^1 =
\pi+\bar{\pi}$, $\pi^2=\pi-\bar{\pi}$ are two distinct minimal functions, showing
the non-extremality of $\pi$ for $R_f(\R/\Z)$. The function plotted in magenta
in \autoref{fig:kzh_2q_move} is such a function~$\bar\pi$.

By \cite[Definition 3.3, Lemma 2.7 and Theorem 4.6]{basu-hildebrand-koeppe:equivariant}, a perturbation function $\bar{\pi}$ is affine linear on the covered intervals, and satisfies $\bar{\pi}(0)=\bar{\pi}(f)=\bar{\pi}(1)=0$, $\Delta\bar{\pi}(x, y)=0$ for any $(x,y)$ such that $\Delta\pi(x,y)=0$. 
Consider the restriction to $\frac{1}{q}\Z$. One can show by linear algebra\footnote{This could also be verified by \sage{simple\underscore finite\underscore dimensional\underscore extremality\underscore test(kzh\underscore 2q\underscore example\underscore 1(), oversampling=1)} which returns \sage{True}.} that the finite dimensional linear system has a unique solution $\bar{\pi}|_{\frac{1}{q}\Z} = 0$. It follows that $\bar{\pi}$ can only be non-zero on the uncovered intervals, i.e., on $[\tfrac{10}{37}, \tfrac{11}{37}]$ and $[\tfrac{14}{37}, \tfrac{15}{37}]$. Furthermore, $\pi|_{\frac{1}{q}\Z}$ is extreme for $R_f(\frac{1}{q}\Z/\Z)$. 

Next, we will show that $\pi|_{\frac{1}{2q}\Z}$ is extreme for $R_f(\frac{1}{2q}\Z/\Z)$.  
Again it can be verified that
\begin{enumerate}[(i)]\setcounter{enumi}{2}
\item $F([\tfrac{10}{37},\tfrac{11}{37}], \{\frac{4}{37}\}, [\tfrac{14}{37},
  \tfrac{15}{37}])$ is an additive edge of~$\Delta\P_{\frac1q\Z}$.
\end{enumerate}
In other words, for $x \in [\tfrac{10}{37},\tfrac{11}{37}]$, we have
$\Delta\pi(x, \frac{4}{37})=0$, and thus $\Delta\bar{\pi}(x, \frac{4}{37})=0$.
Following 
\cite{basu-hildebrand-koeppe:equivariant}, 
for $t \in \R$, define the translation $\tau_t\colon \R \to \R$, $x \mapsto x+t$. The gold-colored single headed arrow in \autoref{fig:kzh_2q_move} indicates the action of translation $\tau_{4/37}$, which sends the interval $[\tfrac{10}{37}, \tfrac{11}{37}]$ to the interval $[\tfrac{14}{37}, \tfrac{15}{37}]$.
Therefore, $\bar{\pi}(\tau_{4/37}(x)) = \bar{\pi}(x) + \bar{\pi}(\tfrac{4}{37})=\bar{\pi}(x)$ for $x \in [\tfrac{10}{37},\tfrac{11}{37}]$, as $\bar{\pi}(\tfrac{4}{37})=0$. In particular, we have $\bar{\pi}(\tfrac{21}{74}) = \bar{\pi}(\tfrac{29}{74})$.
For 
$r\in\R$, define the reflection $\rho_r\colon \R \to \R$, $x \mapsto r-x$. 
The gold-colored double-headed arrow in \autoref{fig:kzh_2q_move} indicates the action of reflection $\rho_{f}$ between the two uncovered intervals, corresponding to the symmetry condition $\pi(x)+\pi(f-x)=\pi(f)=1$ for $x \in [\tfrac{10}{37},\tfrac{11}{37}]$. We have $\bar{\pi}(x) + \bar{\pi}(\rho_{f}(x))= \bar{\pi}(f) = 0$ for $x \in [\tfrac{10}{37},\tfrac{11}{37}]$. In particular, $\bar{\pi}(\tfrac{21}{74}) + \bar{\pi}(\tfrac{29}{74}) = 0$. Therefore, $\bar{\pi}(\tfrac{21}{74}) = \bar{\pi}(\tfrac{29}{74})=0$, the perturbation function $\bar{\pi}$ has values zero at the midpoints of the two uncovered intervals. Since $\bar{\pi}=0$ on the other intervals which are covered, we have $\bar{\pi}|_{\frac{1}{2q}\Z} = 0$. Hence, $\pi|_{\frac{1}{2q}\Z}$ is extreme for $\smash[t]{R_f(\frac{1}{2q}\Z/\Z)}$.
\end{proof}

Our tailored MIP search strategy was to look for functions with properties
(i--iii) from the above proof.  These properties correspond to a prescribed
partial painting and thus can be expressed by fixing some binary variables.
We tried out various pairs of $(q, f)$ for $10
\leq q \leq 40$.  
For $q = 37$ and $f=\frac{25}{37}$, we discovered the function \sagefunc{kzh_2q_example_1}.\footnote{The example can be reproduced using the code in \cite{infinite-group-relaxation-code} as follows. First, call the function \sage{write\underscore lpfile\underscore 2q(q=37, f=25/37, a=11/37, kslopes=4, maxstep=2, m=4)} to generate a MIP problem that maximizes the slope difference $s_4-s_1$. The parameter \sage{a} indicates that the uncovered intervals are $[a-\frac1q, a]$ and $[f-a,f-a+\frac1q]$, and the parameter \sage{m} decides the heuristic choice of $\epsilon = \epsilon' = 1/m$. The MIP problem is written to the file named \sage{mip\underscore q37\underscore f25\underscore a11\underscore 4slope\underscore 2maxstep\underscore m4.lp}. Then, use Gurobi to solve the MIP problem and write the solution to the file named \sage{solution\underscore 2q\underscore example\underscore m4.sol}. Finally, retrieve the function \sagefunc{kzh_2q_example_1} from the solution file by calling \sage{refind\underscore function\underscore from\underscore lpsolution\underscore 2q('solution\underscore 2q\underscore example\underscore m4.sol', q, f, a)}.}

\subsection{Result: Refutation of the generic 4-slope conjecture}

Our search also resolves {\cite[Open Question
  2.16]{igp_survey_arxiv_v1}} by showing that even for functions whose
extremality proof only uses the Interval Lemma, rather than the more general
techniques from \cite{basu-hildebrand-koeppe:equivariant} (translations and
reflections), many slopes are possible.  This is in contrast to the first
5-slope functions\footnote{The functions are available in the electronic
  compendium \cite{electronic-compendium} as
  \sage{hildebrand\underscore{}5\underscore{}slope...}} found by Hildebrand
(2013, unpublished), whose extremality proof requires translating and
reflecting covered intervals.

\begin{figure}[t]
\centering
\includegraphics[height=5cm]{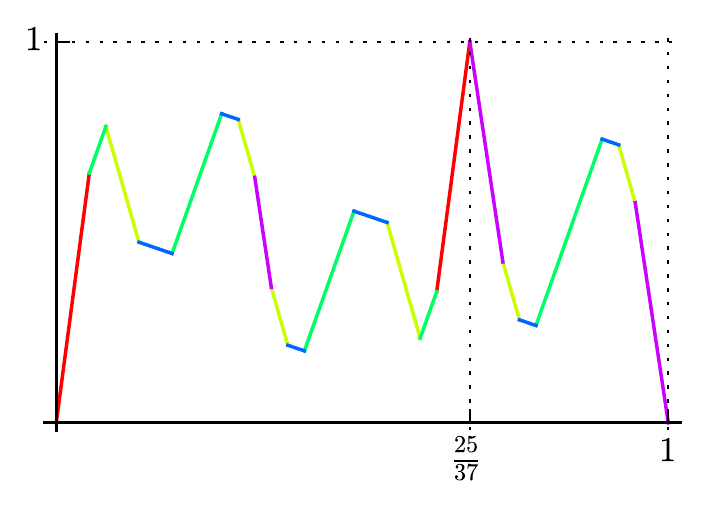}
\quad
\includegraphics[height=5cm]{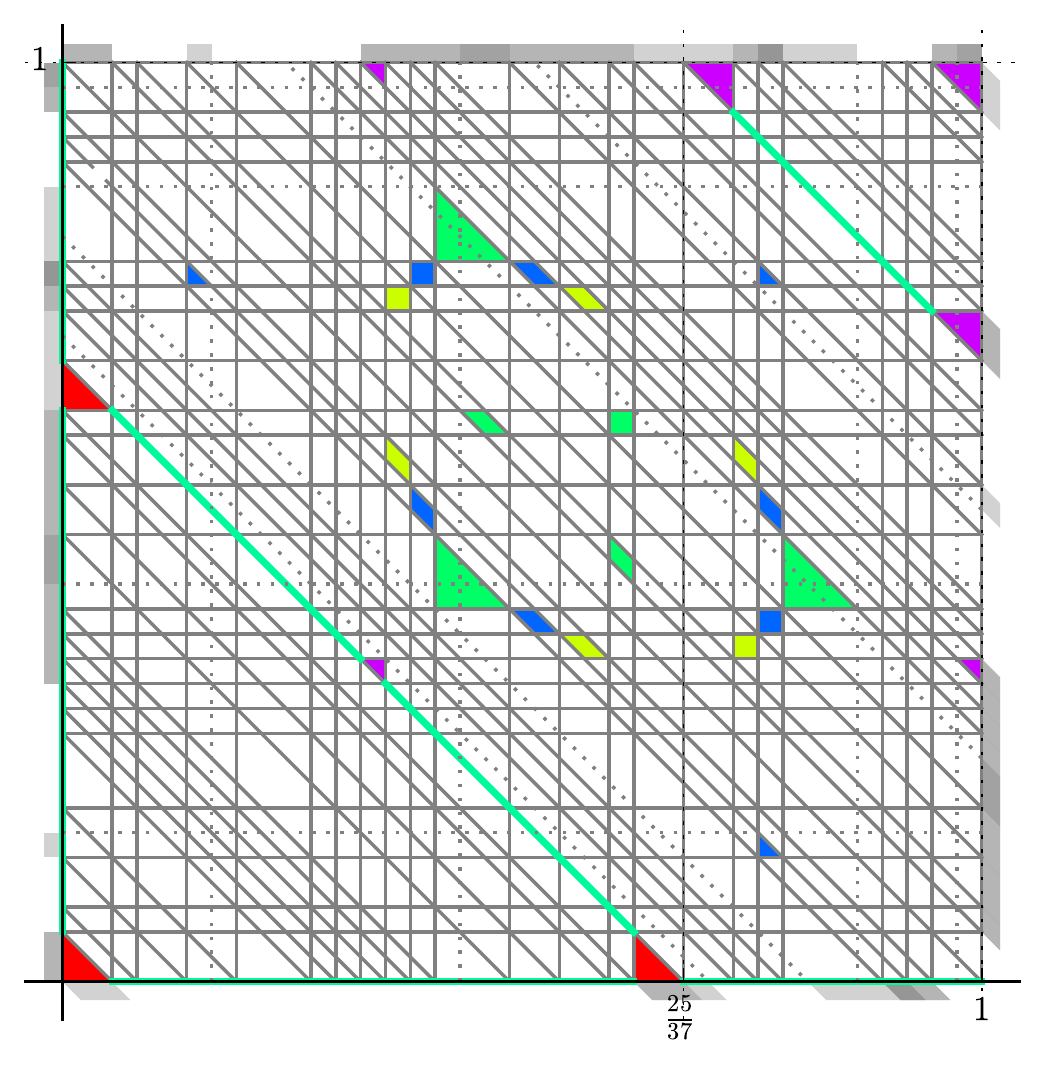}
\caption{The 5-slope extreme function
  \sagefunc{kzh_5_slope_fulldim_1}
  found by our search code (\textit{left}). Its two-dimensional polyhedral
  complex $\Delta\P$ (\textit{right}), as plotted by the command
  \sage{\sagefunc{plot_2d_diagram}(h,colorful=True)}, does not have any
  lower-dimensional maximal additive faces except for the symmetry reflection
  or $x=0$ or $y=0$.} 
\label{fig:kzh_5_slope_fulldim}
\end{figure}

\begin{prop}
  There exists a piecewise linear extreme function $\pi$ of $R_f(\R/\Z)$ with
  more than 4 slopes, such that its additivity domain $E(\pi) := \setcond{(x,
    y)} { \Delta \pi(x, y) = 0}$ is the union of full-dimensional convex sets
  and the lines $x \in \Z$, $y \in \Z$, $x+y \in f+\Z$. 
\end{prop}
See the functions \sagefunc{kzh_5_slope_fulldim_1}
etc.\@, which we have made available as part of \cite{electronic-compendium}.\footnote{These examples can be reproduced using the code in \cite{infinite-group-relaxation-code} as follows. First, call \sage{write\underscore lpfile(q, f, kslopes, m=12, type\underscore cover='fulldim')} with appropriate values of \sage{q, f} and \sage{kslopes}  to generate a MIP problem that maximizes the slope difference. For example, we set \sage{q=37; f=25/37; kslopes=5}. The MIP problem is written to the file named \sage{5slope\underscore q37\underscore f25\underscore fulldim\underscore m12.lp}. Then, use Gurobi to solve the MIP problem. We set the Gurobi parameter \sage{SolutionLimit=1} and call \sage{optimize()} repeatedly, so that the feasible solutions found by Gurobi before reaching the global optimal solution are recorded to the files named \sage{solution\underscore 5slope\underscore fulldim\underscore 1.sol}, etc. Finally, retrieve the functions \sagefunc{kzh_5_slope_fulldim_1}, etc.\@\ from the solution files by calling \sage{refind\underscore function\underscore from\underscore lpsolution('solution\underscore 5slope\underscore fulldim\underscore 1.sol', q, f).}} They are continuous 5-slope extreme functions without any $0$-dimensional
or $1$-dimensional maximal additive faces except for the symmetry reflections
$x+y \in f + \Z$
and the trivial additivities $x\in \Z$, $y\in \Z$.  A graph of the function
\sagefunc{kzh_5_slope_fulldim_1} and a plot of its painting on the
two-dimensional complex are shown in \autoref{fig:kzh_5_slope_fulldim}.

\begin{remark}
  Using our computer-based search we also found extreme functions that do have
  lower-dimensional additive faces, but whose extremality proof does not
  \emph{depend} on those. All intervals are directly covered.
  Examples are provided by 
  the functions \sagefunc{kzh_5_slope_fulldim_covers_1},
  \sagefunc{kzh_6_slope_fulldim_covers_1} etc.\@, which we have made available
  as part of \cite{electronic-compendium}.\footnote{These
    \sage{fulldim\underscore covers} examples can be reproduced in the same
    way as for the \sage{fulldim} examples described in the last footnote,
    except that \sage{type\underscore cover} is set to
    \sage{'fulldim\underscore covers'} when generating the MIP problems. For
    example, the function \sagefunc{kzh_6_slope_fulldim_covers_1} is obtained
    from the MIP problem \sage{6slope\underscore q25\underscore f8\underscore
      fulldim\underscore covers\underscore m12.lp} generated by
    \sage{write\underscore lpfile(q=25, f=8/25, kslopes=6, m=12, type\underscore
      cover='fulldim\underscore covers')}.}  
\end{remark}

\section{Backtracking search}
\label{sec:backtracking-search}
\subsection{Search via covering paintings}
In this section we discuss a new search strategy that addresses the
limitations of the MIP approach of the previous section by using our own
implementation of backtracking tailored to enumerating covering paintings.

\begin{figure}[p]
  \centering
  \includegraphics[width=\linewidth]{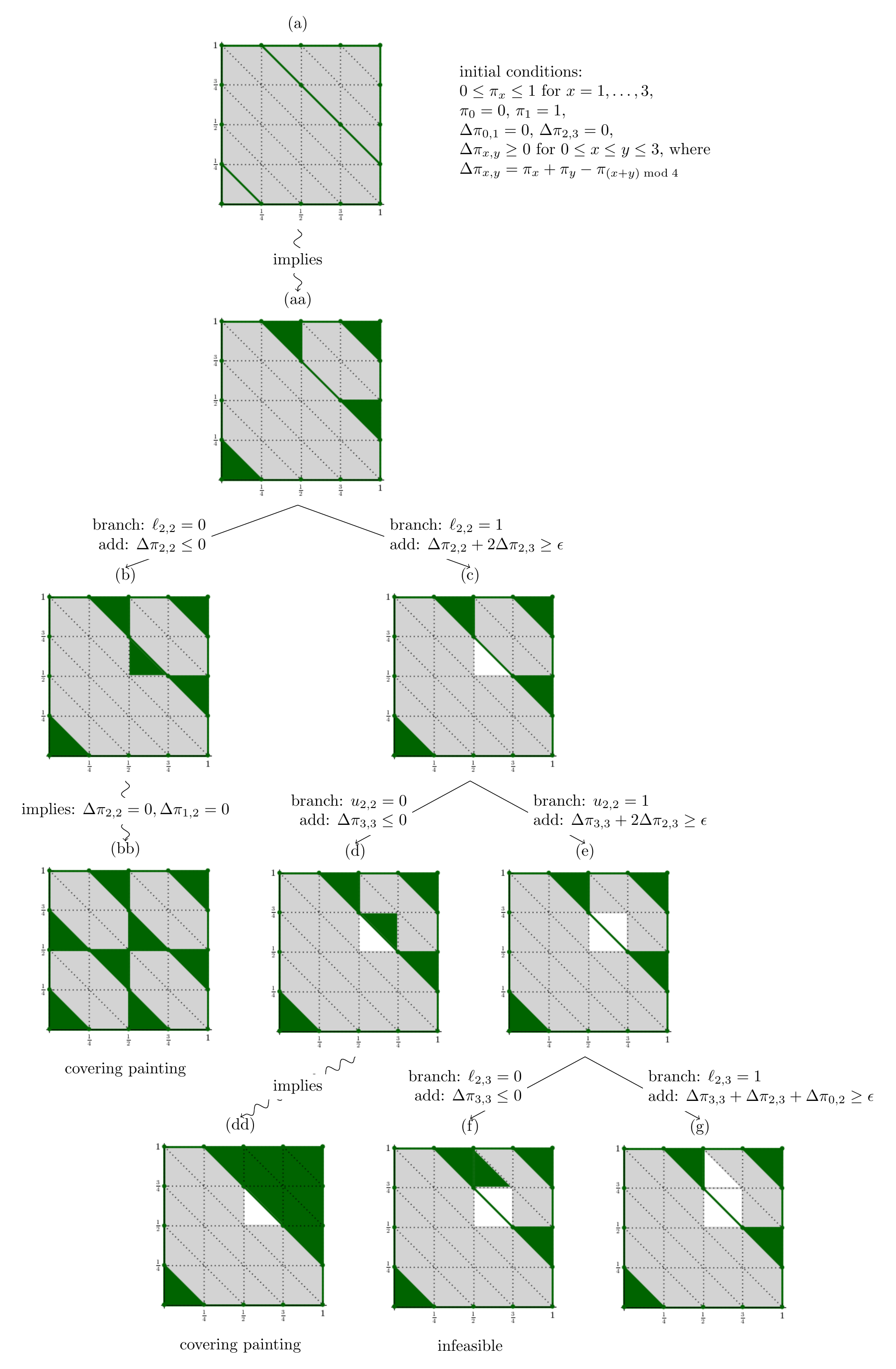}
  \caption{A tree of partial paintings of $\Delta\P_{\frac14\Z}$ with $f=\frac14$.} 
  \label{fig:partial-painting-tree}
\end{figure}

During our backtracking search, we maintain a \emph{prescribed partial painting} as
introduced in \autoref{ppp}.  At the root, the minimality conditions
(\autoref{thm:finite-minimal}) 
\begin{subequations}
\begin{align}
&\pi_0 = 0 \text{ and } 0 \leq \pi_x \leq 1, &\text{ for } x=1, \dots, q-1 \\
&\Delta\pi_{x,y} = \pi_x+\pi_y - \pi_{(x+y)\bmod q} \geq 0, &\text{ for } 0\leq x \leq y \leq q-1 \\
&\pi_{qf} = 1 \text{ and }  \Delta\pi_{x, (qf - x) \bmod  q} = 0, &\text{ for } x=0, \dots, q-1
\end{align}%
\label{eq:init-lp-constraints}%
\end{subequations}%
give the initial prescribed partial painting
on $\Delta\P_{\frac1q\Z}$. See \autoref{fig:partial-painting-tree} for an example where $q=4$ and $f=\frac{1}{4}$. To get an
extreme function, more additivity relations are needed. To achieve this, in
our backtracking search we successively paint some
light grey faces white or dark green, 
until a covering painting is reached. During the painting process, we keep track of the consistency of the colors of the faces 
by using an LP that will be explained later in \autoref{sec:warm-start-lp}. If the current partial painting is infeasible, pruning will be performed.

\autoref{thm:extreme-finite-covered} and \autoref{prop:minimal-polytope} have
the following corollary.
\begin{corollary}
Let $\pi$ be a continuous piecewise linear function with breakpoints in $\tfrac{1}{q} \Z$. If $\pi$ is extreme for $R_f(\R/\Z)$, then there exists a covering painting such that $\pi|_{\frac{1}{q}\Z}$ is a vertex of the polytope 
formed by the minimal functions for $R_f(\frac{1}{q} \Z/ \Z)$
whose additivities correspond to the painting. 
\end{corollary}

The search for extreme functions $\pi$ is thus converted into the search for covering paintings. 

\subsection{Branching rule}
\label{sec:branching-rule}
The search tree has a binary structure. Assume that we are branching on a
triangle $F$ (say the lower triangle whose lower left corner is vertex
$(x,y)$) that is colored (light) grey in the partial painting of the current node. In
terms of \autoref{sec:mip_approach}, the value of $\ell_{x,y}$ is undecided. 
Branching on the node creates two children: one where the triangle $F$ is
painted dark green ($\ell_{x,y}$ is set to $0$) and one where (the interior
of) $F$ is painted white ($\ell_{x,y}$ is set to $1$). 

In the child node with dark green $F$, the following additivity constraints hold. 
\begin{equation}
\Delta\pi_{u,v}=0, \quad \text{ for every vertex }(u, v)\text{ of }F. 
\label{eq:additivity-constraints}
\end{equation}

In the child node with white $F$, the following strict subadditivity relation holds: 
\begin{equation*}
\sum_{\text{vertex } (u, v) \text{ of } F} \Delta\pi_{u,v}>0.
\end{equation*}
Since the strict inequality constraints are not allowed in linear programming, 
we prefer to 
%
%
replace it by 
\begin{equation}
\sum_{\text{vertex } (u, v) \text{ of } F} \Delta\pi_{u,v} \geq \epsilon,
\label{eq:white-face-epsilon}
\end{equation}
where $\epsilon$ is a small positive number. See \autoref{rk:choice-of-epsilon} for the value of $\epsilon$.
%

\subsection{Feasibility and satisfiability checks}
\label{sec:warm-start-lp}
To check the feasibility of a prescribed partial painting and detect its implied additivity
relations, we use linear programming. 
We have investigated two options.  As we mentioned before, SageMath has a good
interface to the Parma Polyhedra Library.  Its implementation of the double
description method provides an attractive interface to checking feasibility and for
testing implied additivities, all in exact arithmetic.  This approach, however,
appears to be quite slow when the dimension of the polytope is large. As a rule of thumb, when the dimension exceeds $9$, it is better to apply the simplex method instead.
Our search code uses the GLPK solver, which is integrated well in SageMath (see
\autoref{sec:lp-ppl-glpk} for details), and allows for warm-starting the simplex method.

We construct the linear optimization problem as follows. The problem has $q$ real variables $\pi_0, \dots, \pi_{q-1}$, and some auxiliary variables $\Delta\pi_{x,y}$ ($0 \leq x \leq y \leq q-1$) that represent the subadditivity slacks. The initial constraints on the variables are 
given in \eqref{eq:init-lp-constraints}.
The subadditivity and additivity constraints are reflected by the bounds of
their slack variables.\footnote{We need to introduce these slack variables
  explicitly due to limitations of warm-starting in the SageMath interface.} 
These bounds will vary along the backtracking process. 
If we walk downwards in the tree, then we append either \eqref{eq:additivity-constraints} or \eqref{eq:white-face-epsilon} to the constraint system; conversely, if we walk upwards in the tree, then we remove the constraint. 

Such changes in the constraint system could affect the feasibility of the problem. Due to the update of the variable bounds, the dual simplex method starting off the last basis is called to check whether the partial painting remains feasible at the current node. If the node is infeasible, the whole sub-tree will be pruned.
For example, the node (f) in \autoref{fig:partial-painting-tree} is infeasible, as $\Delta\pi_{3,3}\leq 0$ and $\Delta\pi_{3,3}\geq \epsilon >0$.

The current constraint system may imply some new additive faces in the
prescribed partial painting. To check if a vertex $(u, v)$ is implied additive, we call the primal simplex method starting off the last basis  
to maximize the objective function $\Delta\pi_{u,v}$. If the optimal value is
$0$, then $\Delta\pi_{u,v} = 0$ and thus the vertex $(u,v)$ is implied
additive (dark green). In terms of \autoref{sec:mip_approach}, $p_{u,v}=0$. Subject to the inclusion constraints \eqref{eq:logical_constraints}, the colors of edges and triangles are updated accordingly.
See the nodes (bb) and (dd) in \autoref{fig:partial-painting-tree} for examples.

\subsection{Heuristic choice of the branching triangle}

By the invariance of the subadditivity condition under exchanging $x$ and $y$,
only the upper left triangular part of the complex $\Delta\P_{\frac1q\Z}$ needs
to be considered for painting. 

In our experiments we found that an exhaustive search of all paintings of the
two-dimensional complex is too expensive even for a moderate size of~$q$.  
Thus, to reach a covering painting quickly, we branch on the colors of the
triangles ($\ell_{x,y}$ and $u_{x,y}$) rather than on the colors of the
vertices ($p_{x,y})$. Furthermore, a heuristic painting strategy is applied:
While picking a light grey triangle to paint dark green or white in branching, we consider
those triangles $F$ whose projections $p_1(F)$ and $p_2(F)$ are currently
uncovered. 
This is of course restrictive, and so we cannot guarantee that our search code will find
all covering paintings and all extreme functions.  However, it has proved to
be a successful heuristic strategy.

At each node, we choose one candidate triangle $F$ (i.e., one $\ell_{x,y}$ or $u_{x,y}$ variable) among all these considered
triangles. It is defined as the smallest triangle in lexicographical order,
whose color has not been branched on yet in the search tree. 

\subsection{Incremental computation}
\label{sec:incremental_computation}

For the purpose of improving the running time efficiency, all computations in the backtracking search, such as updating covered intervals, are done in an incremental manner. 

More precisely, we maintain a list of \textit{connected components} that we have mentioned in \autoref{sec:2d-complex-painting}. When a
new triangle $F$ is painted as additive in the prescribed partial painting,
the components that contain the projection $p_1(F)$ or $p_2(F)$ or $p_3(F)$
are merged into one big component, and all intervals in this new component
become covered. When a new edge (one-dimensional face) $F$ is painted as
additive (dark green), the components that contain its projection intervals are
merged into one big component. If the new component contains an interval that
was covered, then all intervals in this new component are covered. In such a
way, the new covered intervals after adding an additive face can be computed
incrementally from the covered intervals in the previous step. 
For example, the node (a) in \autoref{fig:partial-painting-tree} has three connected components $\mathcal{C}_1 =\{[0,\frac{1}{4}]\}$, $\mathcal{C}_2 =\{[\frac{1}{4}, \frac{1}{2}], [\frac{3}{4},1]\}$ and $\mathcal{U} =\{[\frac{1}{2},\frac{3}{4}]\}$, where the first two are covered and the last one is not. When $\ell_{2,2}$ is set to $0$ in its child (b), $\mathcal{U}$ is merged into $\mathcal{C}_1$, and becomes covered. Thus, the node (b) has two connected components $\mathcal{C}_1' =\{[0,\frac{1}{4}], [\frac{1}{2},\frac{3}{4}]\}$ and $\mathcal{C}_2 =\{[\frac{1}{4}, \frac{1}{2}], [\frac{3}{4},1]\}$ that are both covered.


We mentioned above that a function $\pi$ whose additivities satisfy the prescribed partial painting
has the same slopes on the intervals in each component, see \cite[Remark
3.6]{basu-hildebrand-koeppe:equivariant}. By counting the connected
components, we get an upper bound on the number of slopes that the function
$\pi$ could have.  This allows us to prune subtrees that cannot contain
functions with the desired number of slopes.

The knowledge of connected components is also used at the end of ``vertex
filtering mode'' (\autoref{sec:vertex-filtering-search}), to check efficiently whether all intervals are covered, as follows. If an interval $[\frac{z}{q},\frac{z+1}{q}]$ does not belong to any covered connected component from the partial painting and if none of the values \eqref{eq:directly_covered_constraints} is $0$, then $\pi$ has uncovered interval and hence is not extreme. This saves us from running the full extremality test on $\pi$,  which would consume more time according to \autoref{fig:vertex-filtering-performance}--middle.


\subsection{Heuristic search algorithm}
\label{sec:heuristic_mode}
We summarize the backtracking search via covering paintings in \hyperref[alg:heuristic-mode]{Algorithm~\ref*{alg:heuristic-mode}}. It is referred to as the ``heuristic mode'' search in our code.
\begin{algorithm}
\caption{heuristic mode}
\label{alg:heuristic-mode}
\begin{enumerate}
\item The root node is the initial painting given by the minimality conditions \eqref{eq:init-lp-constraints};
\item Decide for a candidate triangle $F$ of the painting using covered intervals;
\item A node is branched into two child nodes:
\item \label{itm:green-subnode} For the child node in which $F$ is additive
  (dark green), 
	\begin{itemize}
	\item add the new additivity relations \eqref{eq:additivity-constraints} to constraints;
	\item look for implied additive vertices and triangles;
	\item update covered intervals;
	\item if the node is infeasible, backtrack;
	\item if the number of the connected components is less than the desired number of slopes, backtrack;
	\item if a covering painting is found, \textbf{output} it and backtrack;
	\end{itemize}
\item For the sub-node in which $F$ is strictly subadditive (white), 
	\begin{itemize}
	\item add the strict subadditivity relation \eqref{eq:white-face-epsilon} to constraints;
	\end{itemize}
\item Traverse the search tree in depth-first order.
\end{enumerate}
\end{algorithm}

For each covering painting returned by the algorithm, we use the minimality conditions \eqref{eq:init-lp-constraints} and the additivity relations $\Delta\pi_{x,y}=0$ specified by the green vertices $(x,y)$ in the painting, to construct a polytope. (The strict subadditivity constraints $\Delta\pi_{x,y}>0$ corresponding to the green vertices $(x,y)$ in the painting are neglected.)  We enumerate the vertices of the polytope. 
By~\autoref{rk:smaller-polytope} and \autoref{thm:extreme-finite-covered},
the interpolation $\pi$ of a vertex $\pi|_{\frac{1}{q}\Z}$ is extreme for $R_f(\R/\Z)$. 

\subsection{Combined search algorithm}
\label{sec:combined-mode}
The above search algorithms work well for relatively small $q$, but become
inefficient when $q$ is large: 
The vertex filtering search (\hyperref[alg:vertex-filtering-mode]{Algorithm~\ref*{alg:vertex-filtering-mode}})  wastes time on enumerating numerous vertex-functions in high dimension, most of which are non-extreme for the infinite group problem; the heuristic backtracking search via covering paintings (\hyperref[alg:heuristic-mode]{Algorithm~\ref*{alg:heuristic-mode}})
suffers from the combinatorial explosion in branching and the general
performance penalty from using Python.

We propose to combine the vertex filtering search and the heuristic backtracking search to obtain a better performance. 
The combined algorithm starts with branching, but outputs the partial painting and
backtracks at a certain depth before reaching a covering painting. For each
generated partial painting, the algorithm constructs the polytope as described in \autoref{sec:heuristic_mode} and then performs a vertex enumeration as
described in the
vertex filtering search. It remains to determine a good stopping criterion for branching.

Since the vertex enumeration algorithm has a good performance for low-dimensional polytopes, we wish to use the dimension as the stopping criterion. 
However the actual dimension of the polytope given by a painting is unknown, unless it has been constructed, when it is too late. Therefore, instead of the actual dimension, we use expected dimension as the stopping criterion in our code. This expected dimension can be computed efficiently without calling PPL to construct the polytope. 
We set up a $q$-column matrix (\sage{cs\underscore matrix}) to record the equality constraints on $(\pi_0, \dots, \pi_q)$, which are $\pi_0=0$, the symmetry constraints and the additivity constraints specified by the painting.
The matrix is maintained dynamically during the backtracking process. Define the expected dimension to be the co-rank
of the equation system:
$\texttt{exp\underscore dim} := q - \rk(\texttt{cs\underscore matrix})$. 

The algorithm switches from backtracking to vertex enumeration once the
expected dimension becomes smaller than a certain
threshold. \autoref{tab:best-dim-threshold} shows that a value around $11$ is the best empirical threshold for finding an extreme function with many slopes quickly.

The combination of vertex filtering search and heuristic backtracking search produces a more powerful search algorithm, which is called ``combined mode'' in our code. We summarize it as \hyperref[alg:combined-mode]{Algorithm~\ref*{alg:combined-mode}}.

\begin{algorithm}
\caption{combined mode}
\label{alg:combined-mode}
\begin{enumerate}
\item \label{itm:heuristic-in-combined} Run the heuristic search
  (\hyperref[alg:heuristic-mode]{Algorithm~\ref*{alg:heuristic-mode}}), 
with one more stopping criterion added to its step~\ref{itm:green-subnode}:
\begin{itemize}
	\item append the new equations to \sage{cs\underscore matrix};
	\item compute $\texttt{exp\underscore{}dim} := q - \rk(\texttt{cs\underscore{}matrix})$;
	\item if \sage{exp\underscore dim} $\leq$ threshold (empirically,
          threshold $= 11$), \\
          \textbf{output} the partial painting and backtrack;
\end{itemize}
\item For each partial painting returned by
  phase~\ref{itm:heuristic-in-combined}, construct the corresponding polytope
  and run the vertex filtering search
  (\hyperref[alg:vertex-filtering-mode]{Algorithm~\ref*{alg:vertex-filtering-mode},
  steps 3--4}).
\end{enumerate}
\end{algorithm}

\subsection{Results}
Using the combined search (\hyperref[alg:combined-mode]{Algorithm~\ref*{alg:combined-mode}}), our code\footnote{By running \sage{search\underscore kslope\underscore example(k\underscore slopes, q, f, mode='combined')} with various values of \sage{k\underscore slopes, q, f}. For example, \sagefunc{kzh_7_slope_1} can be obtained by setting \sage{k\underscore slopes=7; q=33; f=11}.} was able to find up to 7-slope extreme functions for $q \leq 34$, namely \sagefunc{kzh_7_slope_1} to \sagefunc{kzh_7_slope_4}.

\begin{landscape}
\begin{table}[h]
  \caption{Vertex enumeration in high dimension vs.\@\ Combinatorial explosion in branching}
  \label{tab:best-dim-threshold}
  \begin{minipage}{\linewidth}
  \centering   
  \begin{tabular}[t]{*{4}c@{\hspace{1\tabcolsep}}*{2}c@{\hspace{2\tabcolsep}}*{2}c@{\hspace{2\tabcolsep}}*{2}c@{\hspace{2\tabcolsep}}*{2}c@{\hspace{2\tabcolsep}}*{2}c@{\hspace{2\tabcolsep}}*{2}c@{\hspace{2\tabcolsep}}}
    \toprule 
	
	& \multicolumn{3}{c}{$q=25$}
	& \multicolumn{2}{c}{$q=26$}
	& \multicolumn{2}{c}{$q=27$}
	& \multicolumn{2}{c}{$q=28$}
	& \multicolumn{2}{c}{$q=29$}
	& \multicolumn{2}{c}{$q=30$}
	& \multicolumn{2}{c}{$q=31$} \\
	\cmidrule(lr){2-4} \cmidrule(lr){5-6} \cmidrule(lr){7-8}
	\cmidrule(lr){9-10} \cmidrule(lr){11-12} \cmidrule(lr){13-14}
	\cmidrule(lr){15-16}
	$k$
	& 6 & 6 & 6
	&\hspace{1\tabcolsep} 6 & 6
	&\hspace{1\tabcolsep} 6 & 6
	&\hspace{1\tabcolsep} 6 & 6
	&\hspace{1\tabcolsep} 6 & 7
	&\hspace{1\tabcolsep} 7 & 6
	&\hspace{1\tabcolsep} 7 & 7 \\
	$f$
	& 1 & 7 & 8
	&\hspace{1\tabcolsep} 1 & 9
	&\hspace{1\tabcolsep} 1 & 9
	&\hspace{1\tabcolsep} 1 & 9
	&\hspace{1\tabcolsep} 1 & 10
	&\hspace{1\tabcolsep} 1 & 10
	&\hspace{1\tabcolsep} 1 & 10 \\
	\midrule
    \multicolumn{1}{c}{}
    & \multicolumn{15}{c}{number of $\geq k$-slope solutions} \\
    \midrule 
    
	& 2 & 1 & 1
	&\hspace{1\tabcolsep} 8 & 4
	&\hspace{1\tabcolsep} 14 & 1
	&\hspace{1\tabcolsep} 26 & 17
	&\hspace{1\tabcolsep} 60 & 1
	&\hspace{1\tabcolsep} 3 & 30
	&  &  \\
    \midrule
    \multicolumn{1}{c}{}
    & \multicolumn{15}{c}{running time (s) in vertex filtering search} \\
    \midrule 
    v-enumeration
    & 59 & 40 & 28
    &\hspace{1\tabcolsep} 375 & 322
    &\hspace{1\tabcolsep} 439 & 706
    &\hspace{1\tabcolsep} 3866 & 3806
    &\hspace{1\tabcolsep} 3728 & 3626
    &\hspace{1\tabcolsep} 23642 & 2880
    & & \\
    first solution
    & 86 & 42 & 47
    &\hspace{1\tabcolsep} 378 & 330 
    &\hspace{1\tabcolsep} 440 & 757
    &\hspace{1\tabcolsep} 3873 & 3845
    &\hspace{1\tabcolsep} 3739 & 3650
    &\hspace{1\tabcolsep} 23747 & 2889 
    & & \\
    all solutions
    & 92 & 63 & 51
    &\hspace{1\tabcolsep} 454 & 370
    &\hspace{1\tabcolsep} 555 & 818
    &\hspace{1\tabcolsep} 4211 & 4049
    &\hspace{1\tabcolsep} 4369 & 3958
    &\hspace{1\tabcolsep} 24506 & 3256
    & & \\
    \midrule 
    \multicolumn{1}{c}{threshold}
    & \multicolumn{15}{c}{running time (s) in combined search to find the first $\geq k$-slope solution} \\
    \midrule
    \hphantom{0}5
    & 6 & 122 & 666
    &\hspace{1\tabcolsep} 5 & 1369
    &\hspace{1\tabcolsep} 20 & 1932
    &\hspace{1\tabcolsep} 282 & 2875
    &\hspace{1\tabcolsep} 20 & 7809
    &\hspace{1\tabcolsep} 1884 & 5896
    &\hspace{1\tabcolsep} \hphantom{0}5181 & 35455 \\
    \hphantom{0}6
    & 4 & \hphantom{0}66 & 397
    &\hspace{1\tabcolsep} 3 & \hphantom{0}921
    &\hspace{1\tabcolsep} 14 & 1322
    &\hspace{1\tabcolsep} \hphantom{0}64 & 2084
    &\hspace{1\tabcolsep} 12 & 6278
    &\hspace{1\tabcolsep} 1587 & 1529
    &\hspace{1\tabcolsep} \hphantom{0}4527 & 24243 \\
    \hphantom{0}7
    & 3 & \hphantom{0}32 & 224
    &\hspace{1\tabcolsep} 2 & \hphantom{0}518 
    &\hspace{1\tabcolsep} 11 & \hphantom{0}845
    &\hspace{1\tabcolsep} \hphantom{0}55 & 1292
    &\hspace{1\tabcolsep} 15 & 4807
    &\hspace{1\tabcolsep} 1492 & 1031
    &\hspace{1\tabcolsep} \hphantom{0}5728 & 19043 \\
    \hphantom{0}8 
    & 3 & \hphantom{0}13 & 121
    &\hspace{1\tabcolsep} 5 & \hphantom{0}267 
    &\hspace{1\tabcolsep} 20 & \hphantom{0}641
    &\hspace{1\tabcolsep} 101 & \hphantom{0}779 
    &\hspace{1\tabcolsep} 15 & 4823 
    &\hspace{1\tabcolsep} 3782 & \hphantom{0}449
    &\hspace{1\tabcolsep} 21821 & \hphantom{0}8604 \\
    \hphantom{0}9
    & 1 & \hphantom{00}4 & \hphantom{0}56
    &\hspace{1\tabcolsep} 5 & \hphantom{0}135 
    &\hspace{1\tabcolsep} 20 & \hphantom{0}352
    &\hspace{1\tabcolsep} \hphantom{0}49 & \hphantom{0}516
    &\hspace{1\tabcolsep} \hphantom{0}2 & 3194 
    &\hspace{1\tabcolsep} 2032 & \hphantom{0}242
    &\hspace{1\tabcolsep} 24822 & \hphantom{0}5487 \\
    10
    & 1 & \hphantom{00}4 & \hphantom{0}15
    &\hspace{1\tabcolsep} 4 & \hphantom{00}18
    &\hspace{1\tabcolsep} \hphantom{0}5 & \hphantom{0}121 
    &\hspace{1\tabcolsep} \hphantom{0}47 & \hphantom{0}150
    &\hspace{1\tabcolsep} \hphantom{0}1 & 1460
    &\hspace{1\tabcolsep} \hphantom{0}549 & \hphantom{00}99
    &\hspace{1\tabcolsep} \hphantom{0}8260 & \hphantom{0}2577 \\
    11
    & 2 & \hphantom{00}4 & \hphantom{0}15
    &\hspace{1\tabcolsep} 4 & \hphantom{00}39
    &\hspace{1\tabcolsep} \hphantom{0}5 & \hphantom{00}82 
    &\hspace{1\tabcolsep} \hphantom{0}27 & \hphantom{00}29
    &\hspace{1\tabcolsep} 45 & 1000
    &\hspace{1\tabcolsep} \hphantom{0}271 & \hphantom{00}40
    &\hspace{1\tabcolsep} \hphantom{0}1010 & \hphantom{0}1430 \\
    12
    & & &
    &\hspace{1\tabcolsep} 4 & \hphantom{00}38
    &\hspace{1\tabcolsep} \hphantom{0}5 & \hphantom{00}83
    &\hspace{1\tabcolsep} \hphantom{0}28 & \hphantom{00}29
    &\hspace{1\tabcolsep} 44 & \hphantom{0}932
    &\hspace{1\tabcolsep} \hphantom{0}306 & \hphantom{00}40
    &\hspace{1\tabcolsep} \hphantom{0}2352 & \hphantom{0}1186 \\
    13
    & & &
    & &
    & &
    &\hspace{1\tabcolsep} \hphantom{0}27 & \hphantom{00}28
    &\hspace{1\tabcolsep} 46 & \hphantom{0}928
    &\hspace{1\tabcolsep} \hphantom{0}308 & \hphantom{00}42
    &\hspace{1\tabcolsep} \hphantom{0}1269 & \hphantom{0}3365 \\
    14
    & & & 
    & & 
    & & 
    & & 
    & &
    &\hspace{1\tabcolsep} \hphantom{0}229 & \hphantom{00}41
    &\hspace{1\tabcolsep} \hphantom{0}1227 & \hphantom{0}3637 \\
    \bottomrule  
  \end{tabular} 
  \end{minipage}
\end{table}
\end{landscape}

\section{Targeted search for extreme functions with many slopes}
\label{sec:targeted_search}

\begin{figure}[t]%
\includegraphics[width=.49\linewidth]{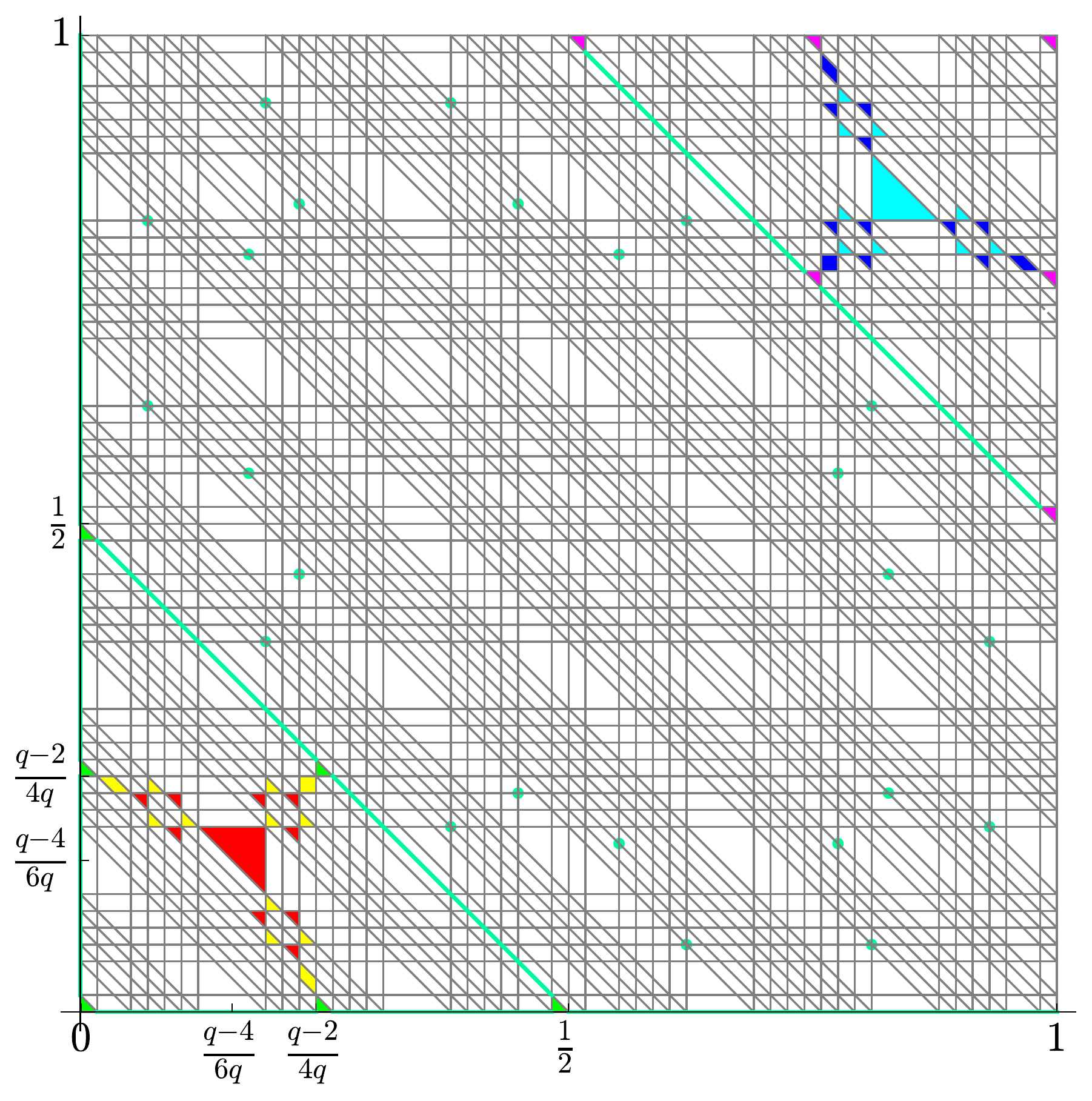}
\includegraphics[width=.49\linewidth]{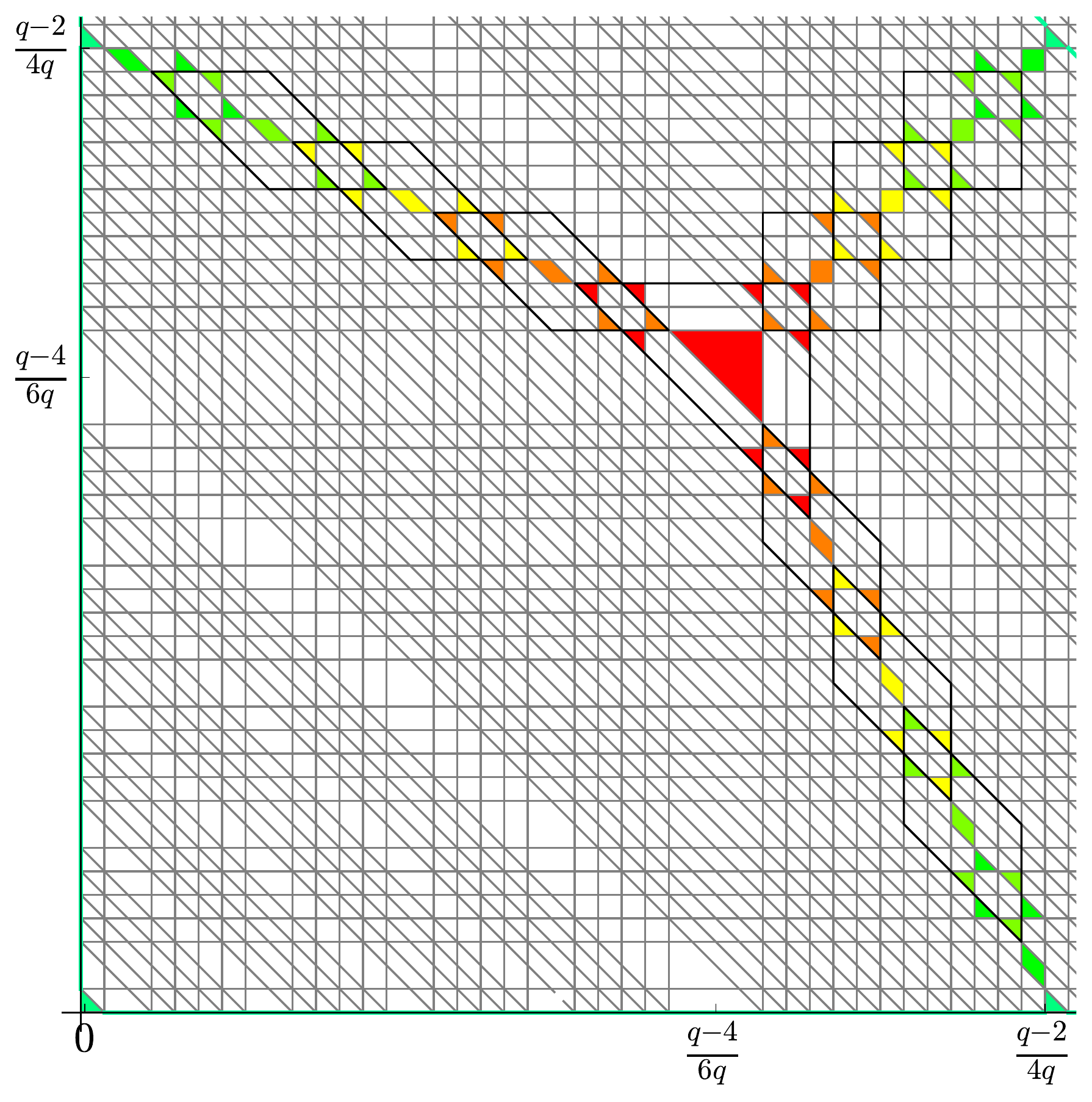}
\caption{Special patterns on the two-dimensional polyhedral complex $\Delta\P_{\frac1q\Z}$. 
\textit{Left,} the $\Delta\P_{\frac1q\Z}$ of the $6$-slope extreme function
\sagefunc{kzh_6_slope_1} with $q=58$. We observe that the additive triangles
are located in the lower left and upper right corners. The function has the
same slopes on the intervals that are projections of the same color additive
triangles. The 6-pointed star patterns appear several times. \textit{Right,} the
lower-left corner of $\Delta\P_{\frac1q\Z}$ of the $10$-slope extreme function
\sagefunc{kzh_10_slope_1} with $q=166$, where we see that the 6-pointed stars
are actually the result of additivity patterns within certain intersecting
quadrilaterals (\emph{black}), which connect like
links of three chains. The detailed structure is described in
\autoref{fig:special_patterns_step}.}%
\label{fig:special_patterns}%
\end{figure}

We observed that many of these newly discovered extreme functions with many
slopes possess the invariance:  $f = 1/2$ and $\pi_i = \pi_{q - i}$ ($ 0 \leq i \leq q/2$). 
Targeting the search to functions with the invariance property allowed us to find more new extreme functions with many slopes whose values of $q$ were twice as large as before.
In addition, their painting on the complex $\Delta\P_{\frac1q\Z}$ often includes special patterns,
as shown by \autoref{fig:special_patterns} for example.
\subsection{Construction of a family of prescribed partial paintings}
We then targeted the search to functions for larger values of~$q$ with prescribed partial paintings that mimic these patterns.
Let $q=36r+22$, where $r \in \Z, r \geq 1$. We construct the prescribed
partial painting on the
complex $\Delta\P_{\frac1q\Z}$ in $2(r+2)$ steps as follows. 

\begin{figure}
    \centering
    \begin{minipage}{.5\textwidth}
        \centering
        \includegraphics[scale=0.5]{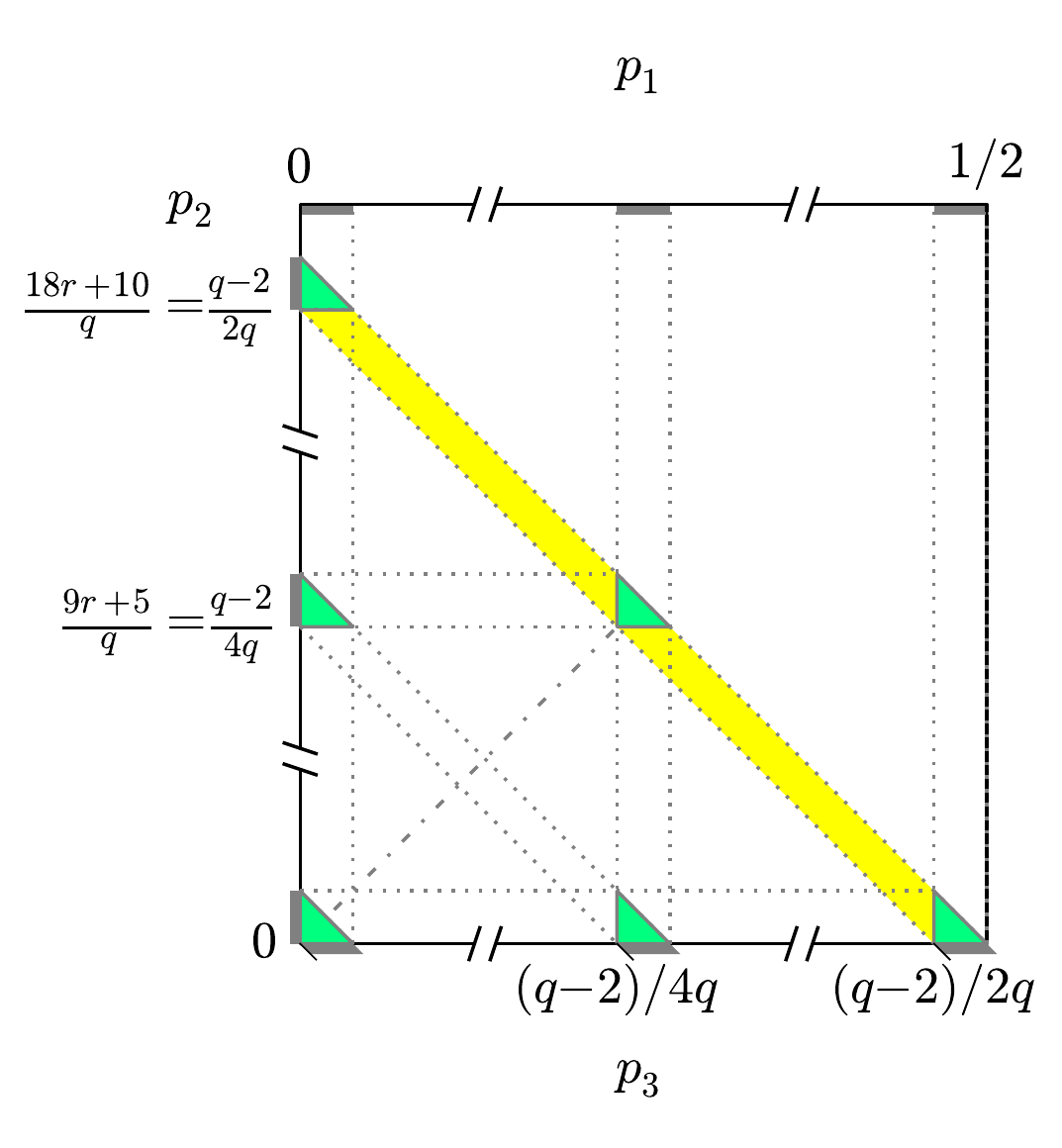} \\
        $\Delta\P^{0}$ of \textit{Step $0$}
    \end{minipage}%
    \begin{minipage}{0.5\textwidth}
        \centering
		\includegraphics[scale=0.5]{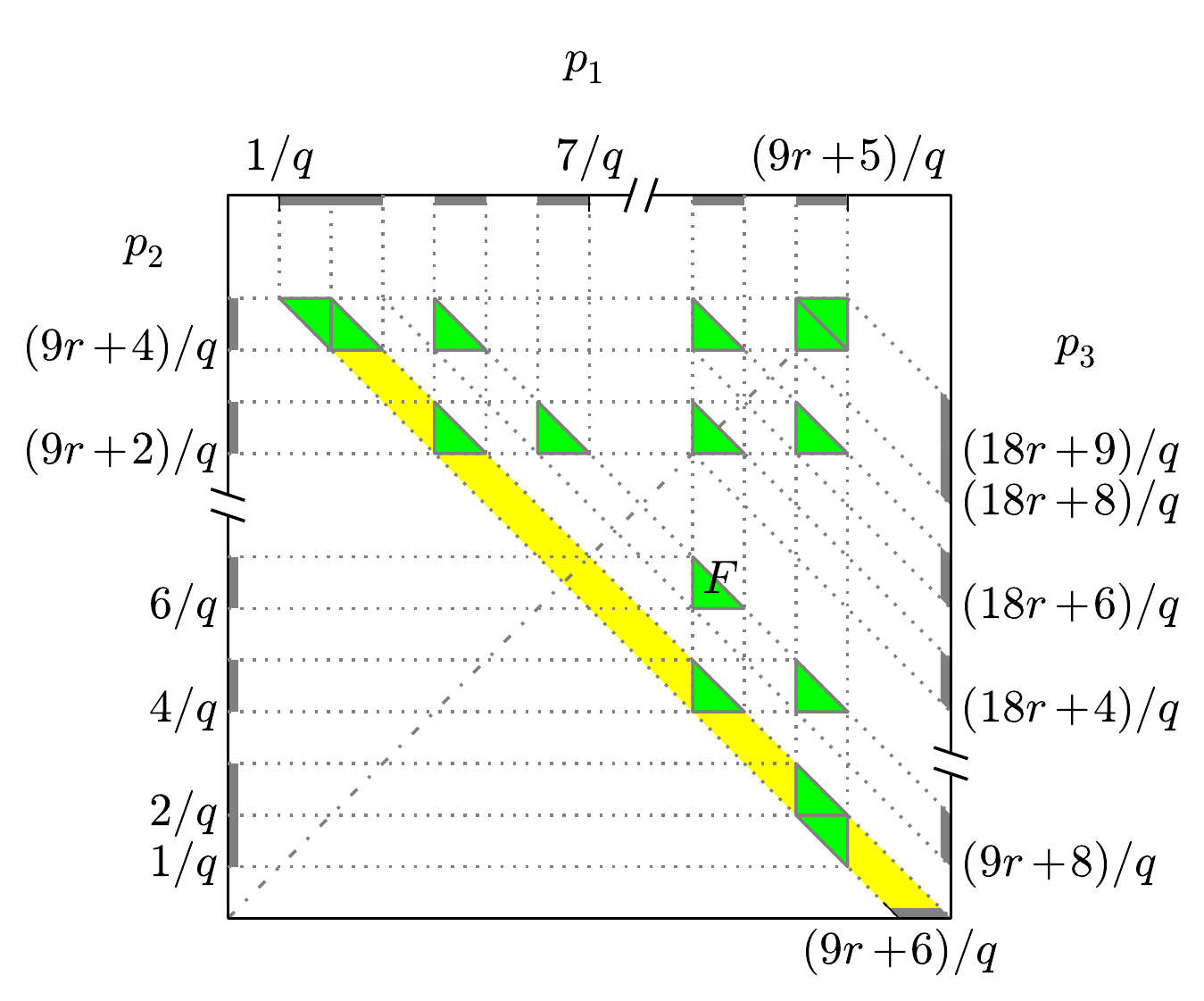} \\
		$\Delta\P^{1}$ of \textit{Step $1$}
    \end{minipage} \\
        \centering
		\includegraphics[scale=0.5]{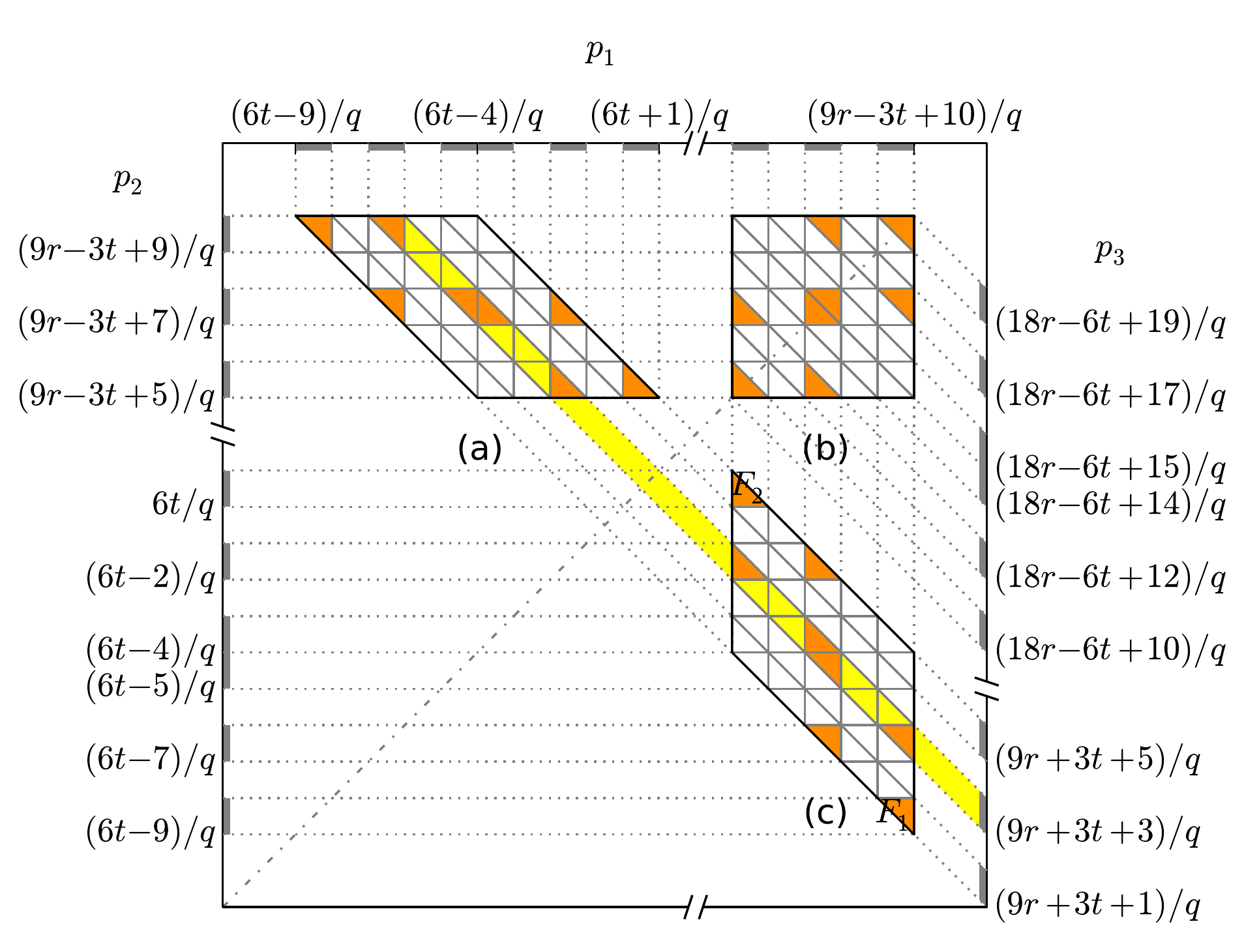} \\
		 $\Delta\P^{t}$ of \textit{Step $t$} \\
		\includegraphics[scale=0.5]{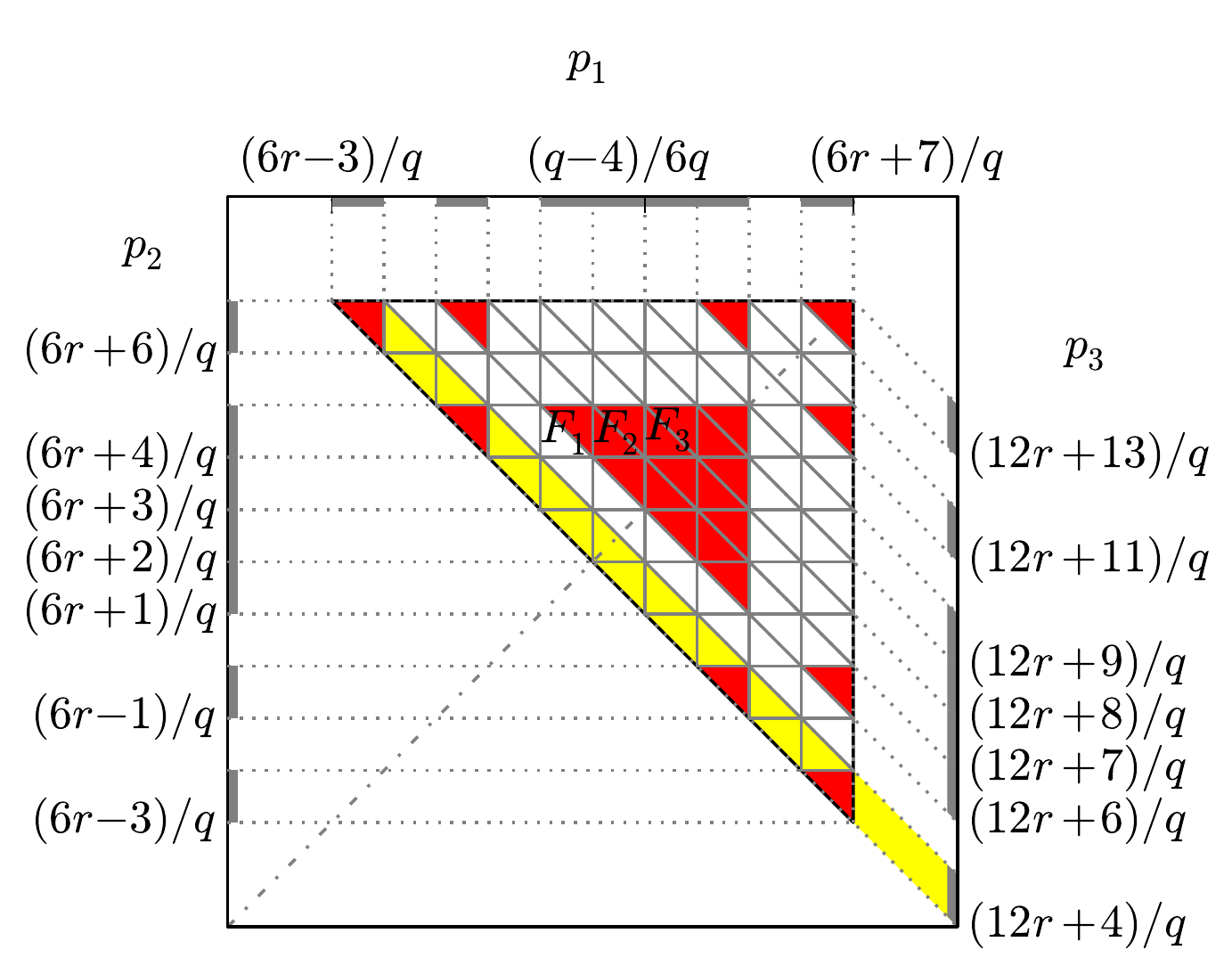} \\
		$\Delta\P^{r+1}$ of \textit{Step $(r+1)$}
    \caption{Prescribed partial paintings.}
    \label{fig:special_patterns_step}
\end{figure}

\textit{Step $0$}: Paint the lower triangles whose lower left corners are the vertices $\ColVec{0}{0}, \ColVec{0}{(q-2)/4q}, \ColVec{0}{(q-2)/2q}, \ColVec{(q-2)/4q}{(q-2)/4q}, \ColVec{(q-2)/4q}{0}$ and $\ColVec{(q-2)/2q}{0}$;
this is illustrated in \autoref{fig:special_patterns_step}--1.

\textit{Step $1$}: Paint the upper triangles whose upper right corners are the vertices  $\ColVec[q]{i}{j}$ for $\ColVec{i}{j}=\ColVec{2}{9r+5}, \ColVec{9r+5}{9r+5}$ and $\ColVec{9r+5}{2}$.
Paint the lower triangles whose lower left corners are the vertices $\ColVec[q]{i}{j}$ for $\ColVec{i}{j} = \ColVec{2}{9r+4}, \ColVec{4}{9r+4}, \allowbreak \ColVec{4}{9r+2}, \allowbreak \ColVec{6}{9r+2}, \ColVec{9r+2}{9r+4}, \ColVec{9r+4}{9r+4}, \ColVec{9r+2}{9r+2}, \ColVec{9r+4}{9r+2},\ColVec{9r+4}{4}, \ColVec{9r+2}{4}$ and $\ColVec{9r+2}{6}$;
see \autoref{fig:special_patterns_step}--2.


\textit{Steps $t=2,3,\dots,r$}:
Paint the parallelogram whose vertices are $\ColVec[q]{i}{j}$ for $\ColVec{i}{j} = \ColVec{6t-9}{9r-3t+10}, \ColVec{6t-4}{9r-3t+10}, \ColVec{6t-4}{9r-3t+5}, \ColVec{6t+1}{9r-3t+5}$ with the orange pattern shown in \autoref{fig:special_patterns_step}--3(a)).
Paint the square whose vertices are $\ColVec[q]{i}{j}$ for $\ColVec{i}{j} = \ColVec{9r-3t+5}{9r-3t+10}, \ColVec{9r-3t+10}{9r-3t+10}, \ColVec{9r-3t+5}{9r-3t+5}, \ColVec{9r-3t+10}{9t-3t+5}$ with the orange pattern shown in \autoref{fig:special_patterns_step}--3(b).
Paint the parallelogram whose vertices are $\ColVec[q]{i}{j}$ for $\ColVec{i}{j} = \ColVec{9r-3t+5}{6t+1}, \ColVec{9r-3t+5}{6t-4}, \ColVec{9r-3t+10}{6t-4}, \ColVec{9r-3t+10}{6t-9}$ with the orange pattern shown in \autoref{fig:special_patterns_step}--3(c).

\smallskip
\textit{Step $(r+1)$}: Paint the triangle whose vertices are$\ColVec[q]{i}{j}$ for $\ColVec{i}{j}=\ColVec{6r-3}{6r+7}, \allowbreak \ColVec{6r+7}{6r+7}, \ColVec{6r+7}{6r+3}$ with the red pattern shown in \autoref{fig:special_patterns_step}--4.

\smallskip
\textit{Step $0$} to \textit{Step $(r+1)$} construct the painting on the lower
left triangular part of the complex $\Delta\P_{\frac1q\Z}$. The painting on the
upper right triangular part of the complex $\Delta\P_{\frac1q\Z}$ will then be determined through a mapping. Specifically, in \textit{Step $t$} for $t=r+2, r+3, \dots, 2r+3$, we paint the triangles whose images under the mapping $\ColVec{x}{y} \mapsto \ColVec{1-x}{1-y}$ are colored in \textit{Step $(2r+3-t)$}. We also paint the diagonal lines $\{\ColVec{x}{y} \colon x+y = \tfrac{1}{2}, 0\leq x\leq \tfrac{1}{2}\}$ and $\{\ColVec{x}{y} \colon x+y = \tfrac{3}{2}, \tfrac{1}{2}\leq x \leq 1\}$, which correspond to the symmetry condition of minimal valid functions.
In the following, we shall refer to the painting constructed as above on $\Delta\P_{\frac1q\Z}$ as the \textit{prescribed partial painting}.


\subsection{Properties of functions satisfying the prescribed partial paintings}
Recall the notion of \textit{connected component} discussed in \autoref{sec:incremental_computation}.
Connected components of a painting on $\Delta\P_{\frac1q\Z}$ are disjoint subsets of $\{[\tfrac{i}{q},\tfrac{i+1}{q}] \colon \allowbreak  i= 0, 1, \dots, q-1\}$. They
satisfy the following properties. If $F$ is an elementary upper or lower triangle whose vertices are colored on the painting, then $p_1(F), p_2(F)$ and $p_3(F)$ are in the same connected component. 
If $F$ is an elementary horizontal, vertical or diagonal edge whose vertices are colored on the painting, then the two intervals among its projections $p_1(F), p_2(F), p_3(F)$ are in the same connected component. In particular, the colored diagonal lines corresponding to the symmetry condition yield that $[\tfrac{i}{q}, \tfrac{i+1}{q}]$ and $[\tfrac{j}{q}, \tfrac{j+1}{q}]$ are in the same connected component, where $i+j=\tfrac{q}{2}-1=18r+10$.

\begin{lemma}
\label{lem:patterns}
Let $q=36r+22$, where $r \in \Z, r \geq 1$.
Let the prescribed partial painting on $\Delta\P_{\frac1q\Z}$ be constructed as above. 
Then all intervals are directly covered. The prescribed
partial painting induces exactly $2(r+2)$ connected components.
\end{lemma} 
The proof of this and the following results appears in \autoref{sec:proofs-of-thm-in-targeted-search}.

Suppose $r\in\Z$, $r\geq 1$, $q=36r+22$ and $f=1/2$. Let $\Pi_r$ be the set of
continuous piecewise linear minimal valid functions $\pi$ with breakpoints in $\frac{1}{q}\Z$,
satisfying in addition the invariance condition $\pi(x) = \pi(1-x)$ for $0
\leq x \leq \frac{1}{2}$ and $\Delta\pi(x,y)=0$ for any $0\leq x, y \leq 1$
such that the point $(x,y)$ is colored in the prescribed partial painting on
$\Delta\P_{\frac1q\Z}$. 

Clearly we can describe the functions~$\pi\in\Pi_r$ in the space of the variables
$(\pi_0, \pi_1, \dots, \pi_q)$; let us denote the corresponding
polytope by~$V_r$.  Due to \autoref{lem:patterns}, we can describe them also
in the space of the slope values $(s_0, s_1, \dots, s_{r+1})$ on the connected
components; let us denote the corresponding polytope by~$S_r$.  (Note that
$s_t = -s_{2r+3-t}$ by the invariance condition, for 
$t=r+2, r+3, \dots, 2r+3$.)  
\begin{lemma}
\label{lem:isomorphism}
There is a linear isomorphism between the polytope $V_r$ in the $(\pi_0, \pi_1, \dots, \pi_q)$ variables and the polytope $S_r$ in the $(s_0, s_1, \dots, s_{r+1})$ variables.
\end{lemma}

The proof (in \autoref{sec:proofs-of-thm-in-targeted-search}) gives an explicit mapping between the values $\pi_i$ and the slopes $s_t$. With the help of these formulas,  we can then show that the slope values are non-increasing.

\begin{lemma}
\label{lem:ordering_of_slopes}
Let $s_0, s_1, \dots, s_{r+1}$ be as above. Then  $s_0 \geq s_1 \geq \dots \geq s_{r+1}$.
\end{lemma} 







The computer-based search can now be run in the space of the slope variables $(s_0, s_1,
\dots, s_{r+1})$, which has the benefit of having a much lower
dimension.  However, if $r$ is large, the search is still nontrivial.  In order to
speed up the search we prescribed extra additivity constraints.  This approach
was successful in finding extreme functions with up to 28 slopes.\footnote{By running \sage{pattern\underscore extreme(r, k\underscore slopes)} with various values of \sage{r} and \sage{k\underscore slopes}. For example, \sagefunc{kzh_28_slope_1} can be obtained by setting \sage{r=21; k\underscore slopes=28}. In order to reduce the running time of the targeted search when $r$ is large, the code \sage{pattern\underscore extreme()} imposes some extra colored vertices on the prescribed partial painting. They correspond to extra additivity constraints on the function $\pi$, which are often satisfied by the previously discovered many-slope extreme functions. Concretely, when $r \geq 16$, we assume that $\Delta\pi_{x,y} = 0$ for $(x,y) = (6r + 5, 36r + 18),
 (6r + 7, 36r + 10),
 (6r + 7, 36r + 12),
 (6r + 10, 36r + 3),
 (6r + 11, 36r),
 (9r - 18, 9r - 18),
 (9r - 12, 9r - 12),
 (9r - 9, 9r - 9),
 (9r - 3, 9r - 3),
 (9r + 3, 9r + 3)$.}

Our search also revealed that in general, we cannot expect the
existence of an extreme point for which the sequence of slope values is
strictly decreasing ($s_0 > s_1 > \dots > s_{r+1}$). 
However, we conclude this section with a weaker conjecture.
\begin{conjecture}
  There exists an extreme point of the polytope~$S_r$ 
  with $\Omega(r)$ different slope values $s_i$. 
\end{conjecture}

We abandoned work on this conjecture in late July 2015 when Basu et al. 
\cite{bcdsp:arbitrary-slopes} announced a different construction that gives extreme functions
with an arbitrary prescribed number of slopes.

%
%
%

\subsection{Result: Extreme functions with many slopes}

The targeted search 
was very successful in finding functions with large numbers of slopes\footnote{We have made the functions available as part of the Electronic Compendium \cite{electronic-compendium} as \sage{kzh\underscore{}7\underscore{}slope...} and \sage{kzh\underscore{}28\underscore{}slope...}, etc.}.  We thus
obtained the following result, which we have stated already in the
introduction.

\begin{theorem}
There exist continuous piecewise linear extreme functions with $2$, $3$, $4$,
$5$, $6$, $7$, $8$, $10$, $12$, $14$, $16$, $18$, $20$, $22$, $24$, $26$, and $28$ slopes.
\end{theorem}




\appendix
\section{Implementation details}

In this appendix, we describe some aspects of our implementation in SageMath
\cite{sage}, an open-source mathematics software system that uses Python and
Cython as its primary programming languages and interfaces with various
existing packages. 

\subsection{SageMath interface for vertex enumeration}
\label{sec:dd-ppl}
PPL is a standard package in SageMath and comes with an efficient 
Cython interface.
Our code uses in particular
the class \sage{C\underscore{}Polyhedron} for computations with closed convex polyhedra.
A polytope 
can be built starting from a system of constraints \sage{cs} of class \sage{Constraint\underscore{}System} via \sage{p = C\underscore{}Polyhedron(cs)}, where the constraint system \sage{cs} is a finite set of linear equality or inequality constraints (class \sage{Constraint}). One calls \sage{p.minimized\underscore{}generators()} to enumerate the vertices of \sage{p}. 
%
PPL also allows for feasibility checks and satisfiability checks (see \autoref{sec:backtracking-search}). 
The feasibility check can be realized by calling
\sage{p.is\underscore{}empty()}.  The satisfiablity check efficiently
tests whether a given inequality or equation~$\sage{c}$ is satisfied by all points in a
polytope~$\sage{p}$
. It is accessed by calling
\sage{p.\allowbreak
  relation\underscore{}with(c).\allowbreak implies(Poly\underscore{}Con\underscore{}Relation.is\underscore{}\allowbreak
  included())}.

Once lrslib \cite{avis1998computational, avis2000revised} has been installed as an optional package in SageMath, it is possible to call the programs \sage{lrs} and \sage{redund} from SageMath. Our code includes a SageMath interface that reads or writes polytopes in the lrslib format. The lrslib command \sage{redund} can thus be used in conjunction with PPL as a preprocessor for vertex enumeration.

\subsection{SageMath interface for linear programming}
\label{sec:lp-ppl-glpk}

\begin{sloppypar}
We pointed out in \autoref{sec:warm-start-lp} the necessity of using an LP
solver with warm-starting capability for the feasibility and satisfiability
checks described above in the high-dimensional case.  
We use the SageMath class \sage{MixedIntegerLinearProgram} as an LP
modeling system. 
Within this framework, a new LP problem \sage{m} can be created by \sage{m =
  MixedIntegerLinearProgram(maximization=True, solver = "GLPK")}, requesting
the GLPK solver as its numerical backend.  In contrast to other backend implementations,
including PPL's rational LP solver, the GLPK backend has the crucial warm
starting capability.
We call \sage{v = m.new\underscore variable(real=True, nonnegative=True)} to define a Python dictionary \sage{v} of non-negative continuous variables for the problem \sage{m}. The upper and lower bound of a variable, say \sage{v[0]}, can be changed via \sage{m.set\underscore max(v[0], max)} and \sage{m.set\underscore min(v[0], min)} respectively. If the variable is unbounded above or below, then one sets \sage{max=None} or \sage{min=None} respectively. The method \sage{m.add\underscore constraint(linear\underscore function, max, min)} sets up a new constraint $\sage{min}\leq\sage{linear\underscore function}\leq\sage{max}$ for the problem \sage{m}. The objective function
of \sage{m} is defined by \sage{m.set\underscore objective(obj)}. For a 
feasibility check,  
we can use 
\sage{obj=None}. 
We request that GLPK use the simplex method to solve the LP via \sage{m.solver\underscore parameter(backend.glp\underscore simplex\underscore or\underscore intopt, backend.glp\underscore simplex\underscore only)}. According to the setting \sage{m.solver\underscore parameter("primal\underscore v\underscore dual", "GLP\underscore PRIMAL")} or \sage{m.solver\underscore parameter("primal\underscore v\underscore dual", "GLP\underscore DUAL")}, the primal or dual simplex method is applied respectively.
We call \sage{m.solve(objective\underscore only=True)} to solve for the optimal value. If it signals a \sage{MIPSolverException}, then the problem is infeasible.
\end{sloppypar}
References:

\section{Limitations of search based on $q \times v$ grid discretization}
\label{sec:limitation-of-grid}

In this section, we discuss limitations of the search based on $q \times v$
grid discretization, an alternative search strategy that was used by Chen
\cite{chen} and Hildebrand (2013, unpublished).

Consider continuous piecewise linear functions $\pi \colon \R/\Z \to [0,1]$, with breakpoints in $\tfrac{1}{q} \Z$ for some $q \in \Z_+$ and $\pi(0) = 0$. Suppose without loss of generality that $f \in \tfrac{1}{q} \Z$. 

As mentioned in \autoref{sec:intro_infinite_literature}, there are two natural ways to discretize the space of functions $\pi$: discretizing function values $\pi_i=\pi(\frac{i}{q})$ for $i \in \{ 0, \dots, q\}$ and discretizing slope values $qs_i$ on $[\frac{i-1}{q}, \frac{i}{q}]$ for $i \in \{ 1, \dots, q\}$. See again \autoref{fig:q_v_grid}. 
The following lemma shows that they are equivalent.
\begin{lemma}
\label{lem:values-slopes-same-denominator}
Let $v$ be a positive integer. The following are equivalent:
\begin{enumerate}
\item $\pi_i \in \frac{1}{v}\Z$ for each $i \in \{ 0, \dots, q\}$.
\item $s_i \in \frac{1}{v}\Z$ for each $i \in \{ 1, \dots, q\}$.
\end{enumerate}
\end{lemma}
\begin{proof}
Since $\pi_0 = \pi_q = 0$ and $ s_i = \pi_i - \pi_{i-1}$ for $i = 1, \dots, q$, the lemma follows.
\end{proof}

\subsection{A lower bound on (a proxy for) arithmetic complexity}

In the following, we investigate the worst-case complexity of the search based
on $q \times v$ grid discretization, by estimating the arithmetic complexity
of extreme functions, i.e., the largest value $v$ needed for any extreme
function $\pi$ with breakpoints in $\frac{1}{q}\Z$.  

We use the notations $d_{\text{ext}}$, $d_{\text{ver}}$, and $d_{\text{bas}}$,
satisfying $d_{\text{ext}} \leq d_{\text{ver}} \leq d_{\text{bas}}$,  
and the constraint system $A \x = \b$ of the polytope
$\Pi_f(\frac{1}{q}\Z/\Z)$, which were introduced in
\autoref{sec:vertex-filtering-results}. It is hard to determine the precise
value of $d_{\text{ext}}$ as a function of $q$, and we are not able to show 
an exponential lower bound for it. For estimating the growth rate
of $d_{\text{ext}}$, we are satisfied with a simplified study, using
$d_{\text{bas}}$ as a proxy. 
We show the following exponential lower bound.

\begin{lemma}
\label{lem:denominators}
Let $q \geq 3$ be an odd positive integer. Let $f \in \frac{1}{q}\Z$, $0<f<1$,
such that $qf$ and $q$ are coprime integers. Let $A \x = \b$, $\x\geq \ve0$ be the constraint
system of \autoref{thm:finite-minimal} written in the standard form. Then
$d_{\mathrm{bas}}$, the maximum absolute value of the determinants of simplex
basis matrices of $A$, is at least $2^{\frac{q-1}{2}}$.
\end{lemma}

\begin{proof}
It suffices to show the existence of a basis matrix $B$ of $A$ with $\mathopen|\det B| \geq 2^{\frac{q-1}{2}}$. To find such a $B$, we first prove the following claim. Because $q$ is odd, the operation of multiplying by $2 \pmod q$ is invertible. For $x \in \{0,1,\dots,q-1\}$, denote the unique $y \in \{0,1,\dots,q-1\}$ satisfying $2y = x \pmod q$ by $x/2$.


\begin{claim}{}\mbox{}\!\!\footnote{Thanks go to Xuancheng Shao for the help in proving this claim.}\;
\label{claim:sequence}
Let $q$ and $f$ be as above. There exists a sequence $(a_0, a_1, \dots,
\allowbreak a_{q-1})$ of integers with $a_0=0$, $a_1=qf$ and $a_2 = qf/2 \pmod q$ such that the following conditions hold:
\begin{enumerate}
\item for odd $i>1$ we have $a_i = a_j/2\pmod q$ for some $j<i$;
\item for even $i>2$ we have $a_i = qf-a_{i-1} \pmod q$.
\item $\{a_0,a_1,\dots,a_{q-1}\} = \{0,1,\dots,q-1\}$.
\end{enumerate}
\end{claim}


\begin{proof}
We construct the sequence as follows. Suppose that $a_0,a_1,\dots,\allowbreak a_k$ are determined for some even $k\geq 2$, such that conditions (1) and (2) are both satisfied for $i\leq k$, and that $a_0,a_1,\dots,a_k$ are all distinct. Let $S = \{a_0,a_1,\dots,a_k\}$. We choose $a_{k+1}$ and $a_{k+2}$ by selecting an element $s \in S$ such that $s/2 \pmod q \notin S$, and then taking $a_{k+1} = s/2 \pmod q$ and $a_{k+2} = qf-s/2 \pmod q$. It suffices to show the existence of such an $s\in S$ whenever $S\neq \Z/q\Z$. 

Suppose that $s/2 \in S$ for every $s\in S$. Since $qf\in S$ and $qf$ and $q$
are coprime, $S$ must contain the coset $qfH = \{qfh:h\in H\}$, where $H$ is the multiplicative subgroup of $(\Z/q\Z)^*$ generated by $2$. In particular, $qf,2qf\in S$. 

By conditions (2) and then (1), we deduce that $S$ also contains $qf-qfH$ and $(qf-qfH)H = qfH-qfH$. Applying this argument repeatedly, we see that $S$ contains $qfH-qfH+qfH-\dots\pm qfH$  for any number of iterations. Since $1,2\in H$, $qfH$ contains $qf$ and $2qf$, and thus any multiple of $qf$ can be written in the form $qfH-qfH+qfH-\dots\pm qfH$. Since $qf$ and $q$ are coprime, we conclude that $S$ contains all of $\Z/q\Z$.
\end{proof}

Define the row vectors $R_0, R_1, \dots, R_{q-1} \in \Z^q$ using the sequence
$a_0,a_1,\dots,\allowbreak a_{q-1}$ constructed above, as follows. Let $R_0$ be the row vector with the only nonzero entry $1$ appearing in the column indexed by $a_0=0$, corresponding to the constraint $\pi_0 = 0$. Let $R_2$ be the row vector with the only nonzero entry $2$ appearing in the column indexed by $a_2$, corresponding to the symmetric constraint $\pi_{a_2}+ \pi_{a_2} = 1$. For $i =1$ and $i = 4, 6, \dots, q-1$, let  $R_i$ be the row vector with two nonzero entries $1$ appearing in the columns indexed by $a_{i}$ and $a_{i-1}$, corresponding to the symmetric constraint $\pi_{a_i} + \pi_{a_{i-1}}= 1$. Finally, for  $i = 3, 5, \dots, q-2$, let $R_i$ be the row vector with nonzero entry $-2$ at index $a_{i}$ and entry $1$ at index $2a_{i} \pmod q$, corresponding to the subadditive constraint $\pi_{a_i}+\pi_{a_i} \geq \pi_{2a_i \bmod q}$.

The basis matrix $B$ is obtained by taking the slack variables for the subadditivity constraints $R_3, R_5, \dots, R_{q-2}$ in $A$ as non-basic variables and others as basic variables. To compute $\det B$, first expand out the columns corresponding to slack variables. We are left with a $q \times q$ matrix $B'$ consisting of the rows $R_0,R_1,\dots,R_{q-1}$, and $\mathopen|\det B| = \mathopen|\det B'|$. See \autoref{ex:determinant-q11} for the case of $q=11,f=3/11$.

To compute $\det B'$, start by expanding along the row $R_0$ containing a unique nonzero entry $1$ and end up with a new matrix with this row and column $a_0$ removed. In the second step, expand along $R_1$, noting that the only nonzero entry remaining in this row is $1$ at column $a_1$. We arrive at a new matrix with this row and column $a_1$ removed. In the third step, expand along $R_2$, noting that the only nonzero entry remaining in this row is $2$ at column $a_2$. We then arrive at a new matrix with this row and column $a_2$ removed. In general, during the ($k+1$)-st step, we expand along the row $R_{k}$ which contains a unique nonzero entry at column $a_{k}$, whose value is either $-2$ or $1$ depending on whether $k\geq 3$ is even or odd.

The computation terminates in $q$ steps. It follows that $\det B'$ is equal to the product of all these unique nonzero entries, $(q-1)/2$ of which are $\pm 2$ and the remaining are $1$. Thus $\mathopen|\det B| =  2^{\frac{q-1}{2}}$.
\end{proof} 

\begin{example} 
\label{ex:determinant-q11}
Consider the case $q=11,f=\frac{3}{11}$. 
By \autoref{claim:sequence},
we have the sequence
\[(a_0, a_1, \dots, a_{10}) = (0, 3, 7, 9, 5, 10, 4, 8, 6, 2, 1).\]
The following matrix $B'$ is a $q \times q$ submatrix of $A$, where the row $R_i$ corresponds to the:
\begin{itemize}
\item constraint $\pi_0 = 0$, for $i = 0$;
\item symmetric constraint $\pi_{a_2} + \pi_{a_2}= 1$, for $i =2$;
\item symmetric constraint $\pi_{a_i} + \pi_{a_{i-1}}= 1$, for $i =1, 4, 6, \dots, q-1$;
\item subadditive constraint $ -2\pi_{a_i} + \pi_{2a_i \text{ mod } q} \leq 0$, for $i = 3, 5, \dots, q-2$
\end{itemize}
of $\Pi_f(\frac{1}{q}\Z/\Z)$.

\[
  B'=\begin{blockarray}{c@{\hspace{5pt}}*{11}r@{\hspace{5pt}}cl}
   & \matindex{0} & \matindex{1} & \matindex{2} & \matindex{3} & \matindex{4} & \matindex{5} & \matindex{6} & \matindex{7} & \matindex{8} & \matindex{9} & \matindex{10} & &\\
    \begin{block}{(c@{\hspace{5pt}}*{11}r@{\hspace{5pt}}c)l}
      & 1 & 0 & 0 & 0 & 0 & 0 & 0 & 0 & 0 & 0 & 0 & & \matindex{$R_0$} \\
      & 1 & 0 & 0 & 1 & 0 & 0 & 0 & 0 & 0 & 0 & 0 & & \matindex{$R_1$} \\
      & 0 & 0 & 0 & 0 & 0 & 0 & 0 & 2 & 0 & 0 & 0 & & \matindex{$R_2$} \\
      & 0 & 0 & 0 & 0 & 0 & 0 & 0 & 1 & 0 &-2 & 0 & & \matindex{$R_3$} \\
      & 0 & 0 & 0 & 0 & 0 & 1 & 0 & 0 & 0 & 1 & 0 & & \matindex{$R_4$} \\
      & 0 & 0 & 0 & 0 & 0 & 0 & 0 & 0 & 0 & 1 &-2 & & \matindex{$R_5$} \\
      & 0 & 0 & 0 & 0 & 1 & 0 & 0 & 0 & 0 & 0 & 1 & & \matindex{$R_6$} \\
      & 0 & 0 & 0 & 0 & 0 & 1 & 0 & 0 &-2 & 0 & 0 & & \matindex{$R_7$} \\
      & 0 & 0 & 0 & 0 & 0 & 0 & 1 & 0 & 1 & 0 & 0 & & \matindex{$R_8$} \\
      & 0 & 0 &-2 & 0 & 1 & 0 & 0 & 0 & 0 & 0 & 0 & & \matindex{$R_9$} \\
      & 0 & 1 & 1 & 0 & 0 & 0 & 0 & 0 & 0 & 0 & 0 & & \matindex{$R_{10}$} \\
    \end{block}
  \end{blockarray}.
\]
\end{example}

\subsection{An upper bound}
We now prove the exponential upper bound stated
in~\autoref{sec:vertex-filtering-results}.

\begin{proof}[Proof of \autoref{lem:denominators-upper-bound}]
To compute $\det B$, we first expand out the columns corresponding to slack
variables as in the proof of \autoref{lem:denominators}. We are left with an
$n \times n$ matrix $B'$, where $n \leq q$. Denote the rows of $B'$ by
$R_0,R_1,\dots,R_{n-1} \in \Z^n$. We distinguish the types of constraints that
the rows $R_0,R_1,\dots,R_{n-1} \in \Z^n$ correspond to in
\autoref{thm:finite-minimal}. There is one constraint $\pi_0 =0$, for which
$\lVert R_i \rVert_2 = 1$; there are $m$ symmetric constraints $\pi_x +
\pi_{(qf-x) \bmod q} =1$, where $m = \frac{q}{2}, \frac{q+1}{2}$ or
$\frac{q+2}{2}$ depending on the parities of $n$ and $qf$, for which $\lVert
R_i \rVert_2 \leq \sqrt{2}$ if $x \neq (qf-x) \bmod q$ and $\lVert R_i
\rVert_2 \leq 2$ if $x = (qf-x) \bmod q$; and there are $n-1-m$ subadditive constraints, for which $\lVert R_i \rVert_2 \leq \sqrt{5}$. Therefore
\begin{equation}
\mathopen|\det B| = \mathopen|\det B'| \leq \prod_{i=0}^{n-1} \lVert R_i \rVert_2 \leq (\sqrt{5})^{n-m-1} (\sqrt{2})^{m+2} \leq 10^{{q}/{4}}.\tag*{\qedhere}
\end{equation}
\end{proof}

\subsection{Conclusion}

Although the question is not conclusively settled, Lemmas~\ref{lem:denominators} and \ref{lem:denominators-upper-bound} indicate that the value $v$ needed in the $q \times v$ grid discretization grows exponentially with $q$. The empirical results of $d_{\mathrm{ext}}$ and $d_{\mathrm{ver}}$ obtained by the vertex filtering search (see \autoref{sec:vertex-filtering-search}) confirm this exponential growth, as shown in \autoref{tab:denomiators} and \autoref{fig:denominators}.

We conclude that the search based on the $q \times v$ grid
discretization for breakpoints and function values (or, for breakpoints and
slope values, by \autoref{lem:values-slopes-same-denominator}) is not suitable
for an exhaustive search if $q$ is large, due to its high worst-case complexity.


\section{Proofs of the theorems in \autoref{sec:targeted_search}}
\label{sec:proofs-of-thm-in-targeted-search}
\begin{proof}[Proof of \autoref{lem:patterns}]
Let $\Delta\mathcal{P}^t$ denote the set of colored triangles in the $t$-th step, for $t=0, 1,\dots, 2r+3$.
Consider their projections $p_1(\Delta\mathcal{P}^t)$, $p_2(\Delta\mathcal{P}^t)$ and $p_3(\Delta\mathcal{P}^t)$.
Let $p(\Delta\mathcal{P}^t) :=  \bigcup_{k=1}^3 p_k(\Delta\mathcal{P}^t)$.
By \autoref{fig:special_patterns_step},
\begin{align*}
& p_1(\Delta\mathcal{P}^0) = p_2(\Delta\mathcal{P}^0) =  p_3(\Delta\mathcal{P}^0) = \{[\tfrac{i}{q},\tfrac{i+1}{q}] : i= 0, 9r+5 , 18r+10\};\\
& p_1(\Delta\mathcal{P}^1) = p_2(\Delta\mathcal{P}^1) =  \{[\tfrac{i}{q},\tfrac{i+1}{q}] : i=1,2,4,6,9r+2, 9r+4\},\\
& p_3(\Delta\mathcal{P}^1) = \{[\tfrac{i}{q},\tfrac{i+1}{q}] : i=9r+6, 9r+8, 18r+4, 18r+6, 18r+8, 18r+9\};\\\displaybreak[1]
& p_1(\Delta\mathcal{P}^t) = p_2(\Delta\mathcal{P}^t) =  \{[\tfrac{i}{q},\tfrac{i+1}{q}] : i=
6t-9, 6t-7, 6t-5, 6t-4, 6t-2, 6t, \\ &\quad 9r-3t+5, 9r-3t+7, 9r-3t+9\}, \\
& p_3(\Delta\mathcal{P}^t) = \{[\tfrac{i}{q},\tfrac{i+1}{q}] : i=9r+3t+1, 9r+3t+3, 9r+3t+5, 18r-6t+10,  \\ &\quad 18r-6t+12, 18r-6t+14, 18r-6t+15, 18r-6t+17, 18r-6t+19\},\\
& \qquad\qquad\text{for }t=2,3,\dots, r; \\\displaybreak[1]
& p_1(\Delta\mathcal{P}^{r+1}) = p_2(\Delta\mathcal{P}^{r+1}) =  \{[\tfrac{i}{q},\tfrac{i+1}{q}] : i=6r-3, 6r-1, 6r+1, \\&\hspace{5cm}6r+2, 6r+3, 6r+4, 6r+6\},\\
& p_3(\Delta\mathcal{P}^{r+1}) = \{[\tfrac{i}{q},\tfrac{i+1}{q}] :
                                     i=12r+4, 12r+6, 12r+7, \\&\hspace{5cm} 12r+8, 12r+9, 12r+11, 12r+13\}.
\end{align*}
The sets $p(\Delta\mathcal{P}^t)$ for $t=r+2,r+3,\dots,2r+3$ can be obtained through the mapping $x \mapsto 1-x$. 

Note that for each $t=0,1,\dots,r+1$, 
$[\tfrac{i}{q}, \tfrac{i+1}{q}] \in p_1(\Delta\mathcal{P}^t) = p_2(\Delta\mathcal{P}^t)$ if and only if $[\tfrac{1}{2}-\tfrac{i+1}{q}, \tfrac{1}{2}-\tfrac{i}{q}]=[\tfrac{18r+10-i}{q},\tfrac{18r+11-i}{q}] \in p_3(\Delta\mathcal{P}^t)$. Therefore, the set $p(\Delta\mathcal{P}^t)$ is stable under the reflection corresponding to the symmetry condition.

We now show that, for each $t=0,1,\dots, r+1$, the intervals in $p(\Delta\mathcal{P}^t)$ are from the same connected component. 

In Step $0$, consider the $3$ green triangles on the yellow diagonal stripe with $p_3 = [\tfrac{q-2}{2q}, \tfrac{1}{2}]$. Since they have the same $p_3$ projection, their $p_2$ projections which form the set $p(\Delta\mathcal{P}^0)$ are from the same connected component. 

In Step $1$, consider the $6$ green triangles on the yellow diagonal stripe with $p_3=[\tfrac{9r+6}{q}, \tfrac{9r+7}{q}]$. 
Their $p_2$ projections $\mathcal{P} := \{[\tfrac{i}{q},\tfrac{i+1}{q}] : i=1,2,4,9r+2, 9r+4\}$ are from the same connected component, say $\mathcal{C}$. 
Let $F$ denote the green lower triangle whose $p_1(F)=[\tfrac{9r+2}{q},\tfrac{9r+3}{q}]$ and $p_2(F)=[\tfrac{6t}{q},\tfrac{6t+1}{q}]$. 
Then $p_2(F) \in \mathcal{C}$ since $p_1(F) \in \mathcal{P} \subseteq \mathcal{C}$.  $p_1(\Delta\mathcal{P}^{1}) = p_2(\Delta\mathcal{P}^{1}) \subseteq \mathcal{C}$. Using the reflection $x \mapsto (\tfrac{1}{2}-x) \bmod 1$ corresponding to the symmetry condition, we have $p(\Delta\mathcal{P}^1) \subseteq \mathcal{C}$. 

In Step $t$ for $t=2,3,\dots, r$, consider the $8$ orange triangles on the yellow diagonal stripe with $p_3=[\tfrac{9r+3t+3}{q}, \tfrac{9r+3t+4}{q}]$. Since they have the same $p_3$ projection, their $p_2$ projections $\mathcal{P} := \{[\tfrac{i}{q},\tfrac{i+1}{q}] : i= 6t-7,\allowbreak 6t-5,\allowbreak 6t-4,\allowbreak 6t-2,\allowbreak 9r-3t+5,\allowbreak 9r-3t+7, 9r-3t+9\}$ are from the same connected component, say $\mathcal{C}$. Let $F_1$ denote the orange upper triangle whose $p_1(F_1)=[\tfrac{9r-3t+9}{q},\tfrac{9r-3t+10}{q}]$ and $p_2(F_1)=[\tfrac{6t-9}{q},\tfrac{6t-8}{q}]$. Let $F_2$ denote the orange lower triangle whose $p_1(F_2)=[\tfrac{9r-3t+5}{q},\tfrac{9r-3t+6}{q}]$ and $p_2(F_2)=[\tfrac{6t}{q},\tfrac{6t+1}{q}]$. Since $p_1(F_1) \in \mathcal{P} \subseteq \mathcal{C}$, $p_2(F_1) \in \mathcal{C}$. Similarly, since $p_1(F_2) \in \mathcal{P} \subseteq \mathcal{C}$, $p_2(F_2) \in \mathcal{C}$. Thus $p_1(\Delta\mathcal{P}^{t}) = p_2(\Delta\mathcal{P}^{t}) \subseteq \mathcal{C}$. Using the reflection $x \mapsto (\tfrac{1}{2}-x) \bmod 1$ corresponding to the symmetry condition, we have $p(\Delta\mathcal{P}^t) \subseteq \mathcal{C}$. 

In Step $r+1$, consider the $4$ red triangles on the yellow diagonal stripe with $p_3 = [\tfrac{12r+4}{q}, \tfrac{12r+5}{q}]$. 
Their $p_2$ projections $\mathcal{P} := \{[\tfrac{i}{q},\tfrac{i+1}{q}] : i=6r-3, 6r-1, 6r+4, 6r+6\}$  are from the same connected component, say $\mathcal{C}$. 
Let $F_k$ denote the red upper triangle whose $p_1(F_k)=[\tfrac{6r+k}{q},\tfrac{6r+k+1}{q}]$ and $p_2(F_k)=[\tfrac{6t+4}{q},\tfrac{6t+5}{q}]$, for $k=1,2,3$. Since $p_2(F_k) \in \mathcal{P} \subseteq \mathcal{C}$, $p_1(F_k) \in \mathcal{C}$ for each $k=1,2,3$. Thus $ p_1(\Delta\mathcal{P}^{r+1}) = p_2(\Delta\mathcal{P}^{r+1}) \subseteq \mathcal{C}$.
Finally, by the reflection $x \mapsto (\tfrac{1}{2}-x) \bmod 1$ corresponding to the symmetry condition, we conclude that all elements of $p(\Delta\mathcal{P}^{r+1})$ are from the same connected component $\mathcal{C}$.

By construction, $p(\Delta\mathcal{P}^0), p(\Delta\mathcal{P}^1), \dots,
p(\Delta\mathcal{P}^{r+1})$ form a partition of the set
$\{[\tfrac{i}{q},\tfrac{i+1}{q}] : i=0, 1, \dots, 18r+10\}$. Recall that $q =
36r+22$. Then by the invariance under the mapping $x \mapsto 1-x$, we have that
$p(\Delta\mathcal{P}^{r+2}), p(\Delta\mathcal{P}^{r+3}), \allowbreak\dots,\allowbreak
p(\Delta\mathcal{P}^{2r+3})$ form a partition of the set
$\{\,[\tfrac{i}{q},\tfrac{i+1}{q}] : i=18r+11, 1, \dots, 36r+21\,\}$. Therefore,
$p(\Delta\mathcal{P}^i) \cap p(\Delta\mathcal{P}^j) = \emptyset$ for any $ i
\neq j, 0 \leq i, j \leq 2r+3$. Since we have considered all additivities
corresponding to the painting that could give rise to merging of components,
it follows that $p(\Delta\mathcal{P}^t)$ for $t=0, 1, \dots, 2r+3$ are
$2(r+2)$ connected components. Furthermore, the intervals
$[\tfrac{i}{q},\tfrac{i+1}{q}]$ for $i=0, 1, \dots, q-1$ are all directly
covered by the prescribed partial painting.
%
\end{proof}

\begin{proof}[Proof of \autoref{lem:isomorphism}]
On the one hand, the function values $(\pi_0, \pi_1, \dots, \pi_q)$ can be expressed in terms of $(s_0, s_1, \dots, s_{r+1})$, as follows.
\begin{align*}
\pi_{6i} &= 6\textstyle\sum_{j=1}^i{s_j} + s_0-2s_1-s_i+2s_{i+1},&& i = 0, 1, \dots, r; \\
\pi_{6i+1} &= 6\textstyle\sum_{j=1}^i{s_j} + s_0-2s_1+2s_{i+1},&& i = 0, 1, \dots, r; \\\displaybreak[1]
\pi_{6i+2} &= 6\textstyle\sum_{j=1}^i{s_j} + s_0-2s_1+3s_{i+1},&& i = 0, 1, \dots, r; \\
\pi_{6i+3} &= 6\textstyle\sum_{j=1}^i{s_j} + s_0-2s_1+4s_{i+1},&& i = 0, 1, \dots, r; \\\displaybreak[1]
\pi_{6i+4} &= 6\textstyle\sum_{j=1}^i{s_j} + s_0-2s_1+4s_{i+1}+ s_{i+2},&& i = 0, 1, \dots, r-1; \\
\pi_{6i+5} &= 6\textstyle\sum_{j=1}^i{s_j} + s_0-2s_1+5s_{i+1}+ s_{i+2},&& i = 0, 1, \dots, r-1; \\\displaybreak[1]
\pi_{9r+5-3i} &= 9\textstyle\sum_{j=1}^r{s_j} - 3\textstyle\sum_{j=1}^i{s_j} \\&\quad\ + s_0-2s_1+7s_{r+1}-s_{i+1},&& i = 0, 1, \dots, r; \\
\pi_{9r+4-3i} &= 9\textstyle\sum_{j=1}^r{s_j} - 3\textstyle\sum_{j=1}^i{s_j} \\&\quad\ + s_0-2s_1+7s_{r+1}-2s_{i+1},&& i = 0, 1, \dots, r; \\
\pi_{9r+3-3i} &= 9\textstyle\sum_{j=1}^r{s_j} - 3\textstyle\sum_{j=1}^i{s_j}
  \\&\quad\  + s_0-2s_1+7s_{r+1}-2s_{i+1}-s_{i+2},&& i = 0, 1, \dots, r-1.
\end{align*}
(These formulas are obtained by integrating the slope values and are easily verified by induction.)
By the symmetry condition, for $t = 0, 1, \dots, 9r+5$, 
\[
\pi_{9r+6+t} + \pi_{9r+5-t} = 18\textstyle\sum_{j=1}^r{s_j} + 3s_0-6s_1+14s_{r+1} = 1.
\]
By the invariance property, for $t = 0, 1, \dots, 18r+22$, $\pi_t = \pi_{q - t}$.

On the other hand, the slope values $(s_0, s_1, \dots, s_{r+1})$ can be expressed in terms of the function values $\pi_i$: $s_0 = \pi_1-\pi_0$, $s_1=\pi_2-\pi_1$, and $s_t = \pi_{6t-8}-\pi_{6t-9}$ for $t=2, 3,\dots, r+1$.
\end{proof}

\begin{proof}[Proof of \autoref{lem:ordering_of_slopes}]
We show the ordering of $(s_0, s_1, \dots, s_{r+1})$ by considering the subadditivity of $\pi$. Since $\pi_1 + \pi_1 \geq \pi_2$, we have $s_0+s_0 \geq s_0+s_1$, and hence $s_0 \geq s_1$. For $i=0,1,\dots,r-1$, the subadditivity condition $\pi_{6i+3}+\pi_{9r+3-3i} \geq \pi_{9r+6+3i}$ implies that
\begin{multline*}
(6\sum_{j=1}^i{s_j} + s_0-2s_1+4s_{i+1}) + (9\sum_{j=1}^r{s_j} - 3\sum_{j=1}^i{s_j} + s_0-2s_1+7s_{r+1}-2s_{i+1}-s_{i+2}) \\
\geq  (18\sum_{j=1}^r{s_j} + 3s_0-6s_1+14s_{r+1}) - (9\sum_{j=1}^r{s_j} - 3\sum_{j=1}^i{s_j} + s_0-2s_1+7s_{r+1}-s_{i+1}).
\end{multline*}
After simplification, we have $s_{i+1} \geq s_{i+2}$.
\end{proof}

\clearpage

\section*{Acknowledgments}
The authors wish to thank Robert Hildebrand for his helpful comments. The authors are thankful for insightful comments from the anonymous referees that helped to improve the paper.
\providecommand\ISBN{ISBN }
\bibliographystyle{../../amsabbrvurl}
{\small
\bibliography{../../bib/MLFCB_bib}
}

\end{document}